\documentclass{amsart}

\usepackage{hyperref}
\hypersetup{breaklinks=true}

\numberwithin{equation}{section}
\theoremstyle{plain}
\newtheorem{theorem}[equation]{Theorem}
\newtheorem{proposition}[equation]{Proposition}
\newtheorem{lemma}[equation]{Lemma}
\newtheorem{corollary}[equation]{Corollary}

\theoremstyle{remark}

\theoremstyle{definition}
\newtheorem{definition}[equation]{Definition}

\newtheorem{question}[equation]{Question}
\newtheorem*{question*}{Question}

\usepackage{graphicx}

\usepackage{url}

    \newtheoremstyle{TheoremNum}
        {1}{1}              %%% space between body and thm
        {\itshape}                      %%% Thm body font
        {}                              %%% Indent amount (empty = no indent)
        {\bfseries}                     %%% Thm head font
        {.}                             %%% Punctuation after thm head
        { }                             %%% Space after thm head
        {\thmname{#1}\thmnote{ \bfseries #3}}%%% Thm head spec
    \theoremstyle{TheoremNum}
    \newtheorem{theoremN}{Theorem}

\title[Poincar\'e inequalities and differentiability]{Characterizing spaces satisfying Poincar\'e Inequalities and applications to differentiability }
\author{Sylvester Eriksson-Bique}

\newcommand{\QED}{
\begin{flushright}
$\square$
~\\ ~\\
\end{flushright}}

\newcounter{prob}
\newcommand{\N}{\ensuremath{\mathbb{N}}}
\newcommand{\NN}{\ensuremath{\mathcal{N}}}
\newcommand{\R}{\ensuremath{\mathbb{R}}}

\newcommand{\scl}{\textbf{Sc}}
\newcommand{\loc}{\textbf{Loc}}

\newcommand{\gap}{\textbf{gap}}

\newcommand{\LIP}{\ensuremath{\ \mathrm{LIP\ }}}
\newcommand{\diam}{\ensuremath{\ \mathrm{Diam\ }}}
\newcommand{\Lip}{\ensuremath{\mathrm{Lip\ }}}
\newcommand{\lip}{\ensuremath{\mathrm{lip\ }}}
\newcommand{\defeq}{\mathrel{\mathop:}=}

\newcommand{\nd}{\ensuremath{\overline{d}}}
\newcommand{\nmu}{\ensuremath{\overline{\mu}}}
\newcommand{\nB}{\ensuremath{\overline{B}}}
\newcommand{\neps}{\ensuremath{\overline{\epsilon}}}
\newcommand{\nC}{\ensuremath{\overline{C}}}
\newcommand{\nX}{\ensuremath{\overline{X}}}
\newcommand{\nK}{\ensuremath{\overline{K}}}
\newcommand{\nD}{\ensuremath{\overline{D}}}

\newcommand{\nS}{\ensuremath{\overline{S}}}
\newcommand{\nZ}{\ensuremath{\overline{Z}}}
\newcommand{\ngam}{\ensuremath{\overline{\gamma}}}

\newcommand{\Mod}{\ensuremath{\mathrm{Mod}}}

\newcommand{\M}{\ensuremath{\mathcal{M}}}

\newcommand{\len}{\ensuremath{\mathrm{len}}}
\newcommand{\undf}{\ensuremath{\mathrm{Undef}}}

\def\XXint#1#2#3{{\setbox0=\hbox{$#1{#2#3}{\int}$ }
\vcenter{\hbox{$#2#3$ }}\kern-.58\wd0}}

\newcommand{\co}{\mskip0.5mu\colon\thinspace}   % Colon for maps.

\def\vint_#1{\mathchoice%
        {\mathop{\kern 0.2em\vrule width 0.6em height 0.69678ex depth -0.58065ex
                \kern -0.8em \intop}\nolimits_{\kern -0.4em#1}}%
        {\mathop{\kern 0.1em\vrule width 0.5em height 0.69678ex depth -0.60387ex
                \kern -0.6em \intop}\nolimits_{#1}}%
        {\mathop{\kern 0.1em\vrule width 0.5em height 0.69678ex
            depth -0.60387ex
                \kern -0.6em \intop}\nolimits_{#1}}%
        {\mathop{\kern 0.1em\vrule width 0.5em height 0.69678ex depth -0.60387ex
                \kern -0.6em \intop}\nolimits_{#1}}}
\def\vintslides_#1{\mathchoice%
        {\mathop{\kern 0.1em\vrule width 0.5em height 0.697ex depth -0.581ex
                \kern -0.6em \intop}\nolimits_{\kern -0.4em#1}}%
        {\mathop{\kern 0.1em\vrule width 0.3em height 0.697ex depth -0.604ex
                \kern -0.4em \intop}\nolimits_{#1}}%
        {\mathop{\kern 0.1em\vrule width 0.3em height 0.697ex depth -0.604ex
                \kern -0.4em \intop}\nolimits_{#1}}%
        {\mathop{\kern 0.1em\vrule width 0.3em height 0.697ex depth -0.604ex
                \kern -0.4em \intop}\nolimits_{#1}}}

\begin{document}

\begin{abstract}
\noindent  We characterize complete RNP-differentiability spaces as those
spaces which are rectifiable in terms of doubling metric measure spaces satisfying
some local $(1, p)$-Poincar\'e inequalities. This gives a full characterization
of spaces admitting a strong form of a differentiability structure in the
sense of Cheeger, and provides a partial converse to his theorem. The proof
is based on a new ``thickening'' construction, which can be used to enlarge
subsets into spaces admitting Poincar\'e inequalities. We also introduce a new
notion of quantitative connectivity which characterizes spaces satisfying local
Poincar\'e inequalities. This characterization is of independent interest, and
has several applications separate from differentiability spaces. We resolve a
question of Tapio Rajala on the existence of Poincar\'e inequalities for the class
of $MCP(K, n)$-spaces which satisfy a weak Ricci-bound. We show that deforming
a geodesic metric measure space by Muckenhoupt weights preserves
the property of possessing a Poincar\'e inequality. Finally, the new condition
allows us to show that many classes of weak, Orlicz and non-homogeneous
Poincar\'e inequalities ``self-improve'' to classical $(1, q)$-Poincar\'e inequalities for some $q \in [1,\infty)$, which is related to Keith's and Zhong's theorem on self-improvement of Poincar\'e inequalities.

\end{abstract}

\maketitle

\tableofcontents

\iffalse

%Editing remarks

\section{Editing remarks and TODO:}

Structure of paper and edits
\begin{itemize}
\item Section for General Poincare, and section for applications. Self-improvement of Orlicz, Semmes and Non-homogeneous. Applications to Euclidean spaces 
\begin{itemize}
\item Make subsections for each application
\item Define the relevant concepts and exhibit $(C,\delta,\epsilon)$ with a Lemma, and state general result in theorem.
\item For Orlicz state connection to $A_p$ weights, and add citation to Orlicz paper.
\end{itemize}
\item Section for Asymptotic Poincare, and applications.
\begin{itemize}
\item Go over proof and look for places to improve. Try to make re-parametrization clearer. 
\item Add remark on asymptotic non-homogeneous.
\end{itemize}
\end{itemize}

\fi

\section{Introduction}

\subsection{Overview} \label{sec:overview}

This paper focuses on the problem of characterizing geometrically
two analytic conditions on a metric space: the existence of a certain measurable
differentiable structure and the property of possessing a Poincar\'e inequality. Since
they were introduced, many fundamental questions about their relationships and
geometric nature have remained open.

The concepts of differentiability and Poincar\'e inequalities involve somewhat different
terminology, history and techniques. Our theorems relate to both of these and some independently interesting applications. Thus, in order to facilitate readability,
we first give a general overview of our main results, followed by some more
detailed discussion on the individual topics.

Poincar\'e inequalities for metric spaces were introduced by Heinonen and Koskela in \cite{heinonen1998quasiconformal}, and have been central tools in the study of such concepts as Sobolev spaces and quasiconformal maps in metric spaces. Spaces satisfying such inequalities, and which are measure doubling, are called PI-spaces. For precise definitions see Definitions \ref{Pispace} and \ref{def:poincare}. While characterizations of a subclass of PI-spaces appeared in \cite{heinonen1998quasiconformal} and \cite{durand2016geometric}, and a general characterization in \cite{keith2003modulus}, it remained a question to find constructions of these spaces and flexible ways of proving these inequalities in particular contexts. For example, see the detailed discussion in \cite{heinonen2007nonsmooth}.

Measurable differentiable structures, on the other hand, were introduced by Cheeger in \cite{ChDiff99}. His main theorem showed that a PI-space possesses a measurable differentiable structure, which permits differentiation of Lipschitz functions almost everywhere. The spaces satisfying such a theorem, without the PI-space assumption, are called \emph{differentiability spaces} (or Lipschitz differentiability spaces \cite{keith04diff}). See below Definitions \ref{diffspace} and \ref{rnpdiffspace}. 

An early question was if the assumptions of Cheeger were in some sense necessary \cite{heinonen2007nonsmooth}. Because positive measure subsets of PI-spaces may be totally disconnected while they remain differentiability spaces \cite{bate2013differentiability}, strictly speaking, a Poincar\'e inequality cannot be necessary. However, it remained a question if a PI-space structure could be recovered in a weaker form such as by taking tangents or by covering the space in some form. This question was related to better understanding the local geometry of PI-spaces. Various authors, such as Cheeger, Heinonen, Kleiner and Schioppa, posed similar questions. A strong form of this question appeared in \cite{cheeger2015infinitesimal}, where it was asked if every differentiability space is PI-rectifiable. Call a metric space \emph{PI-rectifiable} if it can be covered up to measure zero by a countable number of subsets of PI-spaces. See Definition \ref{def:PIrect}.

Our main result fully resolves the PI-rectifiability question for a subclass of RNP-differentiability spaces. Conversely, together with a result of Cheeger and Kleiner \cite{Ch09diff}, it fully characterizes PI-rectifiable spaces.

\begin{theorem}\label{thm:PIrect} A complete metric measure space $(X,d,\mu)$ is a RNP-Lipschitz differentiability space if and only if it is PI-rectifiable and every porous set has measure zero. \footnote{The assumption of porous sets is somewhat technical and is similar to the discussion in \cite{bate2014characterizations}.}
\end{theorem}

RNP-differentiability is \emph{a priori} a stronger assumption than Lipschitz differentiability. In a RNP-differentiability space one can differentiate all Lipschitz functions with values in RNP-Banach spaces, instead of just ones with finite dimensional targets. These spaces were introduced in the pioneering work of Bate and Li \cite{bate2015geometry}, where they showed that such differentiability spaces satisfied certain asymptotic and non-homogeneous versions of Poincar\'e inequalities. Also, an earlier paper by Cheeger and Kleiner \cite{Ch09diff} showed that every PI-space is a RNP-differentiability space.

Recent examples by Andrea Schioppa in \cite{Schioppa2016} show that, indeed, there exist differentiability spaces which are not RNP-differentiability spaces. Or equivalently (by our result), they are not PI-rectifiable. Schioppa's and our work together demonstrate that differentiability of Lipschitz functions depends on the target, and that a sufficiently strong assumption of differentiability is equivalent to possessing a Poincar\'e inequality in some sense. This work exposes an interesting problem of understanding how this dependence on the target is related to the local geometry of the space. 

The proof of Theorem \ref{thm:PIrect} rests on two contributions that are of independent interest. Our starting point is the work of Bate and Li \cite{bate2015geometry}, where they identified a decomposition of the space into pieces with asymptotic non-homogeneous Poincar\'e inequalities. For us, it is more important that these subsets satisfy a quantitative connectivity condition. In order to prove PI-rectifiability, we need to be able to enlarge, or ``thicken'' these possibly totally disconnected subsets into spaces with better connectivity properties. See Theorem \ref{relative} below for a more detailed discussion. 

Once this enlarged space is constructed, one needs to verify that it satisfies a Poincar\'e inequality. This required a way of identifying a quantitative connectivity condition that is easier to establish for the space, and showing that this connectivity condition is equivalent to a Poincar\'e inequality. In other words, we needed a new and weaker characterization of Poincar´e inequalities. We do this by introducing (in Definition \ref{def:conbility}) a novel condition on a metric measure space that we call $(C,\delta,\epsilon)$--connectivity. An interesting feature of this condition is that it is formally very similar to the definition of Muckenhoupt weights, and thus our methods draw a close formal similarity between the theories of Poincar\'e inequalities and the theory of Muckenhoupt weights (see discussion following Definition \ref{def:conbility}). Some analogies between Poincar\'e inequalities and Muckenhoupt weights have already been observed in relation to self-improvement phenomenons by Keith and Zhong \cite{keith2008poincare}. 

In terms of this condition, we show the following. \\

\begin{theorem}\label{thm:contheorem} A $(D,r_0)$-doubling complete metric measure space $(X,d,\mu)$ admits a local $(1,p)$-Poincar\'e inequality for some $ p \in [1,\infty)$ if and only if it is locally $(C,\delta,\epsilon)$--connected for some $\delta,\epsilon \in (0,1)$. Both directions of the theorem are quantitative in the respective parameters. \\
\end{theorem}

The variable $p$ in the above theorem is the exponent in the Poincar\'e inequality, which measures the quality of the inequality. A larger $p$ means worse connectivity. Notable from our perspective is that the characterization is for a $(1,p)$-Poincar\'e inequality for some $p$, and that our characterization applies for any doubling metric measure space. Previous characterizations either assumed Ahlfors regularity \cite{heinonen1998quasiconformal}, or presumed knowledge of the exponent $p$, such as in \cite{keith2003modulus} and \cite{durand2016geometric}. As demonstrated by examples of Schioppa \cite{schioppa2015poincare}, it is possible for this exponent to be arbitrarily large. Thus, applying characterizations from \cite{keith2003modulus} seem difficult in some cases where the exponent is a priori unknown. Further, our formalism avoids direct usage of modulus estimates, and seems easier to apply in our context.

 The new characterization of Poincar\'e inequalities has several applications of independent interests. The first answers affirmatively a question of Tapio Rajala on the existence of Poincar\'e inequalities on certain metric measure spaces with weak synthetic Ricci curvature bounds. These spaces, called $MCP(K,n)$-spaces, were introduced by Ohta in \cite{ohta}. We show that, at least for a large enough exponent, these spaces satisfy a Poincar\'e inequality. One expects that the exponent $p$ could be chosen to be smaller.

\begin{theorem}\label{thm:mcp} If $(X,d,\mu)$ is a $MCP(K,n)$ space, then it satisfies a  local $(1,p)$-Poincar\'e-inequality for $p>n+1$. Further, if $K\geq 0$, then it satisfies a global $(1,p)$-Poincar\'e inequality.
\end{theorem}

We next consider more general self-improvement phenomena for Poincar\'e inequalities.
In a celebrated paper, Keith and Zhong proved that in a doubling
complete metric measure space a $(1, p)$-Poincar\'e-inequality, with $p\in (1,\infty]$, immediately improves
to a $(1, p-\epsilon)$-Poincar\'e inequality \cite{keith2008poincare} for some $\epsilon>0$ depending on the constants in the doubling and Poincar\'e inequality. We ask if somewhat similarly more general Poincar\'e-type inequalities, say Orlicz-Poincar\'e inequalities, also imply some $(1, p)$-Poincar\'e inequalities. 

By a Poincar\'e type-inequality we will refer to inequalities that control the oscillation of a function by some, possibly non-linear, functional of its gradient. Such inequalities have appeared in the work of Semmes \cite{semmescurv}, Bate and Li \cite{bate2015geometry}, Feng-Yu Wang \cite{fengyu08}, Tuominen \cite{tuom07}, Heikkinen \cite{heikkinenself, heiktuom10} and Jana Bj\"orn in \cite{bjornorlicz}. We will here show that on quasiconvex doubling metric measure spaces many types of ``weak'' Poincar\'e inequalities imply a Poincar\'e inequality. In particular, most of the definitions of Poincar\'e inequalities produce the same category of PI-spaces.

\begin{theorem}\label{thm:nohomo} Suppose that $(X,d,\mu)$ is  a $(D,r_0)$-doubling metric measure space and satisfies a local non-homogeneous $(\Phi, \Psi, C,r_0)$-Poincar\'e-inequali\-ty. Then the space $(X,d,\mu)$ is locally $(\overline{C},\delta,\epsilon)$--connected and moreover admits a $(1,q)$-Poincar\'e inequality for some $q \in [1,\infty)$ and some $\overline{C} \in [1,\infty), \delta,\epsilon \in (0,1)$. All the variables are quantitative in the parameters. 
\end{theorem}

In particular, the non-homogeneous Poincar\'e inequalities considered by Bate and Li in \cite{bate2015geometry} improve to $(1, p)$-Poincar\'e inequalities for \emph{some} $p \in [1,\infty)$. The terminology
used in this theorem is defined in section \ref{subsec:application}.
%One could weaken the assumption of quasiconvexity if we allow any function $f$ instead of Lipschitz functions, or by modifying the inequality to the ones considered in \cite{bate2015geometry}. 
We also obtain the following result strengthening the conclusion of Dejarnette \cite{dejarnette2014self} and Tuominen \cite[Theorem 5.7]{tuomthesis}.

\begin{theorem}
\label{thm:Orlicz} Suppose $(X,d,\mu)$ is a doubling metric measure space satisfying a strong $(1,\Phi)$-Orlicz-Poincar\'e inequality in the sense of \cite{bjornorlicz}, then it satisfies also a $(1,q)$-Poincar\'e inequality for some $q \in [1,\infty)$.
\end{theorem}

\noindent \textbf{Remark:} %By repeating the argument from Cheeger in \cite{ChDiff99} we can conclude that any space with a non-homogeneous Poincar\'e inequality is quasiconvex. 
We remark, that in some senses this result is weaker then \cite{dejarnette2014self} and \cite{keith2008poincare}, because we do not control effectively the range of exponents $q$. Further, we remark that in \cite{bjornorlicz} two different types of Poincar\'e inequalities are considered, a strong and a weak one. Since we do not want to distract the reader with a discussion of Luxembourg norms, we remark that the proof and statement also holds with minor modifications for the weak $(1,\Phi)$-Orlicz-Poincar\'e inequalities defined in \cite{bjornorlicz}. The strong Orlicz-Poincar\'e inequality coincides with the ones considered by Heikkinen and Tuominen in, for example, \cite{heikkinenself, heiktuom10, tuomthesis, tuom07}. The inequalities of Feng-Yu Wang \cite{fengyu08} are of a different nature, since they are in fact stronger than regular $(1,2)$-Poincar\'e inequalities. His inequalities are more related to Orlicz versions of Sobolev-Poincar\'e inequalities.

Most notably, the property of being a PI-space can be recovered even if the function $\Phi$ decays arbitrarily fast at the origin. The exponent $q$ of the obtained
Poincar\'e inequality in Theorem \ref{thm:Orlicz} will grow in such cases, but one cannot fully lose a Poincar\'e inequality by such examples. 

Finally, we will show a theorem concerning $A_\infty$-weights on metric measure spaces. These weights can vanish and blow-up on ``small'' subsets of the space, and thus allow flexibility in obtaining weighted Poincar\'e inequalities. This generalizes some aspects from \cite{Franchi1995} concerning sub-Riemannian metrics and vector fields satisfying the H\"ormander condition. For the definition, see Section \ref{sec:corapp} or Definition \ref{def:muckenhoupt}. 

\begin{theorem}\label{thm:weighted} Let $(X,d,\mu)$ be a geodesic $D$-measure doubling metric measure space. If $(X,d,\mu)$ satisfies a $(1,p)$-Poincar\'e inequality and $\nu \in A_\infty(\mu)$, then there is some $q \in [1,\infty)$ such that $(X,d,\nu)$-satisfies a $(1,q)$-Poincar\'e inequality.
\end{theorem}

In the following subsections, we present the above theorems in more detail and explain some historical connections and proof techniques. A reader who is only interested in the characterization of Poincar\'e inequalities and its applications could only read subsections \ref{subsec:charact},\ref{subsec:application}. On the other hand, a reader mostly interested in the PI-rectifiability of RNP-differentiability spaces could simply read subsection \ref{subsec:diff}.

For general exposition of related ideas and concepts, see the expository articles \cite{heinonen2007nonsmooth, semmes2003introduction, bonk2006quasiconformal, kleiner2006asymptotic}. 

All the results in this paper will be stated for complete separable metric measure spaces $(X, d, \mu)$ equipped with a Radon measure $\mu$ with $0 < \mu(B(x, r)) < \infty$ for all balls $B(x,r) \subset X$. Most of the results could also be modified to apply to non-complete spaces.

\subsection{Characterizing spaces with Poincar\'e inequalities} \label{subsec:charact}

Poincar\'e inequalities and doubling measures are useful tools in analysis of differentiable manifolds and metric spaces alike (see Section \ref{notconvention} for the definitions). Following the convention
established by Cheeger and Kleiner, a space with both a doubling measure and a Poincar\'e inequality is referred to as a \emph{PI-space} or $p$-Poincar\'e space (see Definition
\ref{Pispace}). For a general overview of these spaces and some later developments,
see the beautiful expository article by Heinonen \cite{heinonen2007nonsmooth} and the book \cite{heinonen2015sobolev}.

Due to the many desirable structural properties that PI-spaces have, much effort has been expended to understand what geometric properties guarantee a Poincar\'e inequality in some form. By now, several classes of spaces with Poincar\'e inequalities are known. See for example \cite{Laakso2000,mackay2013modulus, ChCAlmost, jerison1986poincare, rajala2012local, bourdon1999poincare, cheeger2013realization, schioppa2015poincare, schioppa2015poincar,keith2003modulus}. The proofs that these
examples satisfy Poincar\'e inequalities are often challenging, and make extensive
use of the geometry of the underlying space. Our first theorem is thus motivated
by the following question, which has also been posed in \cite{heinonen2007nonsmooth}.

\begin{question} Which geometric properties characterize PI-spaces? Find new and weak conditions on a space that guarantee a Poincar\'e inequality.
\end{question}

As an answer to the above question, we introduce a new connectivity, or avoidance, property for a general metric measure space. This condition involves three parameters $(C,\delta,\epsilon)$ and is referred to as $(C,\delta,\epsilon)$--connectivity. The condition is somewhat technical to state, and we defer to Section \ref{sec:PIproof} and Definition \ref{def:conbility} to state it rigorously. In the vaguest sense, it means that for each pair of points $x, y$ and every
$\epsilon$-low density ``obstacle'' set $E$, there exists a ``curve'' of length $Cd(x, y)$ making
``jumps'' of cumulative size at most $\delta d(x, y)$ while avoiding the obstacle set $E$.

It turns out that Definition \ref{def:conbility} fully characterizes doubling spaces admitting Poincar\'e inequalities. The rigorous statement is contained above in Theorem \ref{thm:contheorem}. Note, the doubling is not needed to show that $(C, \delta, \epsilon)$--connectivity is sufficient for a Poincar\'e inequality, because by Lemma \ref{thm:condoubl} doubling is implied by this condition. Similarly, the theorems below could dispense with the assumption of $D$-doubling. However, we add this to explicate the dependence of the parameters on the implied doubling constant. Note that Lemma \ref{thm:condoubl} gives a fairly bad bound for the doubling, and often the doubling constant may be much smaller.

The rigorous definition of connectivity uses curve fragments instead of curves. The usage of curve fragments may seem counter-intuitive and technical at first. However, as a motivation we note that it becomes easier for a curve to avoid an obstacle $E$ if we permit some jumps, or discontinuities. However, a conceptual difficulty arises as we have to allow infinitely many jumps, which leads to the notion of so called curve fragments. These curve fragments are defined on compact subsets of $\R$ instead of intervals. To simplify, the reader might initially imagine a curve which is continuous except for finitely many jumps. Such curves with jumps in relation to Poincar\'e inequalities have appeared implicitly in prior work, such as in \cite{heinonen1998quasiconformal,bonk2005conformal,ChDiff99,durand2016geometric}.

The reason we also need curve fragments stems from the fact that differentiability spaces may be totally disconnected, and thus may only possess such curve fragments with infinitely many jumps. Further, Bate has shown in \cite{bate12diff} that differentiability
spaces possess a rich supply of curve fragments. Moreover, Bate and Li
also used curve fragments to express connectivity in \cite{bate2015geometry}. We also remark that, if
$X$ is quasiconvex, one can avoid the use of curve fragments in Theorem \ref{thm:contheorem}. This leads to a definition in terms of curves, which is discussed in \cite{sylvester:thesis}.

Our connectivity condition is motivated by a similar condition appearing in the work of Bate and Li \cite[Lemma 3.5]{bate2015geometry}. The main difference is that Bate and Li treat this property only in connection to certain classes of differentiability spaces, while here we isolate it as a property of a general metric measure space (see below for more discussion). Another related condition appears in connection to the $\infty$-Poincar\'e inequality discussed in \cite[Theorem 3.1(f)]{durand2016geometric}, but there one must assume $\mu(E) = 0$.

Theorem \ref{thm:contheorem} fully characterizes PI-spaces. Other similar characterizations of spaces with Poincar\'e inequalities appear in the work of Keith \cite{keith2003modulus} and (for some ranges of exponents and with assumption of homogeneity or Ahlfors regularity) in \cite{heinonen1998quasiconformal, durand2016geometric}, but there the characterization is for a fixed $p$, and in terms of modulus estimates. In some cases, one might be interested in Poincar\'e inequalities without knowing \emph{a priori} the value of $p$ sought, or without efficient control on the homogeneity of the space. Also, in general $p$ might be arbitrarily large \cite{schioppa2015poincare}. In such contexts, our characterization seems easier to apply. In particular, we give applications of our characterization below, which prove Poincar\'e inequalities in new contexts.

Next we outline how the connectivity condition implies a Poincar\'e inequality, the other direction being of a different nature and following from Theorem \ref{thm:nohomo}. The proof is crucially based on a general idea of iteration, where connectivity estimates are iterated via maximal-function type estimates to give stronger connectivity properties and ultimately the Poincar\'e inequality. Below, we use this iteration to first prove quasiconvexity in Proposition \ref{thm:quasiconvex}.

Following quasiconvexity, we use the same iteration scheme to obtain a finer notion of connectivity, which we call fine $\alpha$-connectivity. This is done in order to simplify the proof, and in order to quantify more effectively the exponent $p$ appearing in the Poincar\'e inequality. It should be compared to the relationship between $A_\infty$-weights and $A_p$-weights in classical analysis. The main content is that it allows for controlling the size $\delta$ of the jumps in curve fragments $\gamma$ polynomially in terms of the density parameter for the obstacle $\delta$. More precisely, if an obstacle $E$ is of relative size $\tau \in (0,1)$, then the total size of the gaps can be improved to $C\tau^\alpha$. In particular, instead of the jumps being bounded by $\delta$, we force them to decay with the size of the obstacle in a quantitative way. The precise definition is presented in Definition \ref{def:alphacon}. In terms of this notion, the theorem reads as follows.

\begin{theorem}\label{thm:improvement} Assume $\delta, \epsilon \in (0, 1)$ and $r_0 \in (0,\infty),D \in [1,\infty)$. If $(X, d, \mu)$
is a $(C, \delta,\epsilon,r_0)$--connected $(D,r_0)$-doubling metric measure space, then there exist
$\alpha \in (0,1)$ and $C_1,C_2 \in [1,\infty)$ such that the space is also finely $(\alpha, r_0/(8C_1))$-
connected with parameters $(C_1,C_2)$. We can choose 
$\alpha= \frac{\ln(\delta)}{\ln(\frac{\epsilon}{2M})}$, $C_1=\frac{C}{1-\delta}$ and $C_2=\frac{2M}{\epsilon}$, where $M=2D^{-\log_2(1-\delta)+4}$.
\end{theorem}

\noindent \textbf{Remark:} For example, note that $\R^n$ is finely $\frac{1}{n}$-connected. This means that any  set of relative measure $\epsilon$ in the unit ball can be avoided as long as we permit jumps of total size of the order $\epsilon^\frac{1}{n}$. The tightness of this can be seen by choosing an obstacle set E which is a ball centered at a point. Similarly, the space arising from gluing two copies of $\R^n$ along the origin is finely $\frac{1}{n}$-connected. If the Lebesgue measure $\lambda$ on $\R^n$ is deformed by an $A_p$-Muckenhoupt weight $w$, then the space $(\R^n, | \cdot |, wd\lambda)$ is $\frac{1}{pn}$-connected. This follows from the results in \cite[Chapter V]{stein2016harmonic}.

We further remark, that while the connectivity conditions are very similar to Muckenhoupt conditions, a major difference holds. Both conditions admit a ``selfimprovement'' property of the form in Theorem \ref{thm:improvement}. This can be thought of as the statement that $A_\infty$-weights are $A_p$-weights for some $p$. However, Muckenhoupt weights and Poincar\'e inequalities also possess a different type of self-improvement.
An $A_p$-weight automatically belongs to $A_{p-\epsilon}$ for some $\epsilon > 0$, and a $(1, p)$-Poincar\'e
inequality on a doubling metric measure space improves to a $(1, p -\epsilon)$-Poincar\'e
inequality, when $p>1$ (see \cite{keith2008poincare}). However, a finely $\alpha$-connected space may fail to be finely $\alpha+\epsilon$-connected for any $\epsilon>0$. For example, $\R^n$ is finely $\frac{1}{n}$-connected but not finely $\alpha$-connected for any $\alpha > \frac{1}{n}$. Thus, fine connectivity is not an open condition.

Finally, the last step in showing that Definition \ref{def:conbility} implies a Poincar\'e inequality is to iterate the estimate given by Theorem \ref{thm:improvement} to construct curves along which a function has small integral. This involves a summability condition for a geometric series, which leads to the restriction $p > \frac{1}{\alpha}$. The quantitative version of the aforementioned is the following.

\begin{theorem}\label{thm:alphatheorem} Let $(X,d,\mu)$ be a locally $(D,r_0)$-measure doubling locally finely $(\alpha,r_0)$-con\-nected metric measure space with parameters $(C_1,C_2) \in [1,\infty)^2$. Then, for any $p>\frac{1}{\alpha}$ the space satisfies a local $(1,p)$-Poincar\'e inequality at scale $r_0/(16 C_1)$ with constants $(8C_1, C_{PI})$. In short, the space is a PI-space.

We can set $M=2(C_2D^4)^\frac{1}{p\alpha-1}$, $\delta=C_2\left( \frac{ D^4}{M^p} \right)^{\alpha}$, and $C_{PI}=8D^6 \frac{C_1M}{1-M\delta}$. 
\end{theorem}

The range of $p$'s in this theorem is tight in general. Take the space $Y$ arising as the gluing of two copies of $\R^n$ through their origins. The resulting space $(Y, d, \mu)$, where $d$ is the glued metric and $\mu$ the sum of the measures on each component, is finely $\frac{1}{n}$-connected and satisfies a $(1, p)$-Poincar\'e inequality only for $p > n$. On the other hand, for some particular examples, such as $\R^n$ or Ricci-bounded manifolds \cite{ChCAlmost}, we know that the space satisfies a $(1, 1)$-Poincar\'e inequality, but are only finely $\frac{1}{n}$-connected. Also, if this result is applied to Muckenhoupt weights $w \in A_q(\R^n)$, Theorem \ref{thm:alphatheorem} would give a Poincar\'e inequality for $p > nq$, while a Poincar\'e inequality actually holds for the larger range $p > q$ (see e.g. \cite[Chapter V]{stein2016harmonic} or \cite{keniglocal}). Finally, we remark that this theorem is not an equivalence. A $(1, p)$-PI-space might not be $\alpha$-connected for any $\alpha \geq \frac{1}{p}$.

The iteration scheme used to prove quasiconvexity (Proposition \ref{thm:quasiconvex}), improvement of connectivity (Theorem \ref{thm:improvement}) and Poincar\'e inequalities (Theorem 1.9) all are based on an iterated gap-filling and limiting argument. Recall, that the curve fragments guaranteed by Definitions \ref{def:conbility} and \ref{def:alphacon} can contain gaps, or ``jumps''. However, the core idea is to re-use the connectivity condition at the scale of these gaps and to replace them with finer curve fragments. A similar argument appears in \cite[Lemma 3.6]{bate2015geometry}, and earlier in \cite[Proof of Theorem 5.17]{heinonen1998quasiconformal}. In order to prove quasiconvexity, we can always set the relevant obstacle set to be $E = \emptyset$. However, for the application to Theorems \ref{thm:improvement} and \ref{thm:alphatheorem}, we will need to define obstacles at different scales. To do this successfully, we use a maximal function estimate. Such iterative arguments employing maximal functions resemble the proofs of Theorems V.3.1.4 and V.5.1.3 in \cite{stein2016harmonic}.

Finally, while our connectivity condition and conclusion do not use the modulus estimates of Keith \cite{keith2003modulus} and Heinonen-Koskela \cite{heinonen1998quasiconformal}, it is not surprising that at the end we are able to connect our condition to a certain type of modulus estimate. We refer to Theorem \ref{modulusest} below and the discussion preceding it for the precise statement, as it is not relevant for most of our discussion. For similar modulus bounds in other contexts see \cite{durand2016geometric}. In fact, the techniques of this paper are generally useful for obtaining modulus estimates for certain curve families.

\subsection{Applications of the Characterization} \label{subsec:application}

Our Therem \ref{thm:contheorem} can be used in a variety of contexts to establish Poincar\'e inequalities under \textit{a priori} weaker connectivity properties. 

The first result concerns metric spaces with weak Ricci bounds, where Theorem \ref{thm:mcp} states that any MCP-space admits some Poincar\'e inequality. These spaces originally arose following the work of Cheeger and Colding on Ricci limits \cite{ChCAlmost, ChC97}. Many different definitions appeared such as different definitions for $CD(K,n)$ \cite{sturmCDI, sturmCD, lott2009ricci}, $CD^*(K,n)$ \cite{bacher2010localization} and a strengthening $RCD(K,n)$ \cite{ambrosio2014metric, ambrosio2015riemannian}. Ohta also defined a very weak form of a Ricci bound by the \textit{measure contraction property} and introduced $MCP(K,N)$-spaces \cite{ohta}. Spaces satisfying one of the stronger conditions ($CD(K,N)$, $RCD(K,N)$ or $RCD(K,N)^*$) were all known to admit $(1,1)$-Poincar\'e inequalities \cite{rajala2012interpolated}. It was also known that a non-branching $MCP(K,N)$-space would admit a $(1,1)$-Poincar\'e inequality. This was observed by Renesse \cite{Renesse}, whose proof was essentially a repetition of classical arguments in \cite{ChCAlmost}. However, Rajala had conjectured that this assumption of non-branching was inessential. $MCP(K,n)$ spaces are interesting partly because they are known to include some very non-Euclidean geometries such as the Heisenberg group and certain Carnot-groups \cite{Juillet09,rizzi2016measure}. 

We are left with the following open problem. 

\noindent \textbf{Open Question:} Does every $MCP(K,n)$-space admit a local $(1,1)$-Poincar\'e inequality? \\

We next consider self-improvement phenomena for Poincar\'e inequalities, where we have Theorem \ref{thm:nohomo} stating that a large family of weak Poincar\'e inequalities ``self-improve''.

By a Poincar\'e-type inequality we will refer to inequalities that control the oscillation of a function by some, possibly non-linear, functional of its gradient. First, recall the definition of an upper gradient for metric measure spaces by Heinonen and Koskela \cite{heinonen1998quasiconformal}.

\begin{definition}\label{uppergrad} Let $(X,d,\mu)$  be a metric measure space and $f \co X \to \R$ a Lipschitz function. We call a non-negative Borel-measurable $g$ an upper gradient for $f$ if for every rectifiable curve $\gamma\co [0,L] \to X$ we have
$$|f(\gamma(0)) - f(\gamma(L))| \leq \int_0^L g(\gamma(t)) ds_\gamma.$$
\end{definition}

With this definition, we can define a non-homogeneous Poincar\'e inequality.

\begin{definition}\label{nonhomopoincare} Let $(X,d,\mu)$ be a metric measure space. Let $\Phi, \Psi\co [0, \infty) \to [0, \infty)$ be increasing functions with the following properties.

\begin{itemize}
\item $\lim_{t \to 0}\Phi(t)=\Phi(0)=0$ 
\item $\lim_{t \to 0}\Psi(t)=\Psi(0)=0$
\end{itemize}

We say that $(X,d,\mu)$ satisfies a non-homogeneous $(\Phi,\Psi,C,r_0)$-Poincar\'e inequality if for every $2$-Lipschitz\footnote{The constant $2$ is only used to simplify arguments below. Any fixed bound could be used. Also, by replacing the right hand side with the local Luxembourg norm
\[
\Vert g \Vert_{CB,\Phi,\Psi} = 
\inf \Bigg\{\lambda  > 0 \left| \Psi\left(\vint_{B(x,Cr)} \Phi\left(\frac{g}{\lambda} \right) ~d\mu \right)<1\right. \Bigg\},
\]
we could simply assume that $f$ is Lipschitz. This would lead to the notion of ``weak'' Orlicz-Poincar\'e inequality considered in \cite{bjornorlicz}. The proof for Theorem \ref{thm:nohomo} works just as well for these inequalities, and thus the Poincar\'e inequalities with Luxembourg norms on the right hand side also imply $(1, p)$-Poincar\'e inequalities for $p \in [1,\infty)$. For simplicity, we omit this detail. Using Luxembourg norms may be somewhat more natural due to scaling invariance, which does not hold for Definition \ref{nonhomopoincare}. However, these inequalities have not been studied or used in other contexts
than \cite{bjornorlicz}.} function $f\co X \to \R$, every upper gradient $g \co X \to \R$ and every ball $B(x,r) \subset X$ with $r<r_0$ we have
\begin{equation}\label{eq:nohomo}
 \vint_{B(x,r)} |f-f_{B(x,r)}| ~d\mu \leq r \Psi\left( \vint_{B(x,Cr)} \Phi\circ g~d\mu \right).
\end{equation}
\end{definition}

This class of Poincar\'e inequalities subsumes the ones of Heikkinen and Tuominen in \cite{tuomthesis,tuom07,heiktuom10}, the ones considered by Bj\"orn in \cite{bjornorlicz} (which include the ones of Heikkinen and Tuominen), and the Non-homogeneous Poincar\'e Inequalities (NPI) considered by Bate and Li in \cite{bate2015geometry}. In both of these classes, without any substantial additional assumptions we obtain $(1, p)$-Poincar\'e inequalities for some finite $p \in [1,\infty)$ through Theorem \ref{thm:nohomo}. As already discussed in sub-section \ref{sec:overview}, a corollary gives direct results for Orlicz-Poincar\'e inequalities.

The proof of Theorem \ref{thm:nohomo} is based on defining a function $\rho$, which is, roughly speaking, the smallest size of gaps along a curve fragment connecting a pair of points and avoiding a certain set. If the space doesn't have many curve fragments, then this functional oscillates a lot. However, its gradient is concentrated on the a small subset of the space, which makes the right hand side of inequality \eqref{eq:nohomo} small, and forces the connectivity property to hold with some parameters. 

The final application concerns weighted metric measure spaces. We observed above a formal similarity between Definition \ref{def:conbility} and the definition of Muckenhoupt weights. For our purposes we define Muckenhoupt weights as follows.

\begin{definition} \label{def:muckenhoupt} Let $(X,d,\mu)$ be a $D$-measure doubling metric measure space. We say that a Radon measure $\nu$ is a generalized $A_\infty(\mu)$-measure, or $\nu \in A_{\infty}(\mu)$, if $\nu=w\mu$, and there exist $\epsilon, \delta \in (0,1)$ such that for any $B(x,r)$ and any Borel-set $E \subset B(x,r)$
$$\nu(E) \leq \delta \nu(B(x,r)) \Longrightarrow \mu(E) \leq \epsilon \mu(B(x,r)).$$

If $w\mu \in A_\infty(\mu)$ for some locally integrable $w$, then we call $w$ a Muckenhoupt-weight.
\end{definition}

There are different definitions in the literature. These variants and their equivalence is discussed in \cite{kansanen2011strong}. For geodesic metric measure spaces many of them are equivalent. Also, there is the class of strong $A_\infty$-weights introduced by David and Semmes in \cite{david1989strong} (see also \cite{semmes1996nonexistence,semmes1993bilipschitz}). These weights are somewhat different from Muckenhoupt weights, and usually form a sub-class of them \cite{kansanen2011strong}.

For an unfamiliar reader, we remind that Muckenhoupt weights may vanish and blow-up on subsets of the space. See \cite[Section V]{stein2016harmonic} for examples. Such weights are also somewhat flexible to construct, as alluded to in \cite[Section 3.18]{kenigpipher91}. Interestingly enough, Theorem \ref{thm:weighted} shows that deformations by such weights preserve on geodesic spaces the property of possessing a Poincar\'e inequality.

\subsection{Relationship between differentiability spaces and PI-spaces} \label{subsec:diff}

Cheeger \cite{ChDiff99} defined a metric measure analog of a differentiable structure and proved a powerful generalization of Rademacher's theorem for PI-spaces. The spaces admitting differentiation are here called \emph{differentiability spaces}. See below Definition \ref{diffspace}. Our main result, Theorem \ref{thm:PIrect}, is that a subclass of these spaces, called RNP-differentiability spaces, are PI-rectifiable (see \ref{def:PIrect} for a definition). Thus, within this subclass the conditions of Cheeger are both necessary and sufficient.

 %Since it is one of the main results, we restate it here.
%\begin{theoremN}[\ref{thm:PIrect}] A complete metric measure space $(X,d,\mu)$ is a RNP-Lipschitz differentiability space if and only if it is PI-rectifiable and every porous set has measure zero.
%\end{theoremN}

This result builds on earlier work by Bate and Li \cite{bate2015geometry}, where the authors noticed that RNP-differentiability spaces satisfy certain asymptotic and non-homogeneous Poincar\'e inequalities. Even earlier, a number of similarities between PI-spaces and Poincar\'e inequalities were discovered: asymptotic doubling \cite{bate2013differentiability}, an asymptotic lip-Lip equality almost everywhere \cite{bate12diff, schioppa}, large family of curve fragments representing the measure and ``lines'' in the tangents \cite{cheeger2015infinitesimal}. These works developed ideas and gave strong support for some result like Theorem \ref{thm:PIrect}, which however requires a number of new techniques, such as Theorem \ref{thm:contheorem} and Theorem \ref{relative}. 

The work of Bate and Li initiated the detailed study of RNP-differentiability spaces, which permit differentiation of Lipschitz functions with values in RNP-Banach spaces (\textit{Radon Nikodym Property}). The definition is contained below in Section \ref{rnpdiffspace}. For more information on RNP-Banach spaces see \cite{pisier2016martingales}. The question of differentiability of RNP-Banach space valued Lipschitz functions had already arisen earlier in the work of Cheeger and Kleiner where it was shown that a PI-space admits a Rademacher theorem for such functions \cite{Ch09diff}. 

A corollary of Theorem \ref{thm:PIrect} is the following, which also has an easier, more direct, proof presented in the Appendix \ref{appendix}. 
%Li and Bate managed to prove a partial converse by showing that the tangents of abstract RNP-differentiability spaces admit a weak form of Poincar\'e inequality and that the spaces themselves admit a form of non-homogeneous Poincar\'e inequality infinitesimally. This already showed that differentiability was closely connected to Poincar\'e-type inequalities. However, it left open a question, if differentiability spaces were PI-rectifiable, i.e. if they can be covered by positive measure subsets of PI-spaces. This question is resolved here. 
\begin{theorem}\label{thm:tangents}  Let $(X,d,\mu)$ be a RNP-differentiability space. Then $X$ can be covered, up to measure zero, by countably many positive measure subsets $V_i$, such that each $V_i$ is metric doubling, when equipped with its restricted distance,  and for $\mu$-a.e. $x \in V_i$ each space $M \in T_x(V_i)$ admits a $(1,p)$-Poincar\'e-inequality for some $p \in [1,\infty)$. \footnote{The subsets $V_i$ are equipped with the restricted measure and metric.}
\end{theorem}

Here, $T_x(V_i)$ denotes the set of measured Gromov-Hausdorff tangents at $x$ for the space $(V_i,d|_{V_i}, \mu|_{V_i})$. See \cite{gigligromovhaus, keith2003modulus} for the definition of pointed measured Gromov-Hausdorff convergence, and \cite{ChDiff99, cheeger2015infinitesimal} for the definition of tangent spaces.

We give a few remarks on the proof techniques. Our proof of Theorem \ref{thm:PIrect} starts off with citing the result of Bate and Li that decomposes a RNP-differentiability space $(X,d,\mu)$ into parts with asymptotic and non-homogeneous forms of Poincar\'e inequalities, as well as a uniform and asymptotic form of Definition \ref{def:conbility}. We introduce two new ideas to use these pieces. On the one hand, we observe that our connectivity condition implies a Poincar\'e inequality using Theorem \ref{thm:contheorem}, and this already gives a Poincar\'e inequality for tangents of RNP-differentiability spaces (see the appendix \ref{appendix}). To obtain PI-rectifiability we also need the ability to enlarge the pieces used by Bate and Li to satisfy, intrinsically, the connectivity property in \ref{def:conbility}.

The following theorem is used to construct the enlarged space. The terminology used is defined later in Sections \ref{rectifdiff} and \ref{notconvention}. A crucial observation is that the assumptions involve a relative form of doubling and connectivity, and in the conclusion we construct a space with an intrinsic Poincar\'e inequality and an intrinsic doubling property. Thus, it can be used even for general subsets of PI-spaces to enlarge them to PI-spaces (the enlarged space being different from the original space).

\begin{theorem}\label{relative}   Let $r_0>0$ be arbitrary. Assume $(X,d,\mu)$ is a metric measure space and subsets $K \subset A \subset X$ are given, where $A$ is measurable and $K$ is compact. Assume further that $X$ is $(D,r_0)$-doubling along $A$, $A$ is uniformly $(\frac{1}{2},r_0)$-dense in $X$ along $K$, and $A$ with the restricted measure and distance is locally $(C, 2^{-60}, \epsilon, r_0)$--connected along $K$. There exist constants $\overline{C}, \overline{\epsilon}, \nD >0$, and a complete metric space $\nK$ which is locally $(\overline{D}, 2^{-40}r_0)$-doubling and $(\overline{C}, \frac{1}{2}, \overline{\epsilon}, r_0 2^{-330}/C)$--connected, and an isometry $\iota \co K \to \nK$ which preserves the measure. In particular, the resulting metric measure space $\nK$ is a PI-space.
\end{theorem}

For an intuitive, and slightly imprecise, overview of the proof of this construction one can consider the case of a compact subset $K \subset \R^n$. Here $K$ is well-connected when thought of as a subset of $\R^n$, but may not be intrinsically connected. To satisfy a local Poincar\'e-inequality the space must be locally quasiconvex. In order to make $K$ locally quasiconvex we will glue a metric space $T$ to it. The space $T$ is tree-like, and it's vertices correspond to a discrete approximation of $K$ along with its neighborhood. The vertices exist at different scales, and near-by vertices are attached by edges to each other at comparable scales. By using net-points of $K$, and Whitney centers for a neighborhood of $K$, we prevent adding too many points or edges at any given scale or location.  This construction is analogous to that of a hyperbolic filling \cite[Section 2]{bourdon2003} (see also a more recent presentation in English \cite{bonk2014sobolev}). Very similar ideas also appear in the work of Bonk, Bourdon and Kleiner related to problems of quasiconformal maps \cite{bonk2005conformal,bonk2002quasisymmetric,bonkbourdoncombinatorial}.

\subsection{Structure of paper}  \label{subsec:struct}

We first cover some general terminology and frequently used lemmas in section \ref{notconvention}. In section \ref{sec:PIproof}, we introduce our notion of connectivity and prove basic properties and finally derive Poincar\'e inequalities. In section \ref{sec:corapp}, we apply the results in both new and classical settings. Finally, in section \ref{rectifdiff} we apply the results to the study of RNP-differentiability spaces and introduce the relevant concepts. In the appendix, we include a different proof that tangents of RNP-differentiability spaces are PI-spaces and that our connectivity condition is preserved under measured Gromov-Hausdorff-convergence.

As the paper consists of two related, but somewhat independent ideas, one concerning characterizing differentiability spaces and the other concerning characterization of Poincar\'e inequalities, we have tried to separate these in the structure of the paper. A reader interested only in the classification of spaces with Poincar\'e inequalities could read section \ref{notconvention} followed by the main proofs in section \ref{sec:PIproof}.

A reader interested merely in the application could read section \ref{notconvention} and Definition \ref{def:conbility} after which one can read section \ref{sec:corapp} without needing that much from other parts of the paper. However, the proof of Theorem \ref{thm:nohomo} is closely related to Theorem \ref{thm:contheorem}, and thus is included in section \ref{sec:PIproof}.

Finally, a reader who is simply interested in the PI-rectifiability of RNP-diffe\-rentiability spaces and the involved ``thickening construction'' can move after section \ref{notconvention} directly to section \ref{rectifdiff}, which is mostly self-contained except for references to Theorem \ref{thm:contheorem}. \\ \\

\noindent \textbf{Acknowledgments:} The author is thankful to professor Bruce Kleiner for suggesting the problem on the local geometry of Lipschitz differentiability spaces, for numerous helpful discussions on the topic and several comments that improved the exposition of this paper. Kleiner was instrumental in restructuring the proofs in the third and fifth sections and thus helped greatly simplify the presentation. 

The author also thanks a number of people who have given useful comments in the process of writing this paper, such as Sirkka-Liisa Eriksson, Jana Bj\"orn, Nagesvari
Shanmugalingam, Pekka Koskela, Guy C. David and Ranaan Schul. Some results of the paper were heavily influenced by conversations with Tatiana Toro and Jeff Cheeger. An earlier draft of this paper had a more complicated construction used to resolve Theorem \ref{relative}. This construction is here rephrased in terms of a modified hyperbolic filling which is much clearer than the earlier version. This modification was encouraged by Bruce Kleiner, and suggested to the author by Daniel Meyer. We also thank the anonymous referees for numerous comments and corrections. The research was supported by a NSF graduate student fellowship DGE-1342536  and NSF grant DMS-1405899.

\section{Notational conventions and preliminary results}
\label{notconvention}

We will be studying the geometry of complete and separable metric measure spaces $(X,d,\mu)$. Where not explicitly stated, all the measures considered in this paper will be Radon measures.  An open ball in a metric space $X$ with center $x$ and radius $r$ will be denoted by $B(x,r)$, and for $C>0$ we will denote by $CB(x,r)=B(x,Cr)$. Throughout we will assume that $0 < \mu(B(x, r)) < 1$ for every ball $B(x, r) \subset X$.

\begin{definition}
A metric measure space $(X,d,\mu)$, such that $0<\mu(B(x,r))<\infty$ for all balls $B(x,r) \subset X$, is said to be (locally) $(D,r_0)$-doubling if for all $r \in (0,r_0)$ and any $x \in X$ we have 
\begin{equation}
\frac{\mu(B(x,2r))}{\mu(B(x,r))} \leq D. \label{eq:doubl}
\end{equation}

We say that the space is $D$-doubling if this property holds for every $r_0>0$. Further, we simply call a space doubling if there is a constant $D$ such that it is $D$-doubling, and locally doubling if there are constants $(D,r_0)$ such that it is locally $(D,r_0)$-doubling.
\end{definition}

There is also a metric notion of doubling that does not refer to the measure. Our definition is often referred to as measure doubling. However, throughout this paper, except briefly in relation to Theorem \ref{thm:tangents}, we will only use this stronger version, and thus simply say doubling.

\begin{definition}\label{def:asymptdoubl}
A metric measure space $(X,d,\mu)$ is said to be asymptotically doubling if for almost every $x \in X$ we have 
\begin{equation}\label{eq:asdoubl}
\limsup_{r\to 0}\frac{\mu(B(x,2r))}{\mu(B(x,r))} < \infty.
\end{equation}
\end{definition}

We also define a relative version of doubling for subsets.

\begin{definition} \label{subsetdoubling}
A metric measure space $(X,d,\mu)$ is said to be $(D,r_0)$-doubling along $S \subset X$ if for all $r \in (0,r_0)$ and any $x \in S$ we have 
\begin{equation} \label{eq:subsetdoubl}
\frac{\mu(B(x,2r))}{\mu(B(x,r))} \leq D.
\end{equation}
\end{definition}

\begin{definition} \label{def:porous}
A set $S \subset X$ is called porous if there exist constants $c,r_0 \in (0,\infty)$ such that for every $x \in S$ and every $r \in (0,r_0)$ there exists a $y \in B(x,r)$ such that $B(y,cr) \cap S = \emptyset$. A set $S$ is called $\sigma$-porous if there exist countably many porous sets $S_i$ (with possibly different constants $c_i, r_i$), such that 
$$S = \bigcup_i S_i.$$ 
\end{definition}

\begin{definition} \label{density}
Let $(X,d,\mu)$ be a metric measure space, $\epsilon \in (0,1)$ and $A \subset X$ a positive measure subset. A point $x \in A$ is called an $(\epsilon, r_0)$-density point of $A$ if for any $r \in (0,r_0)$ we have
$$1- \epsilon \leq \frac{\mu(B(x,r) \cap A)}{\mu(B(x,r))} \leq 1.$$ 
We say that $A$ is uniformly $(\epsilon,r_0)$-dense along $S \subset A$ if every point $x \in S$ is an $(\epsilon, r_0)$-density point of $A$.
\end{definition}

A map $f\co X \to Y$ between two metric spaces $(X,d_X)$ and $(Y,d_Y)$ is Lipschitz if there is a constant $L$ such that $d_Y(f(x),f(y)) \leq Ld_X(x,y)$. We will denote by $\LIP f$ the optimal constant in this inequality and call it the global Lipschitz constant of a Lipschitz function. Further, for real-valued $f$ we define the two local Lipschitz constants
\begin{equation}
\Lip f(x) \defeq \limsup_{r \to 0} \sup_{y \in B(x,r)} \frac{|f(x)-f(y)|}{r} \label{Lip}
\end{equation}
and
\begin{equation}
\lip f(x) \defeq \liminf_{r \to 0} \sup_{y \in B(x,r)} \frac{|f(x)-f(y)|}{r}. \label{lip}
\end{equation}

A map $f \co X \to Y$ is called bi-Lipschitz if there is a constant $L$ such that $L^{-1} d_X(x,y) \leq d_Y(f(x),f(y)) \leq Ld_X(x,y)$. The smallest such constant $L$ is called the distortion of $f$.

A continuous embedding $f \co (X,\mu) \to (Y,\nu)$ is said to preserve measure, if $f(X)$ is measurable and $f_*(\mu) = \nu|_{f(X)}$. Further, the push-forward of a measure is defined by $f_*(\mu)(A) \defeq \mu(f^{-1}(A))$.

\begin{definition} \label{def:poincare} Let $p \in [1,\infty)$, $C,C_{PI}>0$ be constants. We say that a complete metric measure space $(X,d,\mu)$ equipped with a Radon measure $\mu$ satisfies a $(1,p)$-Poincar\'e inequality with constants $(C,C_{PI})$ if for every $r>0$, every $x \in X$ and every Lipschitz function $f \co X \to \R^n$.
\begin{equation}
\vint_{B(x,r)} |f-f_{B(x,r)}| ~d\mu  \leq C_{PI}r \left( \vint_{B(x,Cr)}\lip f^p ~d\mu\right)^{\frac{1}{p}}.
\end{equation}

Additionally, we say that the metric measure space $(X,d,\mu)$ satisfies a local ($(1,p)$--)Poincar\'e inequality at scale $r_0>0$ if the aforementioned property holds but only for all $r \in (0,r_0)$.  A space satisfies a local $(1,p)$-Poincar\'e inequality if the aforementioned holds for some $r_0>0$.
\end{definition}

\textbf{Remark:} There are different versions of Poincar\'e inequalities and their equivalence in various contexts has been studied in
 \cite{keith2003modulus} and \cite{hajlasz1995sobolev, Heinonen2000}.  The quantity $\lip f$ on the right-hand side could also be replaces by $\Lip f$ (see equations \eqref{lip} and \eqref{Lip} above), and on complete spaces by an  upper gradient (in any of the senses discussed in  \cite{ChDiff99, heinonen1998quasiconformal, durand2012poincare}). For non-complete spaces, which we do not focus on, the issue is slightly more delicate (see counterexamples in \cite{koskela1999}), but can often be avoided by taking completions. Further, as long as the space is complete, we do not need to constrain the inequality for Lipschitz functions but could also use appropriately defined Sobolev spaces. 
 
Note that on the right-hand side the ball is enlarged by a factor $C\geq 1$. Some authors call inequalities with $C>1$ ``weak'' Poincar\'e inequalities, but we do not distinguish between these different terms. On geodesic metric spaces the inequality can be improved to have $C=1$ \cite{hajlasz1995sobolev}.

We use the following notion of PI-space.

\begin{definition} \label{Pispace}
A complete metric measure space $(X,d,\mu)$ equipped with a Radon measure $\mu$ is called a $(1,p)$-PI space at scale $r_0>0$ with doubling constant $D \geq 1$ and Poincar\'e constants $(C,C_{PI})$ if it is $(D,r_0)$-doubling and satisfies a local $(1,p)$-Poincar\'e inequality at scale $r_0$ with constants $(C,C_{PI})$. Further, a space is called a $(1,p)$-PI space or simply a PI-space if there exist remaining constants so that the space satisfies the aforementioned property.
\end{definition}

\begin{definition} \label{curvfrag} A curve fragment in a metric space $(X,d)$ is a Lipschitz map $\gamma\co K \to X$, where $K \subset \R$ is compact.  %For simplicity, we translate the set so that $\min(K)=0$. 
We say the curve fragment connects points $x$ and $y$ if $\gamma(\min(K))=x, \gamma(\max(K))=y$. Further, define $\undf(\gamma)\defeq [\min(K),\max(K)] \setminus K$. If $K=[\min(K),\max(K)]$ is an interval we simply call $\gamma$ a curve.
\end{definition}

As default, and to simplify notation below, we will assume the curve fragment has been translated so that $\min(K)=0$ unless otherwise stated.

Frequently, we will observe that the open set $\undf(\gamma)$ can be expressed as a countable union of maximal disjoint open intervals $(a_i,b_i)$ as $$\undf(\gamma) = [0,\max(K)] \setminus K = \bigcup_{i} (a_i,b_i).$$ We will employ this notation, up to some necessary subscripts, with only brief comments below. These intervals $(a_i,b_i)$ are also called ``gaps'' or ``jumps''. We also define a measure of the size of these jumps
$$\gap(\gamma) \defeq \sum_{i} d(\gamma(a_i),\gamma(b_i)).$$

The length of a curve fragment is defined as
$$\len(\gamma) \defeq \sup_{x_1 \leq \dots \leq x_n \in K} \sum_{i =1}^n d(\gamma(x_{i+1}),\gamma(x_i)).$$
Since $\gamma$ is assumed to be Lipschitz we have $\len(\gamma) \leq \LIP(\gamma)\max(K).$

Analogous to curves we can define an integral over a curve fragment $\gamma\co K \to X$. Denote $\sigma(t) \defeq \sup_{x_1 \leq \dots \leq x_n \in K \cap [0,t]} \sum_{i =1}^n d(\gamma(x_{i+1}),\gamma(x_i))$. The function $\sigma|_K$ is Lipschitz on $K$. Thus, it is differentiable for almost every density point $t \in K$ and for such $t$ we set $d_\gamma(t) = \sigma'(t)$ and call it the \textit{metric derivative} (see \cite{ambrosio2008gradient} and \cite{bate12diff}). We define an integral of a Borel function $g$ as follows

$$\int_\gamma g ~ds \defeq \int_K g(\gamma(t)) \cdot d_\gamma(t) ~dt,$$
when the right-hand side makes sense. This is true for example if $g$ is bounded from below or above. Naturally, if $K$ has Lebesgue measure zero, then the integral vanishes and is useless. Thus, usually we are primarily interested in curve fragments with domains $K$ of positive measure. 

 We will need to take limits of curve fragments, and in order to do so we present
two auxiliary lemmas on reparametrization and compactness.
 %We will need to do some technical modifications of domains of Lipschitz functions and as such we formulate two lemmas useful for that. First, we give a standard re-parametrization result which generalizes to curve fragment the unit-speed parametrization for curves.

\begin{lemma} \label{adjust}Let $(X,d)$ be a complete metric space, $x,y \in X$ points and $\gamma \co K \to X$ a curve fragment connecting $x$ to $y$. There exists a compact $K' \subset [0, \len(\gamma)]$,
an increasing Lipschitz function $\sigma \co K \to K'$ and a $1$-Lipschitz curve fragment
$\gamma'\co K' \to X$ defined by $\gamma'(\sigma(t)) = \gamma(t)$ for $t \in K$. Moreover, this curve fragment
satisfies the following properties.

\begin{enumerate}
 \item $0, \len(\gamma) \in K'$.
 \item $\len(\gamma)=\len(\gamma')$.
 \item $\gap(\gamma)=\gap(\gamma') = |\undf(\gamma')|$.
 \item  Either $\sigma'(t) = 0$ or $1=d_{\gamma'}(\sigma(t))= d_{\gamma}(t)/\sigma'(t)$, for almost every $t \in K$.
 \item For every non-negative Borel function $g$
 $$\int_{\gamma'} g \, ds = \int_\gamma g \,ds.$$
\end{enumerate}
%There is a compact set $K' \subset [0, \len(\gamma)]$, with $0, \len(\gamma) \in K'$, and a $1$-Lipschitz function $\phi\co K \to K'$ and a $1$-Lipschitz curve fragment $\gamma'$ such that $\gamma' \circ \phi=\gamma$ and $\gamma'$ connects $x$ to $y$. Further, the metric derivative is given, almost everywhere, by $d_{\gamma'}\circ \phi \cdot \phi' = d_{\gamma}$. If $(a_i, b_i) \subset [0,\len(\gamma)] \setminus K'$ is a maximal open interval in $[0,\len(\gamma)] \setminus K'$, we have $d(\gamma'(a_i),\gamma'(b_i))=|a_i-b_i|$, and
%$$\len(\gamma') = \len(\gamma).$$
%In particular the curve integrals over $\gamma$ and $\gamma'$ coincide.
\end{lemma}

%\textbf{Remark:} The equality $d_{\gamma'}\circ \phi \cdot \phi' = d_{\gamma}$ is assumed true also when $\phi'=0=d_\gamma$, regardless of if $d_{\gamma'}\circ \phi$ happens to be defined at such a point.  \\

\noindent \textbf{Proof:} First, by replacing the image space $X$ by $\tilde{X}=\gamma(K)$ we can reduce to the case where $X$ is compact. Next, consider the isometric distance embedding $\iota\co X \to C(X)$, where $C(X)$ is the Banach space of continuous functions equipped with the supremum norm on $X$, and $\iota(x) = d(x, \cdot)$. Then, define the gaps of $\gamma$ as $\bigcup_i (a_i,b_i) = [\min(K), \max(K)] \setminus K$, and define the piecewise linearly extended Lipchitz curve $\overline{\gamma}\co [\min(K), \max(K)] \to C(X)$ by
$\overline{\gamma}(t) = \gamma(t)$ for $t \in K$ and
$$\overline{\gamma}(t) = \frac{b_i-t}{b_i-a_i}\iota(\gamma(a_i)) + \frac{t-a_i}{b_i-a_i}\iota(\gamma(b_i))$$
for $t \in (a_i , b_i)$. This curve can be shown to be Lipschitz, and $\len(\overline{\gamma})=\len(\gamma)$. 

Then, let $\tilde{\sigma} \co [\min(K), \max(K)] \to [0,\len(\gamma)]$ be the length-reparametrization, and $\tilde{\gamma} \co [0,\len(\gamma)] \to X$ the length-reparametrization of $\overline{\gamma}$. Here, we use \cite[Lemma 1.1.4]{ambrosio2008gradient}. Define $\sigma \defeq \tilde{\sigma}|_K$, $K' \defeq \sigma(K)$ and $\gamma' \defeq \iota|_{\iota(X)}^{-1} \circ \tilde{\gamma}|_{K'}$. 

Clearly $0, \len(\gamma) \in K'$, and the properties of the metric derivative and invariance of curve integrals follow from \cite[Lemma 1.1.4]{ambrosio2008gradient}. Finally,
$$\undf(\gamma') = [0, \len(\gamma)] \setminus K' = \bigcup_{i, \tilde{\sigma}(a_i) \neq \tilde{\sigma}(b_i)} (\tilde{\sigma}(a_i), \tilde{\sigma}(b_i)).$$
Also, we have $d_i = \tilde{\sigma}(b_i) - \tilde{\sigma}(a_i) = \len(\overline{\gamma}|_{[a_i ,b_i]}) = d(\gamma(a_i), \gamma(b_i))$, since $\overline{\gamma}$ is linear on the interval $[a_i,b_i]$. Then, $d_i = d(\gamma'(\tilde{\sigma}(a_i)), \gamma'(\tilde{\sigma}(a_i)))$
as well from the definition of $\gamma'$. This gives
$$\gap(\gamma') = \sum_{i}d_i = \sum_{i} d(\gamma(a_i), \gamma(b_i)) = \gap(\gamma).$$
Finally, since $d_i = \tilde{\sigma}(b_i) - \tilde{\sigma}(a_i)$, we also have $\gap(\gamma') = |\undf(\gamma')|$.
%Define $\phi \co K \to \R$ by $\phi(t) = \sup_{x_i \in K \cap [0,t]} \sum_{i =1}^n d(\gamma(x_{i+1}),\gamma(x_i))$, where $x_i \leq x_j$ for $i \leq j$. The function $\phi$ is automatically $1$-Lipschitz and increasing and the result follows by considering $\phi(K)=K'$ and defining $\gamma' = \gamma \circ \phi^{-1}$, where $\phi^{-1}$ is defined as the right-inverse. It is not hard to see that $\gamma'$ is $1$-Lipschitz. Moreover, if $\sigma_\gamma \co K' \to \R$ is the total variation of $\gamma'$, which is defined similarly to $\sigma$, we can show $\sigma_\gamma(\sigma(t))=\sigma(t),$ which gives the property for metric derivatives. The remaining properties are easy to check.
\QED

%Sometimes gaps in a curve will need to be enlarged.

If $\gamma \co K \to X$ is a curve fragment, denote by $\Gamma_\gamma = \{(k, \gamma(k))|k \in K\} \subset \R \times X$ its graph. Then, we say that a sequence of curve fragments $\gamma_i\co K_i \to X$ converges to a curve fragment $\gamma\co K \to X$ if $\lim_{i\to \infty} d_H(\Gamma_{\gamma_i} , \Gamma_\gamma) = 0$. Here $d_H$ is the Hausdorff metric for compact sets defined by
$$ d_H(A, B) \defeq \inf \Bigg\{ \epsilon > 0 \Big| A \subset \bigcup_{b \in B} B(b, \epsilon), B \subset \bigcup_{a \in A} B(a, \epsilon) \Bigg\}.$$

\begin{lemma} \label{adjustment}  Fix $L \in (0, \infty)$. Let X be a complete metric space, $S \subset X$ a compact subset and let $\gamma_i \co K_i \to X$ be a sequence of $L$-Lipschitz curve fragments
with $0 = \min(K_i)$, $\max(K_i) \leq L$, and $\gamma(K_i) \subset S$. There exists a subsequence
converging to a Lipschitz curve fragment $\gamma$ with
\begin{equation}
 \gap(\gamma) \leq \liminf_{i \to \infty} \gap(\gamma_i).
\end{equation}
%Let $(X,d)$ be a complete metric space, $x,y \in X$ points and $\gamma \co K \to X$ a 1-Lipschitz curve fragment.  Assume $[0,\max(K)] \setminus K = \bigcup_i (a_i,b_i)$ where $(a_i, b_i)=I_i$ are the maximal disjoint open intervals, and $C>C_i \geq 1$ are constants. Then, there exists a compact $K'$, a continuous and strictly increasing function $\phi \co [0, \max(K)] \to [0, \max(K')]$ such that $\phi(K)=K'$ and $\gamma'=\gamma \circ \phi^{-1}$ defines a 1-Lipschitz curve fragment on $K'$. Also $[0,\max(K')] \setminus K'=\bigcup_i (\phi(a_i), \phi(b_i))$, $C_i|\phi(a_i)-\phi(b_i)| = |a_i - b_i|$ and
%$$\max(K')\leq \max(K) + (C-1)|\undf(\gamma)|.$$
%The curve integrals over $\gamma$ and $\gamma'$ coincide.
\end{lemma}

\noindent \textbf{Proof:} Assume by passing to a subsequence that 
$$\lim_{i \to \infty} \gap(\gamma_i) = \liminf_{i \to \infty} \gap(\gamma_i).$$
Since the collection of all compact subsets of $[0, L] \times S$ forms a complete metric space under the Hausdorff metric $d_H$, we can choose a subsequence of $\Gamma_{\gamma_i}$ that converges in the Hausdorff metric to a set $\overline{\Gamma} \subset [0, L] \times S$. Since $\gamma_i$ are uniformly $L$-Lipschitz, the set $\overline{\Gamma}$ is a graph of a $L$-Lipschitz function\footnote{If it wasn't a graph of a $L$-Lipschitz function, then there would exist a $u,v \in [0,L]$ such that $(u,a), (v,b) \in \overline{\Gamma}$ for distinct $a,b \in S$ and with $L|u-v| < d(a,b)$. But by the definition of Hausdorff convergence of sets, there would exist $(t_{-,i},a),(t_{+,i},b) \in \Gamma_{\gamma_i}$ such that $\lim_{i \to \infty} (t_{-,i},a) = (u,a)$ and $\lim_{i \to \infty} (t_{+,i},b) = (v,b)$. However, then $\lim_{i \to \infty}|t_{-,i}-t_{+,i}|=|u-v|$, and $\lim_{i \to \infty} d(\gamma_i(t_{-,i}), \gamma_i(t_{+,i})) =d(a,b)>0$, which would lead to a contradiction to the $L$-Lipschitz property of $\gamma_i$ for large enough $i$.} and can be expressed as $\overline{\Gamma}  = \Gamma_\gamma$, where $\gamma \co K \to X$ is a $L$-Lipschitz curve fragment defined on $K \defeq \{k| \exists s \in S, (k, s) \in \overline{\Gamma}\}$. Now, if $\undf(\gamma) = \bigcup_j (a_j , b_j)$, we can show that along the subsequence, and any finite collection of gaps, there is a sequence of maximal open intervals $(a^i_{j}, b^i_{j} ) \subset \undf(\gamma_i)$ which converge to $(a_j , b_j)$ as $i \to \infty$. Thus, the
estimate
$$\gap(\gamma) \leq \lim_{i \to \infty} \gap(\gamma_i),$$
is easy to obtain.
%Define $\phi(t) = \int_0^t (1_K + \sum_i C_i 1_{I_i}) ~dt$. Define $K'=\phi(K)$. Then $\max(K')=\phi(\max(K)) \leq |K| + \sum_{i} C_i |I_i| \leq |K| + C |\undf(\gamma)| \leq |K| + |\undf(\gamma)| + (C-1)\undf(\gamma) \leq \max(K) + (C-1)|\undf(\gamma)|$. The rest of the estimates are easy to see, and the invariance of the curve integral follows since $\phi'=1$ almost everywhere on $K$ and from an argument similar to the previous Lemma.
\QED

For a locally integrable function $f\co X \to \R$ and a ball $B=B(x,r) \subset X$ with $\mu(B(x,r))>0$ we define the average of a function as
$$f_B \defeq \vint_{B} f ~d\mu \defeq \frac{1}{\mu(B(x,r))} \int_{B} f(y) ~d\mu(y).$$

We will frequently use the Hardy-Littlewood maximal function at scale $s>0$:
\begin{equation}\label{maxdefloc}
 \M_s f(x) \defeq \sup_{\overset{x \in B(y,r)}{r<s}} \vint_{B(y,r)} |f|(z) ~d\mu(z),
\end{equation}
and the unrestricted maximal function
\begin{equation}\label{maxdef}
\M f(x) \defeq \sup_{x \in B(y,r)} \vint_{B(y,r)} |f|(z) ~d\mu(z), 
\end{equation}
Refer to \cite{stein2016harmonic} for a standard proof of the following result. The $L^1$ norm of a function is denoted by $||f||_{L^1}$. The possible ambiguity of the space or the measure is not relevant for us, as the measure will be evident from the context.

\begin{theorem}\label{maximal} Let $(X,d,\mu)$ be a $D$-measure doubling metric measure space and $s>0$ and $B(x,r) \subset X$ arbitrary, then for any non-negative integrable function $ f $ and $\lambda>0$ we have
\begin{equation} \label{eq:weakestimate}
 \mu\left(\{ \M_s f > \lambda \}\cap B(x,r)\right) \leq D^3\frac{||f 1_{B(x,r+s)}||_{L^1}}{\lambda} 
\end{equation}
and if the space is $D$-doubling, then \eqref{eq:weakestimate} holds for all $s>0$, and moreover we have
$$\mu\left(\{ \M f > \lambda \}\right) \leq D^3\frac{||f||_{L^1}}{\lambda}.$$
%and thus for any $1<p\leq \infty$
%$$||\M_sf||_{L^p(\mu)} \leq ||\M f||_{L^p}\leq C(D,p) ||f||_{L^p}.$$
%If the space is $D$-doubling, 
\end{theorem}

\begin{definition}
A metric space $(X,d)$ is called $L$-quasiconvex if for every $x,y \in X$ there exists a Lipschitz curve $\gamma\co[0,1] \to X$ such that $\gamma(0)=x, \gamma(1)=y$ and $\len(\gamma) \leq Ld(x,y)$. A metric space is locally $(L,r_0)$-quasiconvex if the same holds for all $x,y \in X$ with $d(x,y) \leq r_0$. A metric space $(X, d)$ is called geodesic if it is $1$-quasiconvex.
\end{definition}

%\begin{definition}
%A metric space $(X,d)$ is called geodesic if for every $x,y \in X$ there exists a Lipschitz curve $\gamma\co [0,1] \to X$ such that $\gamma(0)=x, \gamma(1)=y$ and $\len(\gamma) = d(x,y)$.
%\end{definition}

Any $L$-quasiconvex space is $L$-bi-Lipschitz to a geodesic space by defining a new distance

$$\overline{d}(x,y) = \inf_{\substack{\gamma \co [0,1] \to X \\  \gamma(0)=x, \gamma(1)=y}} \len(\gamma).$$
Finally, we will need an elementary fact concerning averages.
\begin{lemma} \label{average}If $(X,d,\mu)$ is a metric measure space and $f$ is a locally integrable function, and $B=B(x,r)$, then for any $a \in \R$

$$\frac{1}{2r}\vint_B |f-f_B| ~d\mu \leq \frac{1}{r}\vint_B |f-a| ~d\mu.$$
\end{lemma}

\noindent \textbf{Proof:} The result follows by twice applying the triangle inequality.

$$\frac{1}{2r}\vint_B |f-f_B| ~d\mu \leq \frac{1}{2r}\vint_B |f-a|+|a-f_B| ~d\mu \leq \frac{1}{r}\vint_B |f-a| ~d\mu$$
\QED

\section{Proving Poincar\'e inequalities}\label{sec:PIproof}

We define a notion of connectivity in terms of an avoidance property. For the definition of curve fragments see Definition \ref{curvfrag}.

%\footnote{According to Lemma \ref{adjust} we could simply assume some Lipschitz bound, but it is easier to normalize the situation. }

\begin{definition} \label{def:conbility} Let\footnote{The definition is only interesting for $\delta<1$.  However, we allow it to be larger to simplify some arguments below. If $\delta\geq 1$ then the curve fragment could consist of two points, that is $K=\{0,d(x,y)\}$ and $\gamma(0)=x,\gamma(d(x,y))=y$. Indeed, any pair of points in any metric space is $(C,\delta,\epsilon)$--connected, when $\delta,C \geq 1$ and $\epsilon$ is arbitrary. Since the length of a curve fragment connecting $x$ to $y$ is at least $d(x,y)$, the definition also must assume $C \geq 1$ to be meaningful. } $0<\delta$, $0<\epsilon <1$ and $C \geq 1$ be given. If $(x,y) \in X\times X$ is a pair of points with $d(x,y)=r>0$, then we say that the pair $(x,y)$ is $(C,\delta, \epsilon)$--connected if for every Borel set $E$ such that $\mu(E \cap B(x,Cr)) < \epsilon \mu(B(x,Cr))$ there exists a Lipschitz curve fragment $\gamma\co K\to X$  connecting $x$ and $y$, such that the following hold.

\begin{enumerate}
\item $\len(\gamma) \leq Cd(x,y)$. 
\item $\gap(\gamma) < \delta d(x,y)$.
\item $\gamma^{-1}(E) \subset \{\min(K), \max(K)\}$.
\end{enumerate}

We call $(X,d,\mu)$ a $(C,\delta,\epsilon,r_0)$--connected space, if every pair of points $(x,y) \in X$ with $0<d(x,y) \leq r_0$ is $(C,\delta, \epsilon)$-connected.

We say that $X$ is (uniformly) $(C,\delta,\epsilon,r_0)$--connected along $S$, if every $(x,y) \in X$ with $x \in S$ and with $d(x,y)\leq r_0$ is $(C,\delta, \epsilon)$-connected in $X$.

If $(X,d,\mu)$ is $(C,\delta,\epsilon, r_0)$--connected for all $r_0$, we simply say that $(X,d,\mu)$ is $(C,\delta,\epsilon)$--connected.
\end{definition}

\noindent \textbf{General remarks:} The set $E$ above will often be referred to as an ``obstacle''. Since we are working with Radon measures on proper metric spaces, to verify the condition we would only need to consider ``test sets'' $E$ which are either all compact, or just open. For open sets this is trivial, since the measure is outer regular and any Borel obstacle $E$ can be approximated on the outside by an open set $E'$. For compact sets, the argument goes via exhausting an open set by compact sets and obtaining a sequence of curve fragments. 

For certain purposes we could also restrict to using curves $\gamma$, and replace the second and third condition by
$$\int_\gamma 1_E ~ds \leq \delta d(x,y).$$
However, while less intuitive to state with curve fragments, that language is necessary for the application to differentiability spaces in section \ref{rectifdiff}, since the spaces we construct may be \emph{a priori} disconnected. Also, it is often easier to construct curve fragments than curves, since they permits certain jumps. The version of this definition with curves is presented in \cite{sylvester:thesis}. \\

\noindent \textbf{Remark on similarity with Muckenhoupt-weights:} Our motivation for using Definition \ref{def:conbility} and for some details of the proofs below stem from the theory of Muckenhoupt weights. One way of seeing this formal similarity is the following analogy. Consider for simplicity the Lebesgue measure $\lambda$ on $\R^n$. Then, as in Definition \ref{def:muckenhoupt} we say that $\mu \in A_\infty(\lambda)$ if there exist, $\delta,\epsilon \in (0, 1)$ such that for any $B(x, r)$ and any Borel-set $E \subset B(x, r)$
$$\mu(E) \leq \epsilon \mu(B(x, r)) \Longrightarrow \lambda(E) < \delta\lambda(B(x, r)).$$
Notice, that we have switched the roles of $\epsilon$ and $\delta$ here compared to the definition before. The definition above is somewhat similar, if we set for the moment $C =1$. As before, we have that $(\R^n , d, \mu)$ is $(1, \delta, \epsilon)$--connected if for every $E \subset B(x, r)$  Borel and any $y \in \R^n$ with $d(x,y)=r$,
$$\mu(E) \leq \epsilon \mu(B(x, r)) \Longrightarrow  \exists \gamma \co K \to X : \gap(\gamma) < \delta r \ \text{plus two other conditions,}$$
and connecting $x$ to $y$.

In this case the scale-invariant condition $\lambda(E) \leq \delta\lambda(B(x, r))$ is replaced by $\gap(\gamma) < \delta r$, and two other conditions as well as an existential quantifier is added. However, in the simple case of the one dimension space $\R$ one directly sees that a $D$-doubling $\mu$ is an $A_\infty(\lambda)$-measure with constants $(\epsilon , \delta)$ if and only if $(\R, | \cdot |, \mu)$ is $(1, \delta', \epsilon')$--connected. The pairs $(\epsilon , \delta), (\epsilon',\delta') \in (0, 1)^2$ depend quantitatively on each other and doubling. Note that, in this case $\lambda(E \cap (x, y)) \leq \gap(\gamma)$ when $\gamma^{-1}(E) \subset \{0, \max(K)\}$. For further discussion, definitions and a beautiful exposition of the theory of Muckenhoupt weights see \cite[Section V]{stein2016harmonic}. \\

\noindent \textbf{Remark on relation to \cite[Lemma 3.5]{bate2015geometry}:} Bate and Li consider the function $\rho_\epsilon^A$, which could be equivalently defined as
\begin{equation} \label{eq:rhofunc}
 \rho_\epsilon^A(x,y) \defeq \inf_{\gamma \co K \to X} \epsilon(\len(\gamma)-\gap(\gamma)) + \gap(\gamma) = \inf_{\gamma} \epsilon\len(\gamma) + (1-\epsilon)\gap(\gamma), 
\end{equation}
where the infimum is taken over all curve fragments connecting $x$ to $y$ whose image is contained in the set $A \cup \{x,y\}$. See also the proof of Theorem \ref{thm:nohomo} below for a related function. Bate and Li show the following for a RNP-Lipschitz differentiability space $X$. For any $\delta_{LB} \in (0,1)$ and at almost every $x \in X$, there exist $D_{x}  \in (0,1)$ and $R_{x}, \epsilon_{0,x}>0$ such that if $y \in X$ is such that $d(x,y) \in (r/2,r)$ for some $r \in (0,R_x)$, and if $A$ is a Borel set such that
$$\frac{\mu(A \cap B(x,\delta_{LB} r/\epsilon_{0,x}))}{\mu(B(x,\delta_{LB} r/\epsilon_{0,x}))} < D_x,$$
then 
$$\frac{\rho^{A^c}_{\epsilon_{0,x}}(x,y)}{d(x,y)} < \delta_{LB}.$$
Equivalently, there exists a curve fragment $\gamma \co K \to X$ connecting $x$ to $y$ with
$$\epsilon_{0,x} \len(\gamma) + (1-\epsilon_{0,x})\gap(\gamma) < \delta_{LB} d(x,y).$$
For this to be meaningful, we must have $\epsilon_{0,x}<\delta_{LB}$, because otherwise $\epsilon_{0,x} \len(\gamma) > \delta_{LB} d(x,y)$.

Comparing this to Definition \ref{def:conbility}, the set $A$ is equivalent to the obstacle set $E$ and the parameter $D_x$ plays the role of the density parameter $\epsilon$. On the other hand, we could use $C=2\delta_{LB}/\epsilon_{0,x}$ as the length parameter, and the size of the gaps $\delta$ is comparable to $\delta_{LB}$, but more precisely should be set as $\delta=\frac{\delta_{LB}}{1-\delta_{LB}}$. We remark, that the density of $A$ is in the ball of radius $ \delta_{LB} r/\epsilon_{0,x} $, but in Definition \ref{def:conbility} we use $Cd(x,y)$. However, these are comparable since $Cd(x,y)/2 \leq \delta_{LB} r/\epsilon_{0,x} \leq Cd(x,y)$. If now $X$ is asymptotically doubling, we can ensure that the density on the ball of radius $Cd(x,y)$ dominates that of the density on $B(x,\delta_{LB} r/\epsilon_{0,x})$, and we can use a $\epsilon$ in Definition \ref{def:conbility} which is less than $D_x$ by a factor proportional to the doubling constant at $x$. However, by \cite{bate12diff} RNP-differentiability spaces are asymptotically doubling, and this causes no problem. 

Finally, we remark, that our avoidance property of $E$, that is $$\gamma^{-1}(E) \subset \{\min(K), \max(K)\},$$ is slightly different from that used by Bate and Li for $A$, that is ${\rm Im}(\gamma) \subset A^c \cup \{x,y\}$. The first implies the latter with the choice of $E=A$, but for the converse we need to assume that $\gamma$ attains $x$ and $y$ only at $\min(K),\max(K)$, respectively. This can be obtained by restricting $\gamma$ to a subset of its domain starting from the last visit to $x$, and ending at the first following visit to $y$.

With these choices and observations, we can rephrase Bate and Li's result as follows.
\begin{lemma}(Bate-Li \cite[Lemma 3.5]{bate2015geometry})\label{lem:rephrased} If $X$ is a RNP-Lipschitz differentiability space, then for any $\delta \in (0,1)$, there exist $R_x,\epsilon_x, C_x>0$, such that any pair $(x,y)$ with $0<d(x,y)<R_x$ is $(C_x, \delta, \epsilon_x)$--connected.
\end{lemma}

We note that the connectivity condition already implies a doubling bound. 

\begin{lemma}\label{thm:condoubl} Let $(X,d,\mu)$ be a complete metric measure space and $(C,\delta,\epsilon)$-con\-nected for some $C\geq 2$ and $0<\delta<1$, then it is also $D$-measure doubling for some $D>1$. Moreover, if $X$ is locally $(C, \delta, \epsilon, r_0)$--connected then it is also $(D, r_0/2)$-doubling.
\end{lemma}

\noindent \textbf{Proof:} We only prove the non-local version as the local one follows by simply
bounding the scales used. Fix $(x,r) \in X \times (0,\infty)$.  We can assume without loss of generality that $B(x,r/2)\neq B(x,r)$. In other words, we can find some $y \in B(x,r) \setminus B(x,r/2)$. Let $s=d(x,y)$. Define $E \defeq B(x,\delta s)$. We will show that $\mu(B(x,\delta s))\geq \epsilon B(x,r)$. 

If $\mu(B(x,\delta s))< \epsilon B(x,r) < \epsilon B(x,Cs)$, then there is a $1$-Lipschitz curve frag\-ment $\gamma \co K \to X$ connecting $x$ to $y$ with $\gap(\gamma) < \delta s$ and with $\gamma^{-1}(B(x,\delta s)) \subset \{\min(K),\max(K)\}$. For simplicity, translate so that $\min(K)=0$. However, let $\overline{\gamma} \co K \to \R$ be a real-valued curve fragment defined as $\overline{\gamma}(t) \defeq d(\gamma(t), \gamma(0))$. Clearly $\gap(\overline{\gamma}) \leq \gap(\gamma) < \delta s$, but also $\overline{\gamma}(0) = 0, \overline{\gamma}(\max(K)) = s$, and $(0, \delta) \cap {\rm Im}(\overline{\gamma}) = \emptyset$ (since $\gamma^{-1}(B(x,\delta s)) \subset \{0,\max(K)\}$). Thus, $\gap(\overline{\gamma}) \geq \delta s$, which is a contradiction.

Thus, we obtain the lower bound for volume, namely $\mu(B(x, \delta r)) \geq \epsilon \mu(B(x, r))$
for all $r > 0$. If $\delta \leq 1/2$, we have $D$-doubling with $D = 1$. If $\delta > \frac{1}{2}$, let $k$ be the smallest positive integer such that $\delta^k \leq \frac{1}{2}$. A $k$-fold iterated argument then gives us
$$\mu(B(x, r/2)) \geq \mu(B(x, r\delta^k )) \geq \epsilon^k \mu(B(x, r)),$$
which gives doubling property \eqref{eq:doubl} with constant $D =\epsilon^k$.
\QED

%with $|\undf(\gamma)| < \delta s$, and with $\gamma^{-1}(E) \in \{0,\max(K)\}$. Now, $\gamma(0)=x$. We know that $(0,\delta s) \cap K \neq \emptyset$, because otherwise $\undf(\gamma)$ would be too big. Choose a $t \in (0,\delta s)$. But $d(\gamma(t),x) \leq t$ by the Lipschitz bound and $\gamma(t) \not\in E$ by assumption. However, the Lipschitz bound gives also $\gamma(t) \in E$, which is a contradiction.  Thus, we obtain the lower bound for volume.

%Next, define $r_1 = \delta s \leq \delta r_0$. Inductively we can define $r_k \leq \delta^k r_0$ such that $\mu(B(x,r_{k}))\geq \epsilon B(x,r_{k-1})$. This gives, by setting $k=1-\frac{1}{\log_2(\delta)}$, that

%$$\mu(B(x,1/2 r)) \geq \mu(B(x,r_k)) \geq \epsilon^{k} \mu(B(x,r)).$$
%This is the desired estimate \ref{eq:doubl}.
%\QED

This lemma implies that the explicit doubling assumption in the theorems below is somewhat redundant. We leave it in order to explicate the dependence on the constants.

The definition of $(C,\delta,\epsilon)$--connectivity involves almost avoiding sets of a given relative size $\epsilon$. One is led to ask if necessarily better avoidance properties hold for much smaller sets. This type of self-improvement closely resembles the improving properties of reverse H\" older inequalities and $A_\infty$-weights \cite[Chapter V]{stein2016harmonic}. It turns out that, since the connectivity holds at every scale, a maximal function argument can be used to improve the connectivity estimate and allows us to take $\delta$ arbitrarily small. We introduce the notion of finely $\alpha$-connected metric measure space.

\begin{definition}\label{def:alphacon} We call a metric measure space $(X,d,\mu)$ finely $\alpha$-connected with parameters $(C_1,C_2)$ if for any $0<\tau<1$ the space is $(C_1,C_2\tau^\alpha, \tau)$--connected. Further, we say that the space is locally finely $(\alpha,r_0)$-connected with parameters $(C_1, C_2)$ if it is $(C_1,C_2\tau^\alpha, \tau, r_0)$--connected for any $0<\tau<1$. If we do not explicate the dependence on the constants, we simply say $(X,d,\mu)$ is finely $\alpha$-connected if constants $C_1, C_2$ exist with the previous properties. Similarly, we define locally finely $(\alpha,r_0)$-connected spaces.
\end{definition}

Clearly any finely $(\alpha,r_0)$-connected space is also $(C,\delta,\epsilon,r_0)$--connected for some parameters but the converse is not obvious. We next establish this converse for doubling metric measure spaces. This is done by an argument of filling in gaps in an iterative way. This will require maintaining some estimates at smaller scales using a maximal function type argument. Since this argument will also be used later to prove the main Theorem \ref{maintheorem}, we will first explain it in the simple context of proving quasiconvexity. %Quasiconvexity could also be proved by techniques of Sean Li and David Bate \cite[Corollary 6.4]{bate2015geometry}.

\begin{proposition}\label{thm:quasiconvex}
 If $(X,d,\mu)$ is a proper locally $(C, \delta, \epsilon, r_0)$--connected %and locally $(D,r_0)$-doubling 
 metric measure space for some $C \in [1, \infty), \delta, \epsilon \in (0, 1), r_0 \in (0, \infty)$, then it is locally $(L,r_0)$-quasiconvex for $L=\frac{C}{1-\delta}$.
\end{proposition}

\noindent \textbf{Proof:} Take an arbitrary pair of points $(x,y) \in X \times X$ and denote $r\defeq d(x,y) \leq r_0$. Further, abbreviate $D_n\defeq (C+\delta C + \delta^2 C + \cdots + \delta^{n} C)$ and $L \defeq D_\infty \defeq \frac{C}{1-\delta}$.  For each $n \geq 0$ we will recursively define Lipschitz curve fragments $\gamma_n \co K_n \to X$ with the following properties.

%We define for each $n\geq 0$ inductively 1-Lipschitz curve fragments $\gamma_n \co K_n \to X$ with the following properties.

\begin{enumerate}
%\item $\min(K_n)=0$, $\max(K_n) \leq D_nd(x,y)$
\item $\len(\gamma_n) \leq D_n d(x,y) \leq Ld(x,y)$.
\item $\gap(\gamma_n) \leq \delta^n d(x,y)$.
\end{enumerate}

Since $(X, d)$ is locally doubling up to scale $r_0$ (by Lemma 3.2), we also obtain
that closed balls of radius less than $r_0/2$ are compact. After using Lemma \ref{adjust} to parametrize the curve by length, then Lemma \ref{adjustment} and the previous conclusion of compactness can be used to extract a subsequential limit curve  $\gamma$ of $\gamma_n$. Technically, this lemma is applied finitely many times on a subdivision of the curve to pieces of length less than $r_0/2$. Since $\gap(\gamma_n) \leq \delta^n d(x,y)$, the limit curve fragment has no gaps and thus is a curve of the desired length connecting $x$ and $y$. %By the length bound and since $d(x,y) \leq \frac{r_0}{2L}$, we know that all the curve fragments $\gamma_n$ lie within such a ball. Thus, we can apply Lemma \ref{adjustment} to obtain the desired curve by taking a sub-sequential limit. The curve fragments can be made 1-Lipschitz by parametrizing by length using Lemma \ref{adjust}.

%Since $(X,d)$ is proper, this suffices to prove the statement by taking a sub-sequential limit\footnote{A more careful argument could remove taking a sub-sequence, and further would allow for finding curves defined on full measure subsets of an interval without completeness. However, these observations are a bit tedious and distract from the main argument.}.

To initiate the recursion, define $K_0\defeq \{0,r\}, \gamma_0(0)=x,\gamma_0(r)=y$. Next, assume $\gamma_n$ has been constructed, and denote $\undf(\gamma_n)= \bigcup_i (a_i, b_i)$, where $(a_i,b_i)$ are maximal disjoint open intervals. We proceed to construct $\gamma_{n+1}$ with the desired
properties. \\

%\noindent \textbf{Dilation step:} We will seek to patch each gap $(a_i,b_i)$ with a new curve fragment but first need to stretch $\gamma_n$ slightly. By Lemma \ref{adjust} and Lemma \ref{adjustment} we can define a new curve fragment $\gamma'_n \co K_n' \to X$ such that the following holds. Denote the maximal open intervals as $[0,\max(K'_n)]\setminus K_n' = \bigcup_i (a'_i, b'_i)$, where $(a_i', b'_i)$ corresponds to $(a_i,b_i)$ via a re-parametrization. Then, we can ensure that $Cd(\gamma'_n(a'_i),\gamma'_n(b_i))=|b'_i-a'_i|$ and further

%\begin{equation}
%\label{estmax1}
%\max(K_n') \leq \max(K_n) + C\delta^n d(x,y)\leq D_{n+1} d(x,y)
%\end{equation}

\noindent \textbf{Filling in:} Next each of the gaps $(a_i,b_i)$ is filled with a new curve fragment $\gamma_i$. Denote by $d_i=d(\gamma_n(a_i), \gamma_n(b_i))$. By the inductive hypothesis
\begin{equation}
\sum_i d_i = \gap(\gamma_n) \leq \delta^n d(x,y). \label{estgap1}
\end{equation}

By using $(C,\delta,\epsilon)$--connectivity with the obstacle $E=\emptyset$ and the pair of points $\gamma_n(a_i), \gamma_n(b_i)$ we find a $1$-Lipschitz curve fragment $\gamma_n^i \co K_n^i \to X$ connecting $\gamma_n(a_i)$ to $\gamma_n(b_i)$ and $\len(\gamma_n^i) \leq Cd_i$ and $\gap(\gamma_n^i) \leq \delta d_i$. By scaling the domain, we can assume $\min(K_n^i) = a_i$ , $\max(K_n^i) = b_i$. Also, by reparametrizing to have constand speed using Lemma \ref{adjust} we can, without loss of generality, assume that the paths $\gamma_n^i$ for $i \in \N$ are uniformly Lipschitz.

%By possibly extending by a constant path we can assume $\max(K_i)  = C d_i =|b_i'-a_i'|$.

 Next define $K_{n+1} \defeq K_n \bigcup K_n^i$ and $\gamma_{n+1}(t)=\gamma_n(t)$ for $t \in K_n$ and $\gamma_{n+1}(t)=\gamma_n^i(t)$ for $t \in K_n^i$. Then $\gamma_{n+1}$ is easily seen to be a Lipschitz curve fragment connecting $x$ to $y$ and defined on the compact set $K_{n+1}$. Also%By Equation \ref{estmax1} above
\begin{eqnarray*}
\len(\gamma_{n+1}) &\leq& \len(\gamma_n) + \sum_{i} \len(\gamma_n^i)  \\
&\leq& D_n d(x, y) + C\delta^n d(x, y) ≤ D_{n+1} d(x, y).
\end{eqnarray*}

Also, we have
$$\undf(\gamma_{n+1}) = \bigcup_i \undf(\gamma_n^i),$$
and thus% by Equation \ref{estgap1}
\begin{eqnarray*}
\gap(\gamma_{n+1}) &=& \sum_i \gap(\gamma_n^i) \\
&\leq& \sum_i \delta d_i \leq \delta^{n+1}d(x,y).
\end{eqnarray*}
This completes the recursive proof.
\QED

The proof of a part of Theorem \ref{maintheorem} is essentially the same. However, we need to first prove an improved version of connectivity, which we stated in the introduction as Theorem \ref{thm:improvement} and restate here. Both arguments will involve the same scheme as above, but additional care is needed in showing that the curve fragments exist at smaller scales. The proofs are structured to use induction instead of recursion.

\begin{theoremN}[\ref{thm:improvement}] Assume $\delta,\epsilon\in(0,1)$, and $r_0 \in(0,\infty),C,D \in [1,\infty)$. If $(X,d,\mu)$ is a $(C,\delta, \epsilon, r_0)$--connected $(D,r_0)$-doubling metric measure space, then there exist $\alpha \in (0,\infty)$ and $C_1, C_2 \in [1,\infty)$ such that it is also finely $(\alpha,r_0/(8C_1))$-connected with parameters $(C_1,C_2)$. We can choose  $\alpha = \frac{\ln(\delta)}{\ln(\frac{\epsilon}{2M} )}$, $C_1 = \frac{C}{1-\delta}$ and $C_2 = \frac{2M}{\epsilon}$, where $M=2D^{-\log_2(1-\delta)+4}$.
\end{theoremN}

\paragraph{\textbf{Remark}:} We do not need it for the proof, but always $0<\alpha \leq 1$. This could be seen by the argument in Lemma \ref{thm:condoubl}, which can be turned around to give a lower bound for $\delta$ in terms of $\epsilon$ and the optimal doubling constant. \\

\noindent \textbf{Proof:} Define $C_1,C_2,M,\alpha$ as above. Throughout the proof denote by $E$ an arbitrary open set in Definition \ref{def:conbility}. We don't lose any generality by assuming that the obstacles which are tested are open. The proof is an iterative construction, where at every stage gaps in a curve are filled in. %and the curve fragments are slightly stretched. 
This argument could be phrased recursively but we instead phrase it using induction. For the following statement abbreviate $D_n \defeq (C+\delta C + \delta^2 C + \cdots + \delta^{n} C)$. %Fix $M \defeq 2D^{-\log_2(1-\delta)+4}$. 

For each $n \in \N$ we will show the following induction statement $\mathcal{P}_n$: If  $x,y \in X$, $d(x,y)=r<r_0/(8C_1)$ and
$$\mu\left(E \cap B\left(x,C_1r \right)\right) \leq \left(\frac{\epsilon}{2M}\right)^n \mu\left(B\left(x,C_1r\right)\right),$$
then there is a Lipschitz curve fragment $\gamma\co K \to X$ connecting $x$ to $y$ such that the following estimate hold:

\begin{enumerate}
%\item $\min(K)=0, \max(K) \leq D_n d(x,y)$
\item $\len(\gamma) \leq D_n d(x,y) \leq C_1d(x,y)$, 
\item $\gap(\gamma) \leq \delta^n d(x,y)$,
\item $\gamma^{-1}(E) \subset \{0, \max(K)\}$.
\end{enumerate}
This clearly is sufficient to establish the claim. Also, note that it results in the given value $\alpha=\frac{\ln(\delta)}{\ln\left(\frac{\epsilon}{2M} \right)}$. \\

\noindent \textbf{Case $n=1$:} Note, $M>D^{-\log_2(1-\delta)+1}$. Then we get

\begin{eqnarray*}\mu(E \cap B(x,Cr))&\leq& \mu(E \cap B(x,C_1r)) \leq \frac{\epsilon}{2M}  \mu(B(x,C_1r)) \\
&\leq& \frac{\epsilon}{2M} D^{-\log_2(1-\delta)+1} \mu(B(x,Cr)) \leq \frac{\epsilon}{2}\mu(B(x,Cr)).
\end{eqnarray*}
We have assumed that the space is $(C,\delta,\epsilon)$--connected and thus the desired Lipschitz curve fragment exists. \\

\noindent \textbf{Induction step, assume statement for $n$, and prove for $n+1$:}
Take an arbitrary $x,y, E$ with the property $d(x,y)=r<r_0$ and $$\mu\left(E \cap B\left(x,C_1r \right)\right) \leq \left(\frac{\epsilon}{2M}\right)^{n+1} \mu\left(B\left(x,C_1r\right)\right).$$
 Define the set 
$$E_M\defeq \left\{\M(1_{E \cap B(x,C_1r)}) \geq \left(\frac{\epsilon}{2M}\right)^{n}\right\} \cap B(x,C_1r).$$
By standard doubling and maximal function estimates in Theorem \ref{maximal}, if $M \geq D^3$, we have 

\begin{eqnarray*}\mu(E_M \cap B(x,Cr))&\leq& \mu(E_M \cap B(x,C_1r)) \leq \frac{D^3 \mu(E\cap B(x,C_1r))}{\left(\frac{\epsilon}{2M}\right)^{n}} \\
&\leq& \frac{\epsilon}{2M} D^{-\log_2(1-\delta)+4} \mu(B(x,Cr)) \leq \frac{\epsilon}{2}\mu(B(x,Cr)).
\end{eqnarray*}
Now, just as in case $n=1$ using $(C,\delta,\epsilon,r_0)$--connectivity, there is a curve fragment $\gamma'\co K' \to X$ connecting $x$ to $y$ with the properties

\begin{enumerate}
%\item $\min(K')=0, \max(K') \leq C d(x,y)$
\item $\len(\gamma') \leq Cd(x,y)$,
\item $\gap(\gamma') \leq \delta d(x,y)$ and
\item $\gamma'^{-1}(E_M) \subset \{0, \max(K)\}$.
\end{enumerate}

The set $\undf(\gamma')$ is open and as such we can represent $ \undf(\gamma')= \bigcup (a_i, b_i)$ with disjoint intervals such that $a_i,b_i \in K'$. We will seek to patch each gap $(a_i,b_i)$ with a new curve fragment which avoids $E$. %, but first need to stretch $\gamma'$ slightly. This will proceed in the same two steps as before. \\
 Denote $d_i \defeq d(\gamma'(a_i), \gamma'(b_i ))$, and thus
\begin{equation}
 \sum_i d_i = \gap(\gamma') \leq \delta d(x, y). \label{eq:diest}
\end{equation}

% \noindent \textbf{Dilation argument:} Dilate the domain using Lemma \ref{adjust} and Lemma \ref{adjustment}. We will simplify notation by using the same symbols for the dilated curve. Then we can assume $D_n d(\gamma'(a_i),\gamma'(b_i))= |b_i-a_i|$, $\min(K')=0$ and 
% 
% \begin{eqnarray}
% \max(K') \leq Cd(x,y) + D_n\delta d(x,y) \leq D_{n+1} d(x,y) \label{estmax}
% \end{eqnarray}
% 
% Define $d(\gamma'(b_i),\gamma'(a_i))=d_i$. By assumption on $\undf(\gamma')$ (before its dilation),
% 
% \begin{equation}
% \sum_i d_i \leq \delta d(x,y) \label{estgap}
% \end{equation}

\noindent \textbf{Filling in:} For each interval with $a_i \neq 0$ and $b_i \neq \max(K')$  both the endpoints $\gamma'(a_i), \gamma'(b_i)$ do not belong to  $E_M$. If $a_i = 0$, then $\gamma'(b_i) \not\in E_M$ and if $b_i = \max(K')$ then $\gamma'(a_i) \not\in E_M$. This is because $\delta<1$, and occurs only for possible terminal intervals. In any case, for one of the points $\gamma'(a_i), \gamma'(b_i) \not\in E_M$. Say, $\gamma'(a_i) \not\in E_M$. Then, by the maximal function bound for the ball $B(\gamma'(a_i), C_1 d_i)$ we have 

$$\mu(E \cap B(\gamma'(a_i), C_1 d_i ) ) \leq  \left(\frac{\epsilon}{2M} \right)^{n} \mu(B(\gamma'(a_i),C_1 d_i)).$$
Thus, by the induction hypothesis we can define Lipschitz curve fragments $\gamma'_i\co K_i \to X$ connecting $\gamma'(a_i)$ to $\gamma'(b_i)$ such that
\begin{enumerate}
\item $\len(\gamma_i')  \leq D_n d_i$,
\item $$\gap(\gamma_i') \leq \delta^{n} d_i \text{ and}$$
\item $\gamma_i'^{-1}(E) \subset \{\min(K_i),\max(K_i)\}$.
\end{enumerate}

Similarly, if $\gamma'(b_i) \not\in E_M$.  In both cases we can dilate and translate so that
$\min(K_i) = a_i$, $\max(K_i) = b_i$, and by Lemma \ref{adjust} we can parametrize the curve fragments by constant speend and thus assume a uniform Lipschitz bound for $\gamma_i$. %By possibly extending by a constant curve, while keeping the gaps the same size, we can assume $\max(K_i)  = D_n d_i$. 

 Next, define $K \defeq K' \bigcup K_i$ and $\gamma(t)=\gamma'(t)$ for $t \in K'$ and $\gamma(t)=\gamma_i'(t)$ for $t \in K_i$. Then $\gamma$ is easily seen to be a Lipschitz curve fragment connecting $x$ to $y$ and defined on the compact set $K$. Clearly from \eqref{eq:diest}
\begin{eqnarray*}
\len(\gamma) &\leq& \len(\gamma') + \sum_i \len(\gamma_i') \leq Cd(x,y) + \sum_i D_n d_i \\
&\leq& (D+\delta D_n) d(x,y) \leq D_{n+1} d(x,y).
\end{eqnarray*}

Also we have
$$\undf(\gamma) = \bigcup_i \undf(\gamma_i'),$$
and by Equation \eqref{eq:diest}
\begin{eqnarray*}
\gap(\gamma) &=& \sum_i \gap(\gamma_i') \\
&\leq& \sum_i \delta^{n} d_i \leq \delta^{n+1}d(x,y).
\end{eqnarray*}
This completes the induction step, since $\gamma^{-1}(E) \subset \{0, \max(K)\}$ is easy to see since
$E \subset E_M$ (since $E$ is open) and since $\gamma^{-1}(E) \subset \bigcup_i (\gamma'_i)^{-1}(E)$.

\QED

The main application of the previous theorems is the establishment of curves, or curve fragments, that avoid a set on which we have poor control of the oscillation of the function. We will present an example of this type of argument here. \\ 

\noindent \textbf{Example:} Fix $M>1$. Let $(X,d,\mu)$ be finely $(\alpha, r_0)$-connected with parameters $(C_1,C_2)$. If $x,y \in X$, $$\left(\vint_{B(x,C_1 r)} g^p ~d\mu \right)^{\frac{1}{p}} < \epsilon,$$ and $d(x,y)=r<r_0$, then there exists a curve fragment $\gamma\co K \to X$ connecting $x$ to $y$, such that 
$$\int_\gamma g  ~ds \leq C_1M \epsilon d(x,y)$$ 
and $\gap(\gamma) <\frac{C_2}{M^{p\alpha}} d(x,y)$. Observe the emergence of $p$ in the size of the gaps in $K$. For $p$ larger the gaps become smaller. \\

\noindent \textbf{Proof:}  Apply the definition of fine $(\alpha,r_0)$-connectedness to the set $E=\{g>M \epsilon \} \cap B(x,C_1r)$, which has measure at most $\frac{1}{M^p} \mu(B(x,C_1 r))$. Since $g \leq M\epsilon$ on the resulting curve fragment $\gamma$ and $\len(\gamma) \leq C_1 d(x,y)$ we obtain the desired integral estimate. \QED

The basic idea below is to use the same argument as in the previous example,
but to replace the sub-level set of $g$ with a sub-level set for the maximal function.
This allows for iteration and maintaining integral bounds during the filling-in stage
of the argument.

\begin{proposition}\label{prop:smallintegral} Assume $(X,d,\mu)$ is locally finely $(\alpha, r_0)$-connected with parameters $(C_1,C_2)$ and $D$-measure doubling. For any $p>\frac{1}{\alpha}$ there exists a constant $C_3$ with the following property. If $x,y \in X$ are points with $d(x,y)=r<r_0/(8C_1)$, and $g$ is a non-negative Borel function with $$\left(\vint_{B(x,2C_1 r)} g^p ~d\mu \right)^{\frac{1}{p}} \leq 1,$$ then there exists a Lipschitz curve  $\gamma\co [0,L] \to X$ connecting $x$ to $y$, such that 
$$\int_\gamma g  ~ds \leq C_3  d(x,y),$$ 
and $\len(\gamma) \leq \frac{C_1}{1-\delta} d(x,y)$. We can set $C_3=\frac{C_1M}{1-M\delta}$, where $M=2(C_2D^4)^\frac{1}{p\alpha-1}$ and $\delta=C_2\left( \frac{ D^4}{M^p} \right)^{\alpha}$.
\end{proposition}

\noindent \textbf{Remark:} By possibly scaling $g$ we can use the statement with a non-unit $L^p$-bounds for $g$. \\ \\

\noindent \textbf{Proof:} Since our spaces are locally compact, one can approximate $g$ from above by a lower semi-continuous function. The assumption of lower semi-continuity is assumed in order to have $g(x) \leq (\M_s h^p(x))^{\frac{1}{p}}$ for every $s>0$. %, and then from below by an increasing sequence of continuous functions. 
This standard argument (see e.g. \cite[Chapter 1]{fuglede1957extremal}, \cite[Proposition 2.27]{heinonen1998quasiconformal}) allows for assuming that $g$ is lower semi-continuous. Fix $M=2(C_2D^4)^\frac{1}{p\alpha-1}$ and choose a $\delta=C_2\left( \frac{ D^4}{M^p} \right)^{\alpha}$. We have $\frac{1}{M}>\delta>0$, and $M\geq 2$. 

We will construct the curve by an iterative procedure depending on $n$. The proof is structured by an inductive statement. %The procedure is the same as in the previous two proofs: a dilation argument followed by gap-filling. 
The main process of filling is very similar to the prior proofs. However, a new issue arises as we need to ensure that the integral estimate holds at the smaller scale. Abbreviate
$$D_n \defeq (C_1M + C_1M (M\delta) + \cdots + C_1 M(M\delta)^{n-1}),$$
$$L_n \defeq (C_1 + C_1 (\delta) + \cdots + C_1 (\delta)^{n-1})$$
and
$$D_\infty \defeq \lim_{n \to \infty } D_n, L_\infty \defeq \lim_{n \to \infty } C_n.$$

\noindent \textbf{Induction statement:} $\mathcal{P}_n:$ Assume $v,w \in X$ are arbitrary with $d(v,w)=r<r_0/(8C_1)$ and $h$ is any lower semi-continuous non-negative function on $X$ with $$\left(\vint_{B(v,2C_1 r)} h^p ~d\mu \right)^{\frac{1}{p}}\leq 1.$$
Then, there is a Lipschitz curve fragment $\gamma \co K \to X$ connecting $v$ to $w$ such that %with $\min(K)=0, \max(K)\leq C_nd(x,y)$,
\begin{equation}
\int_\gamma h  ~ds \leq  D_n d(v,w), \label{intest}
\end{equation}
\begin{equation}
\gap(\gamma) \leq \delta^n d(v,w), \label{gapest}
\end{equation}
\begin{equation}
\len(\gamma) < L_n d(v,w) \label{lengthest}
\end{equation}
and for any maximal open set $(a_i, b_i) \subset \undf(\gamma)$ we have $d(\gamma(a_i), \gamma(b_i))=d_i$ and for at least one $z_i\in \{a_i, b_i \}$
\begin{equation}
\left(\vint_{B(\gamma(z_i),2C_1 d_i)} h^p ~d\mu \right)^{\frac{1}{p}} \leq M^{n}. \label{maxest}
\end{equation}

Once we have shown $\mathcal{P}_n$ for every $n$ we can apply it with $(x,y)=(v,w)$ and $h=g$ to obtain a sequence of curve fragments $\gamma_n$ with $\int_{\gamma_n} g ~ds \leq D_n d(x,y)$, and by reparametrizing to have unit speed by Lemma \ref{adjust} and taking a limit using Lemma \ref{adjustment}  obtain a curve fragment $\gamma$ with $\len(\gamma) \leq C_\infty d(x,y)$ (from \eqref{lengthest}) and $\int_{\gamma} g ~ds \leq D_\infty d(x,y)$ (from \eqref{intest}). Finally, from \eqref{gapest} we can conclude that $\gap(\gamma)=0$, and thus $\gamma$ is in fact a curve. Also,

$$\int_\gamma g ~ds \leq D_\infty d(x, y).$$
follows from \ref{intest} and since $\gamma_n \to \gamma$ and $\gap(\gamma_n) \to 0$. This can be seen, for example, by using the lower semi-continuity of integrals for curves  and the lower semi-continuity of $g$ (see similar arguments
in \cite[Proposition 4]{keith2003modulus}) and then applying it to the extensions of the curve fragments to curves in a Banach space. The fact that $\gap(\gamma_n) \to 0$ implies that the extra contributions from these extensions tend to zero. This argument additionally would need to consider expressing $g$ as a supremum of a family of Lipschitz functions, and using MacShane extensions of such function to the ambient Banach space.

%In the rest of the proof we show by induction $\mathcal{P}_n$ for every natural number $n$. \\ \\
The rest of the proof consists of showing the statement $\mathcal{P}_n$ by induction in $n$.\\ \\

\noindent \textbf{Base case $n=1$:} Denote $r=d(v,w)$. Define $E_M \defeq \left\{ \M_{2C_1 \delta r} h^p > M^p \right\}$. Then $$\mu(E_M \cap B(v, C_1r)) \leq \frac{D^3}{M^p} \mu(B(v, 2C_1r)) \leq \frac{D^4}{M^p} \mu(B(v, C_1r))$$ by Theorem \ref{maximal}. Thus, by fine connectivity, there is a Lipschitz curve fragment connecting $v$ to $w$ with 
$$\len(\gamma) \leq C_1 d(v,w),$$
$$\gap(\gamma) < C_2\frac{\left(D^4\right)^{\alpha}}{M^{\alpha p}} d(v,w) \leq \delta d(v,w),$$
and
$$\gamma(K \setminus \{0, \max(K) \}) \subset \{ \M_{2C_1 \delta r} h^p \leq M^p \}.$$
This immediately gives estimate \eqref{lengthest}. Next, we prove Estimates \eqref{gapest} and \eqref{intest}.

The estimate \eqref{gapest} follows from the choice $\delta = C_2\left( \frac{ D^4}{M^p} \right)^{\alpha}$. By assumption, $\gamma(K \setminus \{0, \max(K) \}) \subset \{ \M_{2C_1r} h^p \leq M^p  \}$ and $h$ is lower semi-continuous. Thus, it follows that $h \circ \gamma \leq M$ on $K \setminus \{0,\max(K) \}$. In particular, we obtain \eqref{intest} by the simple upper bound

$$\int_\gamma h ~ds \leq \len(\gamma) M \leq C_1 M d(v,w).$$

Finally, we show Estimate \eqref{maxest}. Let $(a_i, b_i) \subset [0,\max(K)] \setminus K$ be an arbitrary maximal open interval. Define $d_i \defeq d(\gamma(a_i),\gamma(b_i))$. %By using Lemma \ref{adjust} we can assume $d(\gamma(a_i), \gamma(b_i))=|b_i-a_i|=d_i$. 
Also, either one of $a_i,b_i$ is not in $\gamma^{-1}(E_M)$, since $d_i \leq \delta r < d(v,w)$. Denote it by $z_i$. The maximal function estimate $\M_{2C_1 \delta r}(\gamma(z_i)) h^p \leq M^p$  combined with $d_i \leq \delta r$ gives

$$\left(\vint_{B(\gamma(z_i),2C_1 d_i)} h^p ~d\mu \right)^{\frac{1}{p}} \leq M^{p}.$$

\noindent \textbf{Assume $\mathcal{P}_n$ and show $\mathcal{P}_{n+1}:$} By the case $\mathcal{P}_n$ there exists a Lipschitz curve fragment $\gamma' \co K \to X$ connecting $v$ to $w$ such that
$$\int_{\gamma'} h  ~ds \leq D_n d(v,w),$$
$$\len(\gamma') \leq C_n d(v,w),$$
$$\gap(\gamma') \leq \delta^n d(v,w),$$
and for any maximal open set $(a_i, b_i) \subset \undf(\gamma)$ with $d_i \defeq d(\gamma'(a_i), \gamma'(b_i))=|b_i-a_i|$ and for at least for one $z_i \in \{a_i, b_i\}$

$$\left(\vint_{B(\gamma'(z_i),2C_1 d_i)} h^p ~d\mu \right)^{\frac{1}{p}} \leq M^{n}.$$
Clearly, also
\begin{equation}
 \sum_{i} d_i = \gap(\gamma') \leq \delta^n d(v,w). \label{eq:gapsest}
\end{equation}

\noindent \textbf{Filling in gaps:} Now, similar to the argument in Lemma \ref{thm:improvement} we will fill in the gaps of $\gamma'$. The set $\undf(\gamma')$ is open, and we can represent $ \undf(\gamma')= \bigcup (a_i, b_i)$ with disjoint intervals such that $a_i,b_i \in K'$. %Define $|b_i-a_i|=d_i$. 
By \eqref{maxest} for one $z_i \in \{ a_i, b_i \}$ we have 
$$\left(\vint_{B(\gamma'(z_i),2C_1 d_i)} h^p ~d\mu \right)^{\frac{1}{p}} \leq M^n.$$

By the base case $n=1$ (applied to the re-scaled $h/M^{n}$ and $v=\gamma'(a_i),w=\gamma'(b_i)$) we can define Lipschitz curve fragments $\gamma'_i\co K_i \to X$ connecting $\gamma'(a_i)$ to $\gamma'(b_i)$ such that the following hold
\begin{itemize}
\item \begin{equation}\label{eq:lengami}
        \len(\gamma_i') \leq C_1 d_i.
      \end{equation}
\item \begin{equation}\label{eq:intithest}
        \int_{\gamma_i'} h ~ds  \leq C_1 M^{n+1} d_i.
      \end{equation}
\item \begin{equation}\label{eq:gapgami}
         \gap(\gamma_i') \leq \delta d_i.
      \end{equation}
\item For any maximal open sets $(a^i_j, b^i_j) \subset \undf(\gamma_i')$ we have for $d^i_j \defeq d(\gamma_i'(a^i_j), \gamma_i'(b^i_j))$ and at least for one $z^i_j \in \{a^i_j, b^i_j\}$

$$\M_{2C_1d^i_j} h^p(\gamma(z^i_j)) \leq M^{(n+1)p} .$$
\end{itemize}

\noindent \textbf{Filling in gaps:} First, parametrize $\gamma_i'$ by unit speed using Lemma \eqref{adjust} and translate and scale the domain so that $\min(K_i) = a_i, \max(K_i) = b_i$. This can be done so that $\gamma_i'$ are uniformly Lipschitz. %It follows that $\max(K') \leq C_n d(v,w) + C_1 \delta^n d(v,w)$, which will give the desired length estimate. 
Next define $K \defeq K' \bigcup K_i$, and define a curve fragment $\gamma \co K \to X$ by $\gamma(t)=\gamma'(t)$ for $t \in K$ and $\gamma(t)=\gamma_i'(t)$ for $t \in K_i$. Then $\gamma$ is easily seen to be a Lipschitz curve fragment connecting $v$ to $w$. The domain $K$ is also clearly compact.% We now verify the inequality  %Note that by continuity of $h$ combined with $\gamma_i'(K_i \setminus \{0,\max(K_i)\}) \subset \{ \M_{2C_1d^i_j} h^p \leq M^{(n+1)p}\}$, we get $h \circ \gamma \leq M^{n+1}$ almost everywhere. This gives the following simple bound for the curve integral

%$$\int_{\gamma_i'} h ~ds \leq \len(\gamma_i') M^{n+1} \leq C_1 M^{n+1} d_i.$$

Now, %adding up the individual contributions and 
using \eqref{eq:intithest} and \eqref{eq:gapsest} we obtain \eqref{intest} for $n+1$.
\begin{eqnarray*}
\int_{\gamma} h ~ds &=& \int_{\gamma'} h ~ds + \sum_i \int_{\gamma_i'} h ~ds \\
 &\leq& D_{n} d(v,w) +  C_1\delta^{n} M^{n+1} d(v,w) \\
&\leq& (C_1M + C_1 M(M\delta) + \cdots + C_1 M (M\delta)^{n}) d(v,w). \\
&=& D_{n+1} d(v,w)
\end{eqnarray*}
Also by combining \eqref{eq:lengami} and \eqref{eq:gapsest}, we get \eqref{lengthest} for $n+1$.
$$\len(\gamma) \leq \len(\gamma') + \sum_i \len(\gamma_i') \leq L_nd(v,w) + C_1 \delta^n d(v,w) \leq L_{n+1} d(v,w).$$

Now, prove the estimate for gaps \eqref{gapest} for $n+1$. Note that
$$\undf(\gamma) = \bigcup_i \undf(\gamma_i'),$$
and as such we obtain the desired estimate
\begin{eqnarray*}
\gap(\gamma) &=& \sum_i \gap(\gamma_i') \\
&\overset{\eqref{eq:gapgami}}{\leq}& \sum_i \delta d_i \overset{\eqref{eq:gapsest}}{\leq} \delta^{n+1}d(v,w).
\end{eqnarray*}

Finally, check the Estimate \eqref{maxest} concerning the integral average. Any maximal open interval $I \subset \undf(\gamma)$ is also a maximal undefined interval for $\gamma_i'$. I.e. we can express any such interval as $I=(a^i_j, b^i_j)$ for some maximal open interval $(a^i_j, b^i_j) \subset \undf(\gamma'_i)$. The desired Estimate \eqref{maxest} is equivalent to the desired estimate for $\gamma_i'$.

\QED

The following lemma can be proven in a similar way as \cite[Proposition 1.11]{ChDiff99}.
%follows by estimating the integral over the curve using similar arguments to \cite{bate2015geometry}.

\begin{lemma} \label{lipest} Let $(X,d,\mu)$ be a metric measure space. If $f$ is a Lipschitz function, and $\gamma\co K \to X$  is a $1$-Lipschitz curve fragment, then for every $a,b \in K$
$$|f(\gamma(a))-f(\gamma(b))| \leq \LIP f |\undf(\gamma) \cap [a,b]| + \int_\gamma \Lip f ~ds.$$
\end{lemma}

Finally, we prove the following quantitative version of the sufficiency in Theorem \ref{thm:contheorem}. This also gives a proof of Theorem \ref{thm:alphatheorem}. The necessity is handled in the proof of Theorem \ref{thm:nohomo} below.

\begin{theorem}\label{maintheorem} Let $(X,d,\mu)$ be a $(D,r_0)$-measure doubling locally $(C, \delta, \epsilon,r_0)$--connected metric measure space. There exist $(C_1, C_2)$ such that the space is finely $(\alpha,r_0/(8C_1))$-connected with constants $(C_1,C_2)$ and $\alpha=\log_{\frac{\epsilon}{2D}}(\delta)$. Moreover, there is a constant $C_{PI}$ such that whenever $p>\frac{1}{\alpha}$ and $0<r<\frac{r_0}{16C_1}$, we have for any Lipschitz function $f$ and any $x \in X$
$$\vint_{B(x,r)} |f -f_{B(x,r)}| ~d\mu  \leq C_{PI} r \left( \vint_{B(x,8C_1 r)} \Lip f^{p} ~d\mu \right)^{\frac{1}{p}}.$$
We can set $C_1 = \frac{C}{1-\delta}$, $C_2 = \frac{4D^{-\log_2(1-\delta)+4}}{\epsilon}$ and $C_{PI}=8D^6\frac{C_1 M}{1-M\eta}$, where $\eta = C_2\left(\frac{D^4}{M^p}\right)^\alpha$ and $M = 2(C_2 D^4)^{p\alpha-1}$. 
% is given by Theorem \ref{thm:smallintegral}. 
In particular, the space is a PI-space.
\end{theorem}

\noindent \textbf{Proof:} Denote $\diam(B(x,r))=s \leq 2r$. Without loss of generality assume that $B(x,r) \neq \{x\}$, as otherwise the left hand side of the inequality would vanish and the claim would be trivial. Also note that $B(x,4C_1 r) \subset B(x,4C_1 s) \subset B(x,8C_1 r)$. This follows by an argument by contradiction. Assume that there exists a $w \in B(x,4C_1 r) \setminus B(x,4C_1 s)$. Then, by Proposition \ref{thm:quasiconvex}, there would be a quasigeodesic joining $w$ to $x$. However, then $s=\diam(B(x,r)) \geq \min\{4C_1 s,r\}$, from which we would get $s \geq r$ giving the desired inclusion. The inclusion guarantees, by doubling, that
\[\frac{\mu(B(x,8C_1r))}{\mu(B(x,4C_1s))} \leq D.\]

Pick $z \in B(x,r)$ with $d(z,x) \geq s/3$. Denote,
$$A=\left( \vint_{B(x,4C_1 s)} \Lip f^{p} ~d\mu \right)^{\frac{1}{p}}.$$
First use Theorem \ref{thm:improvement} to conclude that the space is locally finely $(\alpha,r_0/(8C_1))$-connected with constants $(C_1,C_2)$, where $C_1,C_2$ are as given. Then, apply Proposition \ref{prop:smallintegral} for $h=\Lip f/A$ to give a curve $\gamma$ connecting $x$ to $z$ with
$$\int_\gamma \Lip f ~ds \leq C_3 \left( \vint_{B(x,2C_1 d(x,z))} \Lip f^{p} ~d\mu \right)^{\frac{1}{p}} d(x,z) \leq D^4 C_3 A d(x,z),$$
where $C_3=\frac{C_1M}{1-M\eta}$.
Thus, by using Lemma \ref{lipest}
$$|f(z)-f(x)| \leq D^4 C_3 A s.$$

For the last inequality we used doubling. Also, for any $y \in B(x,r)$, we either have $d(x,y) \geq s/6$ or $d(x,z) \geq s/6$. In the first case we get by a similar analysis $|f(y)-f(x)| \leq D^5 C_3 A s$. In the latter we get $|f(z)-f(y)| \leq D^5 C_3 A s$, and $|f(y)-f(x)| \leq 2D^5 C_3 A s$. Thus, in any case we can integrate and get 
$$\vint_{B(x,r)} |f(y)-f(x)| ~d\mu_y \leq 4D^5C_3 A r.$$

Together with Lemma \ref{average} and the choice $a=f(x)$ this completes the proof.

\QED

Similar techniques give a modulus estimate. For any $x,y \in X$, and any $C$, define the collection of curves.
$$\Gamma_{x,y,C} \defeq \{ \gamma \co [0,1] \to X | \gamma(0)=x, \gamma(1)=y, \len(\gamma) \leq Cd(x,y)\}$$

We call a non-negative measurable function $\rho$ admissible for the family $\Gamma_{x,y,C}$ if for any $\gamma \in \Gamma_{x,y,C}$ we have
$$\int_\gamma \rho ~ds \geq 1.$$
The $p$-modulus of the curve family (centered  at $x$, and scale $s$) is then defined as
$$\Mod_p (\Gamma_{x,y,C}) \defeq \inf_{\rho \text{ is admissible for } \Gamma_{x,y,C}} \int \rho^p ~d\mu.$$
For a more detailed discussion on modulus, see \cite{vaisala2006}.

By the previous arguments for a finely $\alpha$-connected space we can show the following lower bound for $p>\frac{1}{\alpha}$.

\begin{theorem}\label{modulusest} If $(X,d,\mu)$ is a complete $(D,r_0)$-measure doubling metric measure space which is also locally finely $(\alpha, r_0)$-connected with parameters $(C_1, C_2)$, then for any $x,y \in X$ with $d(x,y)=r < r_0/ (8C_1)$, any $p>\frac{1}{\alpha}$ we have
$$\Mod_p (\Gamma_{x,y,C_3}) \geq \frac{\mu(B(x,2C_1r))}{2C_3^p r^p},$$
where $C_3$ is as in Proposition \ref{prop:smallintegral}.
\end{theorem}

\noindent \textbf{Proof:} Immediate corollary of Proposition \ref{prop:smallintegral}. If $\rho$ is admissible, and 
$$\int_{X}  \frac{1}{2C_3^pr^p} \rho^p \, d\mu \leq \frac{\mu(B(x,2C_1r))}{2C_3^p r^p},$$
then also
$$ \vint_{B(x,2C_1r)} \rho^p ~d\mu \leq \frac{1}{2C_3^pr^p},$$
and there exists a Lipschitz curve $\gamma$ connecting $x$ to $y$ of length at most $C_3r$ and 
$$\int_\gamma \rho ~ds < \frac{1}{C_3r} C_3r \leq 1,$$
which contradicts the admissibility  of $\rho$. Thus, the original modulus estimate must hold.

\QED 

\noindent \textbf{Remark:} The results of this section are tight in terms of the range of $p$. We give two examples. Taking two copies of $\R^n$ glued along the origin (in fact any subset) we obtain a space that is finely $\frac{1}{n}$-connected and a simple analysis based on modulus estimates tells us that it only admits a $(1,p)$-Poincar\'e inequality for $p>n$. Another example arises from $A_s$ weights $w$ on $\R$. We know that such weights are finely $\frac{1}{s}$-connected \cite{stein2016harmonic}, and thus they admit a Poincar\'e inequality for $p>s$. %These examples suggest that the fine $\alpha$-connectivity measures the stable regime of Poincar\'e inequalities, i.e. those $p$ such that $(1,p)$-inequalities hold on glued spaces or spaces with controlled density changes. 

For purposes of the following proof recall the definition of the upper gradient in Definition \ref{uppergrad} and non-homogenous Poincar\'e inequality in Definition \ref{nonhomopoincare}.  

% \begin{definition} Let $(X,d,\mu)$ be a metric measure space and $C\geq 1$ some constant. Let $\Phi, \Psi\co [0, \infty) \to [0, \infty)$ be increasing functions with the following properties.
% 
% \begin{itemize}
% \item $\lim_{t \to 0}\Phi(t)=\Phi(0)=0$
% \item $\lim_{t \to 0}\Psi(t)=\Psi(0)=0$
% \end{itemize}
% We say that $(X,d,\mu)$ satisfies a non-homogeneous $(\Phi,\Psi,C,r_0)$-Poincar\'e-inequality if for every $2$-Lipschitz\footnote{The constant $2$ is only used to simplify arguments below. Any fixed bound could be used. Also, by adjusting the right-hand side by scaling with the Lipschitz constant of $f$, the inequality could be made homogeneous.} function $f\co X \to \R$, every upper gradient $g \co X \to \R$ for $f$ and every ball $B(x,r) \subset X$ with $r<r_0$ we have
% 
% $$\vint_{B(x,r)} |f-f_{B(x,r)}| ~d\mu \leq r \Psi\left( \vint_{B(x,Cr)} \Phi\circ g~d\mu \right).$$
% \end{definition}

We can now prove Theorem \ref{thm:nohomo}.

%\begin{theoremN}[\ref{thm:nohomo}] Suppose that $(X,d,\mu)$ is $(D,r_0)$-doubling metric measure space and satisfies a local non-homogeneous $(\Phi, \Psi, C,r_0)$-Poincar\'e-inequality. Then $(X,d,\mu)$ is locally $(C,\delta,\epsilon)$-connected and moreover admits a true $(1,q)$-Poincar\'e inequality for some $q>1$. All the variables are quantitative in the parameters.  
%\end{theoremN}

\noindent \textbf{Proof of Theorem \ref{thm:nohomo}:} By adding an increasing positive function to $\Psi,\Phi$ we can assume that $\Psi(t)>0,\Phi(t)>0$ for $t>0$. By weakening the assumption, we can take $\Psi, \Phi$ to be upper semi-continuous. Assume next that for any $2$-Lipschitz $f$ and $g$ an upper gradient for $f$ and any ball $B(x,r)$ with $r<r_0$
$$\vint_{B(x,r)} |f-f_{B(x,r)}| ~d\mu \leq r \Psi\left( \vint_{B(x,Cr)} \Phi(g) ~d\mu \right).$$
Define the right inverses $\xi(t) \defeq \inf \{s | \Psi(s) \geq t \}$. Then we have $\xi(\Psi(t)) \leq t$ and thus
$$\xi \left( \frac{1}{r} \vint_{B(x,r)} |f-f_{B(x,r)}| ~d\mu \right) \leq \vint_{B(x,r)} \Phi(g) ~d\mu.$$

Since $\lim_{t \to 0} \Phi(t) =0$ we can choose a function $\sigma(t)$ such that for any $0<s \leq \sigma(t)$ we have $$\Phi(s) \leq t.$$

We will show local $(B,\frac{1}{2},\epsilon)$--connectivity for $0<\epsilon$ small enough for a specific $B$. Our choice will be $$B = \max\left(2C,\left[\sigma\left(2^{-1}\xi\left(\frac{1}{20D^5} \right) \right)\right]^{-1} \right).$$ Suppose $X$ is not $(B,\frac{1}{2},\epsilon,\frac{r_0}{8B})$--connected. If this were the case, we could choose $x,y \in X$ with $r = d(x,y)<r_0(8B)^{-1}$, and an obstacle $E$ such that
$$\frac{\mu(E \cap B(x,Br))}{\mu(B(x,Br))} \leq \epsilon,$$
and such that for any Lipschitz curve fragment $\gamma\co K \to B(x,Br)$ connecting $x$ to $y$ such that $\gamma^{-1}(E) \subset \{\min(K), \max(K)\}$ with $\len(\gamma)\leq Br$ we'd have $\gap(\gamma) >\frac{r}{2}$. Define
$$\rho(z) \defeq \inf_{\gamma\co K \to X} \frac{1}{2B} \len(\gamma) + \gap(\gamma)(E) .$$

Here, the infimum is taken over all compact $K \subset \R$ and all Lipschitz curve fragments $\gamma\co K \to X$ connecting $x$ to $z$ with $\gamma^{-1}(E) \subset \{\min(K), \max(K)\}$. It is easy to see that $\rho$ is $2$-Lipschitz, and has upper gradient $g=1_E + \frac{1}{B}$. The Lipschitz bound follows for $v,w \in X$ by the following argument. For simplicity assume $v,w \not\in E$. Assume $\gamma_{x,v} \co K \to X$ is any curve connecting $x$ to $v$. Then, we can define another curve fragment by $\gamma_{x,w} \co K \cup \{\max(K) + d(v,w)\} \to X$ by setting $\gamma_{x,w}|_K = \gamma_{x,v}$, and $\gamma_{x,w}(\max(K) + d(v,w)) = w$. This connects $x$ to $w$, and we have
$$\rho(w) \leq \frac{1}{2B}\len(\gamma_{x,w}) + \gap(\gamma_{x,w})(E) \leq \frac{1}{2B}\len(\gamma_{x,v}) + \gap(\gamma_{x,v})(E) + (1+ \frac{1}{2B})d(v,w).$$
Infimizing over all curve fragments $\gamma_{x,v}$, and switching the roles of $v$ and $w$ in the previous argument, gives the desired Lipschitz bound for $\rho$. If either $v,w \in E$, then the argument is slightly more subtle, as we must excise a neighborhood of $\max(K)$ and restrict $\gamma_{x,w}$ to the set $K \cup \{\max(K) + d(v,w)\} \setminus (\max(K)-\tau, \max(K)+\tau)$, and let $\tau \to 0$.

Also, $\rho(0)=0$. By assumption, there is no curve fragment $\gamma$ with $\gamma^{-1}(E) \subset \{\min(K), \max(K)\}$ connecting $x$ to $y$ with $\len(\gamma)\leq Br$,  and $\gap(\gamma)(E)  \leq \frac{r}{2}$, thus either of these inequalities fails, and we have
$$\rho(y) \geq \frac{r}{2}.$$
By the Lipschitz-condition, we have
$$\rho(z) \leq \frac{r}{10},$$
for $z \in B(x,\frac{r}{20})$ and
$$\rho(z) \geq \frac{2r}{5}$$
for $z \in B(y,\frac{r}{20})$. Thus, by a doubling estimate $\frac{1}{2r}\vint_{B(x,2r)} |f-f_{B(x,2)}| d\mu \geq \frac{1}{20D^5}$. On the other hand
\begin{eqnarray*}
\vint_{B(x,2Cr)} \Phi(g) ~d\mu &\leq & \Phi\left(\frac{1}{B}\right) + \Phi(1) \frac{\mu(E \cap B(x,Br))}{\mu(B(x,2Cr))} \\
                                               &\leq & \Phi\left(\sigma\left(\xi\left(\frac{1}{20D^5} \right)/2\right) \right) + D^{\log_2(B/C) +2} \Phi(1) \epsilon \\
                                                &\leq & \xi\left(\frac{1}{20D^5} \right)/2+ D^{\log_2(B/C) +2} \Phi(1) \epsilon. \\
\end{eqnarray*}

Combining all the observations, we would get
$$\xi\left(\frac{1}{20D^5}\right) \leq  \xi\left(\frac{1}{20D^5}\right)/2 + D^{\log_2(B/C) +2} \Phi(1) \epsilon.$$
If $\epsilon <  \frac{\xi\left(\frac{1}{20D^5}\right)}{2\Phi(1)D^{\log_2(B/C) +2}}$, this inequality fails deriving a contradiction and thus proving the conclusion.

\QED

% \begin{theoremN}[\ref{thm:contheorem}] A $(D,r_0)$-doubling complete metric measure space $(X,d,\mu)$ admits a local $(1,p)$-Poincar\'e-inequality for some $p \in [1, \infty)$ if and only if it is locally $(C,\delta,\epsilon)$-connected for some $\delta,\epsilon \in (0,1)$. Both directions of the theorem are quantitative in the respective parameters.
% \end{theoremN}

Finally, we conclude this section with a proof of Theorem \ref{thm:contheorem}. \\

 \noindent \textbf{Proof:} The converse statement was proved in Theorem \ref{maintheorem}. Thus, we only prove that a measure-doubling space admitting a $(1,p)$-Poincar\'e-inequality for some $p \in [1,\infty)$ also admits $(C,\delta,\epsilon)$--connectivity for some constants. This is a particular case of the more general Theorem \ref{thm:nohomo} by setting $\Phi(x)=x^p$ and $\Psi(x)=C_{PI}x^{\frac{1}{p}}$. Here, $C_{PI}$ is the constant in the Poincar\'e inequality. However, there is an issue related to whether the Poincar\'e inequality is assumed for Lipschitz functions or all continuous functions and their upper gradients. We use the latter for non-homogeneous Poincar\'e inequalities. In the case of a complete metric space these are equivalent by \cite[Theorem 2]{keith2003modulus}\footnote{The application of \cite[Theorem 2]{keith2003modulus} is somewhat subtle. Keith assumes a global Poincar\'e inequality and global doubling, while here we work with local Poincar\'e ienqualities and local doubling. Worse  yet, our space is complete but only locally compact at a definite scale. Keith's proof can be fully localized to apply in this setting in the sense that the scales $r_{i}$ of his four claims, $i=1, \dots 4$, depend quantitatively on each other. Presenting this in detail here would be distracting and wouldn't offer any new ideas or mathematics. We are also not aware of a reference where this step is done in detail. However, \cite[Theorem 8.4.1]{heinonen2015sobolev} presents an argument that is easily localizable. The crucial step involves a pair of points $d(x,y)=r$ and requires using certain bounded length paths which can be contained in a ball $B(x, Cr)$ for some $C$, which is precompact by local doubling if $r$ is sufficiently small. Each of the other estimates there, except for quasiconvexity, is explicitly local. Further, the proofs use of \cite[Theorem 8.1.7]{heinonen2015sobolev} can also be localized since the estimates are only needed up to a certain scale. The quasiconvexity proof in \cite[Theorem 8.3.2]{heinonen2015sobolev} can also be localized since the construction of a quasigeodesic between to points $x,y$ with $d(x,y)=r$ involves the Poincar\'e inequality only at scales comparable to $r$, and the pre-compactness of some ball $B(x, Cr)$, which can be obtained from the local doubling bound. The proof of local quasiconvexity could also be derived fairly directly from the proof in \cite{bjornsemilocal}, although they use a different Poincar\'e inequality as an assumption. We leave a detailed examination and explicit quantification of these arguments to the interested reader.}.
 %We need that our local $(1,p)$-Poincar\'e inequality for Lipschitz functions implies the $(1,p)$-Poincar\'e inequality for continuous functions and their upper gradients. The $(1,p)$-Poincar\'e inequality at a scale $r_0$ implies quasiconvexity up to scale $r_1\sim r_0$ (by Lemma \cite[Lemma 9]{keith2003modulus}), and thus a modulus estimate of the family of curves connecting pairs of points with respect to bounded $\rho$. This modulus bound implies a modulus bound for curves of bounded length, since long curves have small modulus. Finally, if the pair of points are at a scale $r_2 \sim r_1$, then we get that all of these curves lie in a proper subset of $X$ (since $X$ is locally doubling and complete, and thus locally compact). Then, the propositions of \cite{keith2003modulus} and \cite[Proposition 2.17]{heinonen1998quasiconformal} can be used to give a modulus bound with respect to all non-negative $\rho$ for this family relative to a compact ball which contains it. Finally, this modulus bound corresponds to a capacity bound \cite[Proposition 2.17]{heinonen1998quasiconformal}, and can then be translated to a pointwise Poincar\'e inequality. By the argument in  \cite[Lemma 5.15.]{heinonen1998quasiconformal}, this gives a Poincar\'e inequality. Each of these steps involves only using properness, doubling, or the Poincar\'e inequality at a definite multiple of the initial scale. However, a completely formal proof would need to rewrite these steps with taking care of the associated scales.}. 
\QED

\section{Corollaries and Applications} \label{sec:corapp}

If $X$ is a complete and proper metric space, we denote by $\Gamma(X)$ the space of geodesics parametrized by the interval $[0,1]$. This space can be made into a complete proper metric space by using the distance
$d(\gamma,\gamma')=\sup_{t \in [0,1]} d(\gamma(t), \gamma'(t)).$ Further, define for each $t \in [0,1]$ the map $e_t \co \Gamma(X) \to X$ by $e_t(\gamma) \defeq \gamma(t)$.

\begin{definition} Let $\rho$ be a decreasing function with $\rho(1) \geq 1$. A proper geodesic metric measure space $(X,d,\mu)$ is called a measure $\rho(t)$-contraction space, if for every $x \in X$ and every $B(y,r) \subset X$, there is a Borel probability measure $\Pi$ on $\Gamma(X)$ such that $e_0(\gamma) = x$, $\Pi$-almost surely, and $e_1^*\Pi = \frac{\mu|_{B(y,r)}}{\mu(B(y,r))}$, and for any $t \in (0,1)$ we have
$$e_t^*\Pi \leq \frac{\rho(t)}{\mu(B(y,r))} d\mu.$$
\end{definition}

This class of metric measure spaces includes the MCP-spaces by Ohta \cite{ohta}, as well as $CD(K,N)$ spaces, and any Ricci-limit, or Ricci-bounded space. In fact, this class of spaces is more general than the MCP-spaces considered by Ohta. Until now, it was unknown if such spaces possess Poincar\'e inequalities. We show as a corollary of Theorems \ref{thm:contheorem} and \ref{thm:alphatheorem} that this indeed holds. The definition could be modified to permit a probability measure on $C$-quasiconvex paths, but we don't need this result here.

\begin{corollary} Any measure $\rho(t)$-contraction space $(X,d,\mu)$ admits a $(1,p)$-Poin\-ca\-r\'e inequality for some $p \in [1,\infty)$.
\end{corollary}

\textbf{Proof:} Clearly any such space is $\rho(1/2)$-doubling. We will show that the space is $\left(2, \delta, \frac{\delta}{100\rho(1/2)^6 \rho(\delta/3)}\right)$--connected. Consider any $x,y \in X$, and $d(x,y)=r$, and a Borel set $E \subset B(x,2r)$ with $\mu(E) \leq \frac{\delta}{50\rho(1/2)^6 \rho(\delta/3)} \mu(B(x,2r))$. It is sufficient to find a Lipschitz curve $\gamma \co [0,L] \to X$ connecting $x $ to $y$ with $L \leq 2r$ and $\int_\gamma 1_E \, ds < \delta d(x,y)$, since then a curve fragment can be defined by restricting $\gamma$ to an appropriate compact set $K \subset [0,L] \setminus \gamma^{-1}(E)$. 

Consider a midpoint $z$ on a geodesic connecting $x$ to $y$ and the ball $B(z,r/2)$. By $\rho(1/2)$-doubling and the assumption we have $\mu(E) \leq \frac{\delta}{50\rho(\delta/3)} \mu(B(z,r/2))$.

Apply the definition of a $\rho$-contraction space to the point $x$ and the ball $B(z,r/2)$. This constructs a probability measure $\Pi$ on curves from $x$ to $B(z,r/2)$, such that
$$e_t^*\Pi \leq  \frac{\rho(t)}{\mu(B(z,r/2))}  d\mu.$$

Note that $\gamma(\delta/3) \in B(x,\frac{\delta r}{3})$ $\Pi$-almost surely. Thus,
\begin{eqnarray*}
 \int_\Gamma \int_{\delta/3}^1 1_{E}(\gamma(t))\, dt \, d\Pi_\gamma
 &\leq&  \frac{\rho(\delta/3)}{\mu(B(z,\frac{r}{2}))} \mu(E) \\
&\leq&  \frac{\delta }{50}.
\end{eqnarray*}

Thus, with probability strictly bigger than $\frac{1}{2}$ we have $\int_{\delta/3}^1 1_{E}(\gamma(t))\, dt  \leq \frac{\delta}{20}$. Let $\Gamma_x$ be the collection of these geodesics. For any such geodesic $\gamma_x \in \Gamma_x$ we have $\len(\gamma_x) \leq r$ and  
$$\int_{\gamma_x} 1_E\, ds \leq r \int_0^1 1_E(\gamma_x(t)) \, dt \leq \left(\frac{\delta}{3} + \frac{\delta}{20} \right) r < \frac{\delta}{2}r.$$
Let $S_x$ be the set of their end points in $B(z,r/2)$. Then $\mu(S_x)/\mu(B(z,r/2)) =  e_1^*(\Pi)(\Gamma_x)> \frac{1}{2}$. Similarly, by symmetry such a set $S_y$ with respect to $y$ can be constructed. By the volume estimates $S_x \cap S_y \neq \emptyset$. We can thus find a common point $w \in B(z, r/2)$, and curve fragments from $x$ to $w$ and from $w$ to $y$. These curves, when concatenated, give a Lipschitz curve $\gamma \co [0,2] \to X$ connecting $x$ to $y$ with $\len(\gamma) \leq 2r$ and
$$\int_\gamma 1_E \, ds < \delta r.$$
Finally, using Borel regularity of Lebesgue measure on the real line we can find a compact set $K$ in the complement of $\gamma^{-1}(E)$ of almost maximal measure. Then $\gamma|_{K \cup \{0,2\}}$ is our desired curve fragment.
%This curve fragment once dilated by $4r$ results in the desired $1$-Lipschitz curve fragment.

\QED

The conclusion of Theorem \ref{thm:mcp} follows from the following corollary.

\begin{corollary} If $\rho(t)=t^{-n}$, and $(X,d,\mu)$ is a $\rho$-contraction space, then it satisfies a Poincar\'e inequality for $p>\frac{1}{n+1}$. Since $MCP(0,n)$-spaces satisfy this assumption we conclude that they admit such a Poincar\'e inequality. The same holds for every $MCP(K,n)$-space, but with a local Poincar\'e inequality.
\end{corollary}

\textbf{Proof:} Since for any $\delta>0$ the space is
$$\left(2, \delta, \frac{\delta}{100\rho(1/2)^6 \rho(\delta/3)}\right)-\text{connected.}$$
Plug in $\delta = \tau^{\frac{1}{n+1}}$ and $\rho(t)=t^{-n}$. This gives that the space is finely $\alpha$-connected for $\alpha=\frac{1}{n+1}$. The Poincar\'e inequality follows from Theorem \ref{thm:alphatheorem}. For $K>0$ the same result is obvious and $K<0$ we modify it slightly to obtain that the space is locally finely $(\alpha,r_0)$-connected (see \cite{ohta} for the modified $\rho$).
\QED

\noindent \textbf{Remark:} It would seem that this result should be true for $p>n$, and probably $p \geq 1$, but our methods do not yield these sharper results. 

We next prove a result on the existence of Poincar\'e inequalities on spaces deformed by Muckenhoupt-weights (recall Definition \ref{def:muckenhoupt}). In a geodesic doubling space any generalized $A_\infty(\mu)$-measure is doubling. The geodesic assumption below can be slightly weakened but counter-examples can be constructed if $(X,d,\mu)$ is just doubling.
\begin{lemma}\label{doubling} Let $(X,d,\mu)$ be a $D$-measure doubling geodesic metric measure space. Then any Radon measure $\nu \in A_{\infty}(\mu)$ is doubling.
\end{lemma}

\noindent \textbf{Proof:} Let the Muckenhoupt parameters, as in Definition \ref{def:muckenhoupt}, corresponding to $\nu$ be $\delta,\epsilon$. By \cite[Lemma 3.3]{colding1998liouville}\footnote{Proof works for geodesic metric spaces, although it is only stated for manifolds. See also \cite[Proposition 6.12]{ChDiff99}.} for $\eta\leq 2^{\frac{\log_2(1-\epsilon)}{\log_2(1+D^{-5})}-2}$.
$$\mu(B(x,r) \setminus B(x,(1-\eta)r)) < (1-\epsilon)\mu(B(x,r)).$$

Thus $\mu(B(x,(1-\eta)r))>\epsilon\mu(B(x,r))$, and from the Muckenhoupt condition \ref{def:muckenhoupt}, we get $\nu(B(x,(1-\eta)r))>\delta\nu(B(x,r))$. Iterating this estimate $N=-1/\log_2(1-\eta)+1$ times, we get
$$\nu(B(x,r/2)) \geq \delta^N \nu(B(x,r)).$$

This gives the desired doubling property.

\QED

We need another technical Lemma from the paper of Kansanen and Korte \cite{kansanen2011strong}. This is a metric space generalization of the classical property of self-improvement for $A_\infty$-weights in Euclidean space \cite{stein2016harmonic}.

\begin{lemma}\label{expand}(\cite[Corollary 4.19]{kansanen2011strong}) Let $(X,d,\mu)$ be a $D$-measure doubling geodesic metric measure space. Then for any Radon measure $\nu \in A_{\infty}(\mu)$ and for any $0<\tau<1$ there exists a $0<\delta<1$, such that for any $B(x,s)$ and any $E \subset B(x,s)$
$$\nu(E) \leq \delta \nu(B(x,s)) \Longrightarrow \mu(E) \leq \tau \mu(B(x,r)).$$

\end{lemma}

\begin{theorem} If $(X,d,\mu)$ is a geodesic PI-space and if $\nu \in A_{\infty}(\mu)$, then $(X,d,\nu)$ is also a PI-space
\end{theorem}

\textbf{Proof:} By Lemma \ref{expand} we see that $(X,d,\nu)$ is doubling, and from Theorem \ref{thm:contheorem} we get that  $(X,d,\mu)$ is $(C,\frac{1}{2},\epsilon_\mu)$--connected for some $C,\epsilon_\mu$. Further by Lemma \ref{expand} we have a $\epsilon_\nu$ such that for any Borel $E \subset B(x,r) \subset X$
$$\nu(E) \leq \epsilon_{\nu} \nu(B(x,r)) \Longrightarrow \mu(E) \leq \epsilon_\mu \mu(B(x,r)).$$

Thus, since any obstacle $E$ with volume density $\epsilon_\nu$ with respect to $\nu$ will be an obstacle with volume density $\epsilon_\mu$ with respect to $\mu$, it is easy to verify that $(X,d,\nu)$ is $(C,\delta,\epsilon_\nu)$--connected. All the parameters can be made quantitative. 
\QED

\section{PI-rectifiability, asymptotic connectivity and differentiability} \label{rectifdiff}

\subsection{Main theorems}

In order to state our main theorem, we need to define the relevant notions. We are interested in understanding when a given metric measure space is PI-rectifiable.

\begin{definition} \label{def:PIrect}A metric measure space $(X,d,\mu)$ is PI-rectifiable if there is a decomposition into countably many measurable sets $U_i,N \subset X$ such that
$$X = \bigcup_i U_i \cup N,$$
where $\mu(N)=0, \mu(U_i)>0$ and there exist isometric and measure preserving embeddings $\iota_i \co U_i \to \overline{U_i}$ with $(\overline{U}_i, \overline{d}_i, \overline{\mu}_i)$ PI-spaces (with possibly very different constants).
\end{definition}

We define a differentiability space as in \cite{keith04diff, ChDiff99}.
\begin{definition}\label{diffspace} A metric measure space $(X,d,\mu)$ is called a differentiability space (of analytic dimension $\leq M$) if there exist countably many measurable $U_i ,N \subset X$ and Lipschitz functions $\phi_i\co U_i \to \R^{n_i}$ (with $n_i \leq M$) such that $\mu(N)=0$, 
$$X=\bigcup_i U_i \cup N,$$
and such that for any Lipschitz function $f\co X \to \R$, for every $i$ and $\mu$-almost every $x \in U_i$, there exists a \textit{unique} linear map $D^{\phi_i}f(x)\co \R^{n_i} \to \R$ such that
$$f(y) = f(x) +D^{\phi_i}f(x)(\phi_i(y)-\phi_i(x))+o(d(x,y)).$$ 
\end{definition}

 A stronger definition is obtained by assuming differentiability of Lipschitz functions with certain infinite dimensional targets \cite{bate2015geometry}.

\begin{definition}\label{rnpdiffspace} A metric measure space $(X,d,\mu)$ is a RNP-(Lipschitz) differentiability  space (of analytic dimension $\leq M$) if there exist countably many measurable $U_i,N \subset X$ and associated Lipschitz functions $\phi_i\co U_i \to \R^{n_i}$ (with $n_i \leq M$) such that $\mu(N)=0$, 
$$X=\bigcup_i U_i \cup N,$$
and such that for any Banach space $V$ with the Radon-Nikodym property and any Lipschitz function $f\co X \to V$, for every $i$ and $\mu$-almost every $x \in U_i$ there exists a \textit{unique} linear map $D^{\phi_i}f(x)\co \R^{n_i} \to V$ such that
$$f(y) = f(x) +D^{\phi_i}f(x)(\phi_i(y)-\phi_i(x))+o(d(x,y)).$$ 
\end{definition}

\begin{definition} A Banach space $V$ has the Radon-Nikodym property if every Lipschitz function $f:[0,1] \to V$ is almost everywhere differentiable.
\end{definition}

For several equivalent and sufficient conditions, consult for example \cite{pisier2016martingales} or the references in \cite{cheeger2007characterization}.

\begin{definition} A metric measure space $(X,d,\mu)$ is called asymptotically well-connected if for all $\delta \in (0,1)$ and almost every point $x$ the exist $r_x>0$, $\epsilon_x \in (0,1)$ and $1 \leq C_x$ such that for any $y \in X$ such that $d(x,y)<r_x$ the pair $(x,y)$ is $(C_x,\delta,\epsilon_x)$--connected.
\end{definition}

We can now state a general theorem that characterizes PI-rectifiable metric measure spaces. This is proven later in this section.

\begin{theorem}\label{generalrect} A complete metric measure space $(X,d,\mu)$ is a PI-rectifiable if and only if it is asymptotically doubling and asymptotically well-connected.
\end{theorem}

Using this we can prove the PI-rectifiability result.

\begin{theoremN}[\ref{thm:PIrect}] A complete metric measure space $(X,d,\mu)$ is a RNP-Lipschitz differentiability space if and only if it is PI-rectifiable and every porous set has measure zero.
\end{theoremN}

\noindent \textbf{Proof of Theorem \ref{thm:PIrect}:} By \cite{bate12diff} RNP-Lipschitz differentiability spaces, being Lipschitz differentiability spaces, are asymptotically doubling. Further, by \cite[Lemma 3.5]{bate2015geometry} which is rephrased above in Lemma \ref{lem:rephrased} we have that a RNP-Lipschitz differentiability space is asymptotically well-connected. See also the discussion preceding Lemma \ref{lem:rephrased} for some subtle points about the differences in terminology in \cite{bate2015geometry}. As for the converse, we refer to the main result in \cite{Ch09diff}, where it is shown that a PI-space is a RNP-Lipschitz differentiability space. It is then trivial to conclude that a PI-rectifiable space is also a RNP-Lipschitz differentiability space if every porous set has measure zero. See the discussion at the end of the introduction in \cite{bate2014characterizations} as well as  the arguments in \cite{bate2013differentiability} for additional details. 
\QED

The proof of Theorem \ref{generalrect} is based on a general ``thickening Lemma''. To see this, we outline the proof. By using measure theory arguments, the asymptotic connectedness can be used to produce subsets such that the space $X$ satisfies some doubling and connectivity estimates uniformly along such sets. Then these subsets are enlarged, or ``thickened'' to   improve a relative form of connectivity to an intrinsic form of connectivity. This step is included in the following theorem.

\begin{theoremN}[\ref{relative}]  Let $r_0>0$ be arbitrary. Assume $(X,d,\mu)$ is a metric measure space and subsets $K \subset A \subset X$ are given, where $A$ is measurable and $K$ is compact. Assume further that $X$ is $(D,r_0)$-doubling along $A$, $A$ is uniformly $(\frac{1}{2},r_0)$-dense in $X$ along $K$, and $A$ with the restricted measure and distance is locally $(C, 2^{-60}, \epsilon, r_0)$-connected along $K$. There exist constants $\overline{C}, \overline{\epsilon}, \nD >0$, and a complete metric space $\nK$ which is locally $(\overline{D}, 2^{-40}r_0)$-doubling and $(\overline{C}, \frac{1}{2}, \overline{\epsilon}, r_0 2^{-330}/C)$--connected, and an isometry $\iota \co K \to \nK$ which preserves the measure. In particular, the resulting metric measure space $\nK$ is a PI-space.
\end{theoremN}

The enlarged space $\nK$ is obtained by attaching a tree-like metric space $T$ to $K$. We remark, that there is no unique ``thickening'' constructing, but it is our belief that gluing a tree-like metric space is the easiest way to obtain the desired conclusions. Additionally, the previous Lemma, when applied to subsets $K \subset X$  of PI-spaces $X$, produces new examples of PI-spaces from old ones which consist of possibly disconnected sets $K$ with tree-like graphs glued onto them. % where the possible disconnectedness of $K$ is repaired by the glued space $T$.

Our primary goal in constructing $T$ is increasing connectivity and making $\nK$ quasiconvex. That goal could be attained by attaching an edge for every pair of points in $K$. This, however, would fail to be doubling. Thus, we need more care in constructing $T$. The main issues are as follows. %Our construction is focused on the following issues.

\begin{enumerate}
\item The resulting space needs to be a complete locally compact metric measure space. In particular, the added intervals need measures associated with them.
\item The measures of the added intervals need to be controlled by $\mu$ in order to preserve doubling. 
\item Doubling needs to be controlled. We shouldn't add too many intervals at any given scale and location. 
%\item The intervals should themselves be ``simple'' PI-spaces, i.e. the measures should be related to Lebesgue measure.
\item A curve fragment in $A$ between points in $K$ should be replaceable, up to small measure, by  a possibly somewhat longer curve fragment in the glued space. 
\end{enumerate}  

A natural construction to obtain these goals arises from a modification of a so called hyperbolic filling. This construction first appeared in \cite{bourdon2003}, and later in \cite{bonk2014sobolev}.

We will first indicate how Theorem \ref{generalrect} follows from Theorem \ref{relative}.

\noindent \textbf{Proof of Theorem \ref{generalrect}:} By general measure theory arguments, we can construct sets $X = \bigcup_i V_i \cup N = \bigcup_{i,j} K_i^j \cup N'$, such that each $V_i$ are measurable, $K_i^j \subset V_i$ are compact, $\mu(N)=\mu(N')=0$ and the following uniform estimates hold.

\begin{itemize}
\item $X$ is  $(D_i,r_i)$-doubling along $V_i$.
\item $X$ is $(C_i, 2^{-60}, \epsilon_i,r_i)$--connected along $V_i$ for some $\epsilon_i < 1$.
\item $V_i$ is $(1-\epsilon_i/2, r^j_i)$-dense in $X$ along $K_i^j$.
\end{itemize}

For the definitions of these concepts see the beginning of Section \ref{notconvention}.

Since the set $V_i$ is uniformly $(1-\epsilon_i/2, r_i^j)$-dense in $X$ along $K_i^j$, directly we obtain that $V_i$ with its restricted measure and metric is $(C_i, 2^{-60}, \epsilon_i/2,\min(r_i,r_i^j)/C_i)$--connected along $K_i^j$. To see this, verify definition \ref{def:conbility} for every obstacle $E$ by adjoining the complement $X \setminus V_i$ to any set $E$ considered. By the density bound we see that $E \cup X \setminus V_i$ satisfies a slightly worse density bound, and thus is a permissible obstacle.

Theorem \ref{relative} can now be applied to $A=V_i$, $K=K_i^j$ and $\epsilon=\epsilon_i/2, r_0=\min(r_i, r_i^j)/C_i = r_0^{i,j}, C = C_i$ and $D = D_i$ to obtain isometric and measure preserving embeddings $\iota \co K_i^j \to \nK_i^j$ to metric measure spaces $(\nK_i^j, d_i^j, \mu_i^j)$. Further, by the same theorem, there exist positive constants $\overline{D}_{i,j}, \overline{C}_{i,j}$ and $\overline{\epsilon}_{i,j}$ such that the metric measure spaces $(\nK_i^j, d_i^j, \mu_i^j)$ are locally $(\overline{D}_{i,j}, r_0^{i,j} 2^{-60})$-doubling and locally $(\overline{C}_{i,j}, 2^{-1}, \overline{\epsilon}_{i,j}, r_0^{i,j}2^{-320})$--connected. Thus, by Theorem \ref{thm:contheorem} each $\nK_i^j$ is a PI-space. This completes the proof of PI-rectifiability.

\QED

\subsection{Construction}
\label{construction}

Continue next with the proof of Theorem \ref{relative}. First, we present the construction, which is followed by a separate subsection containing the proofs of the relative connectivity and doubling properties. \\ \\

\noindent \textbf{``Thickening'' construction for Theorem \ref{relative}:} Assume for simplicity that $D\geq 2$. Re-scale so that $r_0=2^{20}$. We will define $\nK \defeq K \cup T$, where $T$ arises as the metric space induced by a metric graph $(V,E)$, and $K$ naturally sits at the boundary of $T$. The vertices will arise from two types of points, net points of $K$ and gap points that correspond to center points in $A \setminus K$.

First define the $2^{-j}$-nets $N_j$ of $K$ for $j \geq 0$, i.e. maximal collections of points such that $d(n,m)\geq 2^{-j}$ for $n,m\in N$. In order to make some notation below consistent define $N_{-1} \defeq \emptyset$. Define these inductively such that $N_{j} \subset N_{j+1}$, and define  the collection of all net points $N\defeq \bigcup_j N_j$. Further, define the scale function as $\scl_N(n)\defeq 2^{-i}$ if $n \in N_i \setminus N_{i-1}$

Using the Vitali covering Lemma \cite{stein2016harmonic,federer}, choose ``gap'' points $g_k^m \in (\NN_1(K) \setminus K)\cap A$, indexed by $(k ,m) \in U \subset \N_0 \times \N_0$ , such that
\begin{eqnarray}
(\NN_1(K) \setminus K) \cap A &\subset&  \bigcup_{(k,m) \in U} B\left(g_k^m, 2^{-k-10}\right), \label{whitneycover}\\
2^{-k-1} < &d(g_k^m, K)&\leq 2^{-k} \label{whitneydist}
\end{eqnarray}
and for $(k,m) \neq (k',m')$ with $(k,m), (k',m') \in U$
\begin{equation}
B\left(g_k^m, 2^{-k-15}\right) \cap B\left(g_{k'}^{m'}, 2^{-k'-15}\right)=\emptyset. \label{whitneydisjoint}
\end{equation}  

Recall that $\NN_1(K) \defeq \{x \in X | d(x,K) < 1\}$. This type of covering is also known as a Whitney-covering. Given such a covering, define the gap points at scale $2^{-k}$ by $G_k\defeq \left\{g_k^m | (k,m) \in U \right\}$ and $G \defeq \bigcup_k G_k$ the collection of all gap points.  We will define a scale function by $\scl_G(g)\defeq 2^{-k}$ for $g \in G_k$.

Define the vertex sets as $V \defeq V_G \cup V_N$, where $V_G$ arises from $G$ and $V_N$ arises from $N$. Vertices will be defined as pairs $(x,r)$ with $x \in X$ and $r>0$. First take
$$V_G \defeq \left\{(g,2^{-k}) | g \in G_k \right\}.$$
These will also be referred to as bridge-points because they are used to repair the gap corresponding to the center $g$ at scale $2^{-k}$.

Next, define for each $n \in N$ the maximal scale at which there is a near-by gap point.
\begin{equation}
l(n) \defeq \max\left\{2^{-k} \left| \exists g \in G_k : d(g,n) \leq 2^{4-k}\text{ and } \scl_G(g) =2^{-k}\leq \scl_N(n)  \right. \right\}. \label{lenscale}
\end{equation}
If there is no such point, we will define $l(n) \defeq 0$. We define the vertex set arising from the net points as
$$V_N \defeq \left\{(n, 2^{-k}) | n \in N, 2^{-k} \leq l(n)\right\}.$$

In particular, for every $n \in N$ for which $l(n)>0$ we will have infinitely many pairs $(n,2^{-k}) \in V_N$, and moreover if $(n,2^{-k}) \in V_N$, then also $(n,2^{-k'}) \in V_N$ for all $k' \geq k$. In a way that is made precise below, this gives a way of connecting $(n,2^{-k}) \in V_N$ to $n \in K$.

Define the edge set $E$ as follows. Two distinct vertices $(x,r), (y,s) \in V$ are connected\footnote{A minor technical point is that we use tuples. Thus, really the resulting space is a directed graph. We could also use undirected graphs, but it will simplify notation below to allow both edges. }, i.e. $e=((x,r), (y,s)) \in E$ if 
\begin{equation}
d(x,y) \leq 2^4(r+s), \label{edgecond1}
\end{equation}
and 
\begin{equation}
\frac{1}{2} \leq r/s\leq 2.  \label{edgecond2}
\end{equation}

We define a simplicial complex $T$ with edges corresponding to $e \in E$ and vertices corresponding to $v \in V$. In this geometric realization use $I_e$ to represent the interval corresponding to $e$. For each $v \in V$ we represent the corresponding point in $T$ with the same symbol. The simplicial complex becomes a metric space by declaring the length of the edge $e=((a,s),(b,t))$ and corresponding interval $I_e$ to be 
\begin{equation}
|e| \defeq 2^4(s+t). \label{edgelength}
\end{equation}
 A metric $d_T$ is induced by the path metric on the simplicial complex. Since $T$ may not be connected, some pairs of points $v,w \in T$ may have $d_T(v,w)=\infty$.

When a point $x$ lies on an interval $I_e$, we will often abuse notation and say $x \in e$.  Also, we will use the word edge to refer either to the symbol $e$ associated to it, or the interval $I_e$ in the geometric realization. Define the set which will be our metric space as $\nK \defeq K \cup T$. We define a symmetric function $\delta$ on $\mathcal{A} \defeq K \times K \cup \bigcup_{n,(n,s) \in V_N} \{(n,(n,s)),((n,s),n)\} \cup T \times T$, which is a subset of $\nK \times \nK$ by
\begin{itemize}
\item $\forall x,y \in K: \delta(x,y) \defeq d(x,y)$,
\item $\forall n \in N, (n,s) \in V_N : \delta(n,(n,s)) \defeq \delta((n,s),n) \defeq 3 \cdot 2^4 s$ and
\item $\forall v,w \in T : \delta(v,w) \defeq d_T(v,w)$.
\end{itemize}
Pairs $(a,b) \in \nK \times \nK$ are called admissible if $(a,b) \in \mathcal{A}$. Then we define $\nd$ to be the largest metric on $\nK$ which is dominated by $\delta$ for all admissible pairs.

The distance can be given explicitly by a minimization over discrete paths. Let $x,y \in \nK$ and $\sigma=(\sigma_0, \dots, \sigma_m)$ be a sequence of points in $\nK$. The variable $m$ is called the length of the path. We call such a sequence a discrete path in $\nK$, and say that it connects $x$ to $y$, if $\sigma_0=x, \sigma_m=y$. We call it admissible, if for each consecutive points $\sigma_{i}, \sigma_{i+1}$ the pair $(\sigma_i, \sigma_{i+1}) \in \mathcal{A}$. %, either they both lie in $T$ or $K$, or one of them is equal to $n$, for some $n \in N$, and the other is equal to $(n,s)$ for $(n,s) \in V_N$. In other words, we require $\delta(\sigma_i, \sigma_{i+1})$ to be defined by one of the cases above. 
Denote by $\Sigma^n_{x,y}$ the space of all admissible discrete paths of length $n$ that connect $x$ to $y$. With these definitions, the distance becomes
\begin{equation}
\nd(x,y) \defeq \inf_n \inf_{\sigma \in \Sigma_{x,y}^n} \sum_{i=0}^{n-1} \delta(\sigma_i, \sigma_{i+1}). \label{distdef}
\end{equation}

%Instead of using the path metric on $T$, we define $d_T$ to be the smallest metric on $T$ such that for two vertices $(x,r), (y,s) \in V$:

%$$d_T((x,r),(y,s))=d(x,y)+r+s,$$

%and for any two points $x,y$ on a common edge $I_e$ we use the Euclidean metric 

%\begin{equation}
%d_T(x,y)=|x-y| \label{edgedist1}.
%\end{equation}

%We remark, that the path metric is not desirable on $T$ since it may not even be connected. The previous conditions uniquely defines the metric $d_T$. For points on separate edges, $x\in e_1=((a_1,r_1), (b_1,s_1))$ and $y \in e_2=((a_2,r_2), (b_2,s_2))$ the distance can be computed as follows. Let $x_e$ be the corresponding co-ordinate on $I_{e_1}=[0,|e_1|]$, and $y_e$ the corresponding co-ordinate on $I_{e_2}=[0,|e_2|]$, where the endpoints have been identified with the corresponding vertices in $T$.

%\begin{eqnarray} d_T(x,y) = \min & \{ &  x_e + y_e + d(a_1,a_2)+r_1+r+2, \nonumber \\ 
%    && x_e + |e_2|-y_e + d(a_1, b_2)+r_1+s_2, \nonumber \\          
%    &&   |e_1|-x_e + y_e + d(b_1,a_2)+s_1+r_2, \nonumber \\
%    && \left. |e_1|-x_e + |e_2|-x_e + d(b_1, b_2)+s_1+s_2 \ \ \label{edgedist2}  \right\rbrace 
%    \end{eqnarray}

 A measure on each interval $I_{e}=[0,|e|]$, which is associated to the edge $e=((x,r),(y,s))$, is defined as weighted Lebesgue measure:
\begin{equation}
\mu_e \defeq \frac{\mu(B(x,r)) + \mu(B(y,s))}{|e|} \lambda, \label{measuredef}
\end{equation}
where $\lambda$ is the Lebesgue measure on $I_{e} \subset \R$. 

The total measure  is $\mu_T \defeq \sum_{e \in E} \mu_e$, where each measure $\mu_e$ is extended to be zero on the complement of $I_e$, and the measure $\nmu$ on the new space $\nK$ is defined as $\mu_T$ on $T$ and $\mu|_K$ on $K$. In other words
$$\nmu \defeq \mu_T + \mu|_K.$$

%The three conditions define a unique minimal metric satisfying the three conditions. The only case where the metric is not explicitly defined is when $x \in T$ is part of an edge $e=((a,r), (b,s))$, and $y \in K$, in which case we can define

%$$\nd(x,y) = \min\{d(x,(a,r))+d(a,y)+r, d(x,(b,s))+d(b,y)+s \}.$$

%The triangle inequality is trivial to verify. We also introduce the notation of distance of $x \in \nK$ to edge $e=((a,s),(b,t)) \in E$

%$$\nd(e,x) = \max\{\nd((a,s),x), \nd((b,t),x)\}.$$

\subsection{Proof that $(\nK, \nd,\nmu)$ is a PI-space}

Our goal is to show that the space constructed in the previous subsection satisfies the properties from Theorem \ref{relative}. This involves some preliminary lemmas.

\begin{definition} Let $L> 1$. A metric space $(X,d)$ is said to be uniformly $(L,r_0)$-perfect if for any $B(x,r)$ with $0<r<r_0$ the following holds
$$B(x,r_0) \setminus B(x,r) \neq \emptyset \Longrightarrow B(x,r) \setminus B(x,r/L) \neq \emptyset.$$

The space $X$ is said to be uniformly $(L,r_0)$-perfect along a subset $S \subset X$ if the same holds for any $x \in S$ and $0<r<r_0$.
\end{definition}

\begin{lemma} Let $r_0,C,\epsilon>0$ be arbitrary. If $A$ is a metric measure space and $A$ is $(C, \delta, \epsilon,r_0)$--connected along $K \subset A$ for some $0<\delta < \frac{1}{8}$, then $A$ is uniformly $(7/5,r_0)$-perfect along $K$.
\end{lemma}

\noindent \textbf{Proof:} Choose an arbitrary $x \in K$ and $0<r<r_0$. Assume $B(x,r_0) \setminus B(x,r) \neq \emptyset$. Let $s = \inf_{y \in A \setminus B(x, 5r/7)} d(x,y)$. By assumption $s<r_0$. If $s < r$, we are done. For the sake of contradiction, assume that $s \geq r$. Choose $y$ such that $d(x,y)< 8s/7$ and $d(x,y)<r_0$. By the definition of $s$ we know that $B(x,s) \setminus B(x,6s/7)$ is empty. Then, by the connectivity assumption (and choosing $E=\emptyset$), there is a Lipschitz curve  fragment $\gamma \co K' \to A$ connecting $x$ and $y$ with $\gap(\gamma)< \frac{1}{8}d(x,y) \leq \frac{s}{7}$. However,  since $B(x,s) \setminus B(x,6s/7)$ is empty, we can obtain a contradiction that $\gap(\gamma) \geq \frac{s}{7}$ by a similar argument as in Lemma \ref{thm:condoubl} and using the real-valued the curve fragment $\gamma' \co K \to \R$ given by $\gamma'(k)=d(\gamma(k),x)$.
\QED
%Assume $a = \sup_{t \in K', \gamma(t) \in B(x,6s/7)} t$ and $b = \inf_{a<t \in K', \gamma(t) \in A \setminus B(x,s)}$. Then 

%$$d(x,y)/8 < s/7\leq d(\gamma(a),\gamma(b)) \leq |b-a|.$$ 
%Further, the interval $(a,b) \subset \R \setminus K'$. Namely, if $t \in K' \cap (a,b)$, we have $\gamma(t) \in B(x,s)$ (by definition of $b$), but then automatically $\gamma(t) \in B(x,6s/7)$, which contradicts the definition of $a$. Then we conclude that $|b-a| \leq |\undf(\gamma)|< d(x,y)/8$, and this is a contradiction. 

This lemma is mainly used to give upper bounds for volumes in a doubling metric measure space.

\begin{lemma}\label{volume}Let $r_0,D>0$ be arbitrary. Assume that $X$ is a complete metric space, $X$ is $(D,r_0)$-doubling along $A \subset X$ and $A$ is $(C, \delta, \epsilon,r_0)$--connected along $K \subset A$ for some $0<\delta < \frac{1}{8}$. Then if $x \in K$, $0<r<r_0/2$ and $A \cap B(x,r_0)\setminus B(x,r) \neq \emptyset$,
\begin{equation}
\mu(B(x,r/2)) \leq (1-\sigma)\mu(B(x,r))
\end{equation}
for $\sigma = \frac{1}{D^4}$.
\end{lemma}

\noindent \textbf{Proof:} By the previous Lemma we get that $A$ is uniformly $(7/5,r_0)$-perfect along $K$. Then, using perfectness for the ball $B(x,7r/8)$ and the assumption there is a point $y \in A \cap B(x,7r/8) \setminus B(x,5r/8)$. For such a $y$ we have $B(y,r/8) \subset B(x,r) \setminus B(x,r/2)$. Also, by doubling 
$$\mu(B(y,r/8)) \geq \frac{\mu(B(x,r))}{D^4},$$
which gives
$$\mu(B(x,r/2)) \leq \mu(B(x,r)) - \mu(B(y,r/8)) \leq \left(1-\frac{1}{D^4} \right) \mu(B(x,r)).$$

\QED

We will also need a lemma concerning concatenating curve fragments.

\begin{lemma}\label{connectivity} Assume $(X,d,\mu)$ is a $(D,r_0)$-doubling metric space. If $1>\epsilon, \delta>0$, $L,C \geq 1$ are fixed and $p_i \in X$, for $i=1, \dots, n$, are points such that each pair $(p_{i}, p_{i+1})$ is  $(C,\delta, \epsilon)$--connected and $\sum_{i=1}^{n-1} d(p_i,p_{i+1})\leq Ld(p_1, p_n)$ , then $(p_1,p_n)$ is $$(L C,2L\delta, \epsilon D^{\log_2(\delta)-\log_2(L)-\log_2(n)-6})\text{--connected}.$$ 
\end{lemma}
\noindent \textbf{Proof:} Denote $r=d(p_1,p_n)$, and assume $E \subset B(p_1,LCr)$ with 
$$\mu(E \cap B(p_1,LCr)) \leq \epsilon D^{\log_2(\delta)-\log_2(L)-\log_2(n)-6} \mu(B(p_1,LCr)).$$

Define $l(1)=0$ and $l(j)=\sum_{i=1}^{j-1} d(p_i,p_{i+1})$ for $j=2, \dots, n$. Also, define the set of intervals $\mathcal{I}\defeq \{[l(i),l(i+1)] | i=1, \dots, n-1\}$ and $\mathcal{G} = \{I \in \mathcal{I} | |I| \geq \frac{\delta r}{n}\}$. For each $I=[l(i),l(i+1)] \in \mathcal{G}$ we have $d_i=d(p_i,p_{i+1})\geq \frac{\delta r}{n}.$ In particular, by doubling

$$\mu(E \cap B(p_i,Cd_i)) \leq \epsilon D^{\log_2(\delta)-\log_2(L)-\log_2(n)-6} \mu(B(p_1,LCr)) \leq \epsilon \mu(p_i, Cd_i).$$

Thus, we can define $\gamma_I \co K_I \to X$ to be a Lipschitz curve fragment connecting $p_i$ to $p_{i+1}$ with $K_I \subset I$, $\len(\gamma_I) \leq Cd_i$, $\gap(\gamma_I) < \delta d_i$ and $\gamma_I(I) \cap E \subset \{p_i, p_{i+1}\}$. Further, assume by a slight dilation and translation that $\min(K_I)=\min(I), \max(K_I)=\max(I)$. Define $K \defeq \{l(j)\} \cup \bigcup_{I \in \mathcal{G}} K_I$. Next, define a Lipschitz curve fragment $\gamma \co K \to X$ by setting $\gamma(l(i))=p_i$ and $\gamma(t) = \gamma_I(t)$ for $t \in I \in \mathcal{G}$. Clearly $\len(\gamma) \leq LC r$, and $\undf(\gamma) \subset \bigcup_{I \in \mathcal{I} \setminus \mathcal{G}} I \cup \bigcup_{I \in \mathcal{G}} \undf(\gamma_I)$. Thus,

\begin{align*}
 \gap(\gamma) &\leq \sum_{I=[l(i),l(i+1)] \in \mathcal{I} \setminus \mathcal{G}} d(p_{i+1},p_i) + \sum_{I=[l(i),l(i+1)] \in \mathcal{G}} \gap(\gamma_I) \\
 & <\delta d(x,y) + \delta \sum_{i=1}^{n-1} d_i \leq 2L\delta r.
\end{align*}

Since $p_i$ might belong to $E$, we might need to remove small neighborhoods of these points to satisfy the third condition in \ref{def:conbility}. This modification of $\gamma$ is trivial.

\QED 

Further, since our space $\nK$ arises by gluing a tree to $K$, the following lemma is useful.

\begin{lemma}\label{interval} Assume $(X,d,\mu)$ is a $(D,r_0)$-doubling metric space, and $I \subset X$ is isometric to a bounded interval in $\R$ of length at most $r_0/2$. Assume further that the restricted measure is given as $\mu|_I=c\lambda$ where $\lambda$ is the induced Lebesgue measure on $I$ and $c>0$ is some constant, and that for any sub-interval $[a,b] \in I$ we have $B((a+b)/2, (b-a)/2)\subset I$. Then for any $\delta>0$ and any $x,y \in I$ we have that $(x,y)$ is $(1,\delta, \delta (2D)^{-2})$--connected.
\end{lemma}

\noindent \textbf{Proof:} Let $x,y \in I$. Denote by $r=d(x,y)$ and by $z$ the midpoint of $x$ and $y$ on the interval $I$. Take an arbitrary Borel set $E$ with

$$\frac{\mu(B(x,r) \cap E)}{\mu(B(x,r))} \leq \delta / (2D)^2. $$

We will connect the pair of points $x,y$ by the geodesic segment $\gamma \co J \to X$ where $J \subset I$ is the sub-interval defined by $x$ and $y$. Note that $J=B(z,r/2)$. The curve $\gamma$ may intersect with $E$, but

\begin{eqnarray*}
|\gamma^{-1}(E)| &=&  \lambda(J \cap E) = \frac{1}{c} \mu(J \cap E) \\
                 & = & \lambda(J) \frac{\mu(J \cap E)}{\mu(J)} \leq  r \frac{\mu(E \cap B(x,r))}{\mu(B(z,r/2))} \\
                 & \leq & r D^2 \frac{\mu(E \cap B(x,r))}{\mu(B(x,r))} \leq r D^2 \delta / (2D)^2 = \delta r / 4.
\end{eqnarray*}
On the second line we used the fact that $J=B(z,r/2) \subset B(x,r)$.

Thus, by restricting $\gamma$ to an almost full measure subset in the complement of the set $\gamma^{-1}(E)$, we obtain the desired curve fragment.

\QED

Occasionally, we will need to vary the constant $C$ in Definition \ref{def:conbility}. Thus, we use the following.

\begin{lemma}\label{changeC} Assume $(X,d,\mu)$ is a $(D,r_0)$-doubling metric space and $0<\delta,\epsilon$, $C\geq 1$ and $K \in \N$ be given constants. If $(x,y)\in X \times X$ with $d(x,y) \leq r_02^{-K}/C $ is $(C,\delta,\epsilon)$--connected then it is also $(2^KC, \delta, \epsilon D^{-K-1})$--connected. 
\end{lemma}

\noindent \textbf{Proof:}  Denote by $r=d(x,y)$. Take an arbitrary Borel set $E$ with
$$\frac{\mu(B(x,2^K Cr) \cap E)}{\mu(B(x,2^K Cr))} \leq \epsilon D^{-K-1}.$$

Then, doubling also gives 
$$\frac{\mu(B(x,Cr) \cap E)}{\mu(B(x,Cr))} \leq \epsilon. $$

Thus, the result follows from the assumption and Definition \ref{def:conbility}.

\QED

\noindent \textbf{Proof that the construction in Subsection \ref{construction} works for Theorem \ref{relative}:} \\

\noindent \textbf{Preliminaries:} Continue using the notation of Subsection \ref{construction}. For future reference, we compute some basic estimates and define some related notation. Define the scale function  $\scl \co V \to \R$ which is given by $\scl((x,r))=r$  and location function $\loc \co V \cup K \to K$ which assigns for $(x,r) \in V$ the value $\loc((x,r))=x$ and restricts to identity on $K$. For a vertex $v=(x,r) \in V$ we denote the set $E_v=\{e  \in E | e \text{ is incident to } v\}$.  Throughout this section we will denote by $|S|$ the size of a set $S$, where the measure is evident from the context. If $S$ is finite, $|S|$ means the number of elements in that set. For edges and intervals $e$ we also use $|e|$ to denote their lengths.

\begin{lemma} \label{numest1}
For all $v \in V$ we have 
\begin{equation}
|E_v| \leq 4D^{25}. 
\end{equation}

\end{lemma} 

\noindent \textbf{Proof:} Denote $v=(x,r) \in V$. To see this consider $V_v^N \defeq \{w | w \in V_N, (v,w) \in E \}$ and $V_v^G \defeq \{w | w \in V_G, (v,w) \in E\}$. Then $|E_v| \leq 2|V_v^N| + 2|V_v^G|$. For $w \in V_v^G$, then $\scl(w) \geq r/2$, and $d(\loc(w),\loc(v)) \leq 2^6 r$. In particular, by the disjointedness property in Equation \eqref{whitneydisjoint}, the points $\loc(w)$ have pairwise distance at least a $r2^{-16}$ in $B(\loc(v),2^{6}r) \cap A$. Thus, their number is bounded by $D^{25}$, i.e. $|V_v^G| \leq D^{25}$. 

Similarly, if $w \in V_v^N$ we have $\scl(w) \geq r/2$. In particular, $l(\loc(w)) \geq r/2$. Thus, the pairwise distances of distinct $\loc(w)$ are at least $r/2$. Also, $$d(\loc(w),\loc(v)) \leq 2^6 r.$$ The number of such $w$ is bounded by $|V_v^N| \leq D^9$. We get by summing up the  estimates for $V_v^G$ and $V_v^N$ that
$$|E_v| \leq 4D^{25}.$$

\QED

\begin{lemma} \label{estimatesrel} The following estimates are true on $\nK$.

\begin{enumerate}
\item For every $x,y \in K \subset \nK$: $\nd(x,y) = d(x,y)$.
\item For every $e \in E$ and every $x,y \in I_e \subset T$ we have $\nd(x,y)=d_T(x,y)$.
\item For every $(x,r) \in V$ there is a $n \in N$ such that $e=((n,r),(x,r)) \in E$. In particular, $\nd((x,r),n) \leq 2^7r$.
\item If $(x,r) \in V$ and $y \in K$, then $\max\{r, d(x,y)\} \leq \nd((x,r),y) \leq d(x,y) + 2^8 r$.
\item If $(x,r), (y,s) \in V$ are distinct, then $\max\{d(x,y), r+s\} \leq \nd((x,r), (y,s)) \leq d(x,y) + 2^9(r+s)$.
\item If $x \in e=(v,w) \in E$ and $y \in \nK \setminus e$, then $\nd(x,y)=\min\{ \nd(x,v) + \nd(v,y), \nd(x,w) + \nd(w,y) \}$.
\item If $v \in V$, and $e \in E_v$, then $|e| \leq 2^7 \scl(v)$ and $\mu(B(\loc(v), \scl(v))) \leq \nmu(I_e) \leq 2D^7\mu(B(\loc(v), \scl(v)))$.
\item For any $n \in N$ and any $0<r\leq l(n)$ we have
$$\sum_{\substack{2^l \leq r}}\sum_{e \in E_{(n,2^l)}} \nmu(I_e) \leq 2^{4}D^{37} \mu(B(n,r)).$$
\item There is a function $\rho \co N \to G \cup N$ such that $\rho(n)=n$ if $l(n)=0$, and if $l(n)>0$, then $\rho(n) \in G_{l(n)}$ and  
$$d(\rho(n),n) \leq 2^4l(n).$$
Further, for any $g \in G$ we have $|\rho^{-1}(g)| \leq D^{8}$.
\end{enumerate}
\end{lemma}

\noindent \textbf{Proof:} We will proceed in numerical order. \\

\noindent \textbf{Estimate 1:} Take arbitrary $x,y \in K$. It is obvious that $\nd(x,y) \leq d(x,y)$. For an arbitrary $\epsilon> 0$ we can find a $n \in \N$ and $ \sigma \in \Sigma^n_{x,y}$ such that
$$\nd(x,y) + \epsilon \geq \sum_i \delta(\sigma_i, \sigma_{i+1}).$$

We can assume by possibly making the sum on the right smaller that for all $i=0,\dots,n$ it holds that $\sigma_i \in K$ or $\sigma_i \in V \subset T$. Thus, there is a function $f \co \{0, \dots, n \}\to K$, given by $f(i)=\sigma_i$ if $\sigma_i \in K$ and $f(i) = \loc(\sigma_i)$ if $\sigma_i \in V$. By the definition of $\delta$ preceding \eqref{distdef} we have for $i=0, \dots, n-1$
$$d(f(i), f(i+1)) \leq \delta(\sigma_i, \sigma_{i+1}).$$
By the triangle inequality
$$\nd(x,y) + \epsilon \geq \sum_i d(f(i), f(i+1)) \geq d(x,y).$$

Since $\epsilon>0$ is arbitrary $\nd(x,y) \geq d(x,y)$, which shows $\nd(x,y) = d(x,y)$. \\ \\

\noindent \textbf{Estimate 2:} The proof, that we have for any points $x,y \in T$ on a common edge $e \in E$ the equality $d_T(x,y) = \nd(x,y)$, is a trivial application of the definition of $\delta$ and observing that a discrete path that does not directly connect the points $x,y$ in $e$ will traverse an entire edge adjacent to both end points of $e$, or two separate edges adjacent to each end point of $e$. The length of such a path is larger than $|e| \geq d_T(x,y)$. \\ \\

\noindent \textbf{Estimate 3:} Next, take an arbitrary $(x,r) \in V$. We wish to find a near-by $n \in N$. Either $x \in N$ or $x \in G$. In the first case define $n \defeq x$, and by the definition of $\nd$ we have
$$\nd(n,(x,r)) \leq 3\cdot 2^4 r \leq 2^6 r.$$

Assume that $x=g \in G$. Then by  \eqref{whitneydist} there is a $k \in K$ such that $d(x,k) \leq r$. Let $n \in N$ be a net point closest to $k$ with $r \leq \scl(n)$. Then $d(n,k) \leq 2r$. Thus, by the triangle inequality $d(g,n) \leq 4r \leq 4 \scl(n)$, and by \eqref{lenscale} we have $l(n) \geq r$. Then both $(n,r)$ and $(g,r)$ are vertices in $V$ joined by an edge. Further, $\nd((n,r), (g,r)) = 2^5r$ by Estimate 2 above and \eqref{edgelength}. By the triangle inequality and what we just observed $\nd(n,(g,r)) \leq 2^7r$. \\ \\

\noindent \textbf{Estimate 4:} Take a $y \in K$ and $(x,r) \in V$. Choose $n \in N$ such that $\nd((x,r),n) \leq 2^7r$. We have $d(y,n) \leq d(x,n) + d(y,x) \leq 2^7r + d(y,x)$. Thus, $\nd((x,r), y ) \leq 2^8r + d(x,y)$. 

In order to obtain the desired lower bound, consider the height function given by $h \co \nK \to \R$ and defined as $h(x)\defeq 0$ for $x \in K$ and by $h(v)\defeq r$ for $v=(x,r) \in V$. On each edge extend $f$ linearly. Now for any $x,y \in \nK$ for which $\delta$ was defined (see discussion preceding \eqref{distdef}) we get $|h(x)-h(y)| \leq \delta(x,y)$. Thus, for any discrete admissible path $\sigma$ connecting $x$ to $y$ we obtain
$$|h(x)-h(y)| \leq \sum_i \delta(\sigma_i, \sigma_{i+1}).$$

Taking an infimum over $\sigma$ we get $|h(x)-h(y)|\leq \nd(x,y)$. A particular case of this gives when $(x,r) \in V$ and $y \in K$
$$r=|h((x,r))-h(y)| \leq \nd((x,r),y).$$

The lower bound $d(x,y) \leq \nd((x,r), y)$ is proven similarly, to the corresponding one in Estimate 5.
\\

\noindent \textbf{Estimate 5:} Let $(x,r), (y,s) \in V$ be distinct. Using Estimate 3 we can find $n_x, n_y \in K$
 such that $\nd((x,r), n_x) \leq 2^7 r$ and $\nd((y,s), n_y) \leq 2^7 s$. Further, by the calculations for Estimate 3 we get $d(n_x, x) \leq 2r$ and $d(n_y, y) \leq 2s$. A trivial application of the triangle inequality then gives $\nd((x,r), (y,s)) \leq 2^9(r+s) + d(x,y)$. 
 Next consider the lower bound. Let $\sigma$ be any discrete admissible path that passes through vertices in $V$ and points in $K$ and connects $(x,r)$ to $(y,s)$. The discrete path must traverse an edge adjacent to $(x,r)$, and one adjacent to $(y,s)$. These edges might be the same, or different, but in either case one can verify that their total length is at least $r+s$, and thus  $r+s \leq \nd((x,r), (y,s))$. Next, with the triangle inequality we obtain
$$d(x,y) \leq \sum_i d(\loc(\sigma_i), \loc(\sigma_{i+1})) \leq \sum_{i=0}^{n-1} \delta(\sigma_i, \sigma_{i+1}).$$

Taking the infimum gives $d(x,y) \leq \nd(x,y)$, which completes the proof. \\

\noindent \textbf{Estimate 6:} Trivial from the definition in \eqref{distdef}. \\ 

\noindent \textbf{Estimate 7:} These estimates follow directly from doubling and \eqref{edgecond1} and \eqref{edgecond2}. \\ 

\noindent \textbf{Estimate 8:} Take an arbitrary $n \in N$ and $r \leq l(n)$. Let $r \leq L=2^k \leq l(n)$ such that $L \leq 2r$. For any $s \leq l(n)$ and any $e \in E_{(n,s)}$ using Estimate 7 of this lemma we obtain
\begin{eqnarray}
\mu_e(I_e) &\leq& 2D^7\mu(B(n,s)). \label{edgest}
\end{eqnarray}

Apply Lemma \ref{numest1} to get
\begin{eqnarray*}
\sum_{\substack{2^l \leq r}} \sum_{e \in E_{(n,2^l)}} \nmu(I_e) &\leq& \sum_{\substack{2^l \leq L}} \sum_{e \in E_{(n,2^l)}} \nmu(I_e)  \\ 
 &\leq& \sum_{\substack{2^l \leq L}} 8D^{32} \mu(B(n,2^l)). 
\end{eqnarray*}

The last sum can be estimated using Lemma \ref{volume} to get
$$ \sum_{\substack{2^l \leq L}} 2^4 D^{32} \mu(B(n,2^l)) \leq  2^4 D^{36}\mu(B(n,L))\leq 2^4D^{37} \mu(B(n,r)).$$

Here we applied Lemma \ref{volume}, where the required non-emptyness of $B(n,r_0) \setminus B(n,s)$ can be guaranteed once $s \leq l(n)/2< r_0$. Note that $2^4 l(n)<r_0$, and $l(n)/2\leq d(n,g) \leq 2^{4}l(n)$ for some $g \in G$ by definition of $l(n)$ in Equation \eqref{lenscale}, thus $g \in B(n,r_0) \setminus B(n,s) \neq \emptyset$ for $s\leq l(n)/2$. \\ \\

\noindent \textbf{Estimate 9:} In the case $l(n)=0$, define $\rho(n)\defeq n$. By the definition of $l(n)$ in \eqref{lenscale}, for each $n \in N$ such that $l(n)>0$ we can find a $g \in G_{l(n)}$ such that $d(g,n) \leq 2^4 l(n)$. In this case, define $\rho(n)\defeq g$ by choosing one of them. This defines a function $\rho \co N \to G \cup N$. Next, if $g \in G$ we can use the definition of $l(n)$ to obtain two properties. Firstly, for any $n \in \rho^{-1}(g)$ we have the estimate $d(g,n) \leq 2^4 \scl(g)$. Secondly, if $n,m \in \rho^{-1}(g)$ are distinct then $d(n,m) \geq \scl(g)$ (since by Equation \eqref{lenscale} we have $n,m \in N$ and $\scl_N(n),\scl_N(m) \geq \scl(g)$). Thus, by doubling their total number can be bounded by $|\rho^{-1}(g)| \leq D^{8}$.

\QED

Next, we show completeness by taking an arbitrary Cauchy sequence $x_i$ and finding a limit point. If infinitely many of the $x_i$ lie in $K$, or a single edge $I_e$, the result follows trivially. Thus, we can assume by passing to a sub-sequence that each $x_i \in e_i =((a_i, r_i), (b_i, s_i))$ for distinct edges $e_i \in E$. We can assume $r_i,s_i \to 0$, since by the compactness of $K$ and the doubling of $A$ we obtain that for any given $\epsilon>0$ there are only finitely many edges with $r_i,s_i > \epsilon$. By estimate \eqref{whitneydist} we have $d(a_i, K)\leq r_i$, and thus by Estimate 4 in Lemma \ref{estimatesrel} there is a $k \in K$ such that 
$$\nd((a_i,r_i), k_i) \leq 2^{10} r_i.$$

On the other hand by \eqref{edgecond1} and the second and fifth estimate in Lemma \ref{estimatesrel} we get
$$\nd((a_i,r_i), x_i) \leq \nd((a_i,r_i), (b_i,s_i)) \leq d(a_i,b_i) + 2^9(r_i + s_i).$$

In particular, $\nd(x_i,k_i) \leq 2^9(r_i + s_i) + 2^{10}r_i.$ Thus, from the fact that $x_i$ is Cauchy, we can conclude that $k_i$ is also a Cauchy sequence in $K$ and has a limit point $k_\infty$. It is now easy to see that $\lim_{i \to \infty} x_i = k_\infty$. 

Our desired metric measure space is $(\nK,\nd,\nmu)$. The remainder of the proof shows that it is a PI-space.  We have already observed that $\nK$ is complete. Since $K$ and $T$ are totally bounded with respect to $\nd$ it is not hard to see that $\nK$ is also totally bounded and thus compact. The measure is a sum of finite Radon measures and thus as long as it's bounded it will itself be a Radon measure. This boundedness follows from estimates below. It remains to show the doubling and connectivity properties. The last step is to apply Theorem \ref{thm:contheorem} to conclude that the resulting space is a PI-space. 

In order to distinguish between metric and measure notions on $\nK$ and those on the original space $X$ we will often use a line above the symbol. For example for $x \in K$, the metric ball in $X$ is denoted $B(x,r)$, and the ball in $\nK$ is denoted $\nB(x,r)$. There is no risk of confusing this with the closure of the ball because we will not be using that in any of the arguments below. \\

\noindent \textbf{Doubling:} The doubling condition \eqref{eq:doubl} depends on two parameters $(x,r)$ and we need to check four cases depending on the location and scale. This is done by estimating the volume of various balls from above and below. \\

\noindent \textbf{Case 1: $x \in K$ and  $0<r<\frac{1}{2}$:}

\begin{itemize}
\item Lower bound for $x \in K, 0<r<1$:  We will use two disjoint subsets of $\nB(x,r)$: $\nB_K$ (part in $K$) and $\nB_T$ (corresponding to a portion of the glued intervals). In precise terms we set
$$\nB_K \defeq \nB(x,r) \cap K$$
and
$$\nB_T \defeq \bigcup_{\substack{v \in V_G, e \in E_v \\ \nd(v,x)<r2^{-10} }} I_e.$$

If $y \in \nB_T$, then there is a $v \in V$ such that $y \in e \in E_v$ and $\nd(v,x) < r2^{-10}$. By case 4 in Lemma \ref{estimatesrel}, we have $\scl(v) \leq \nd(v,x) <r2^{-10}$. Then, by Cases 3, 6 and 7 of Lemma \ref{estimatesrel} we have $\nd(x,y) \leq |e| + \nd(x,v) \leq r2^{-3} + r2^{-10} < r$.  Thus
\begin{equation}
\nB_K \cup \nB_{T} \subset \nB(x,r). \label{subset1}
\end{equation}

The sets $\nB_K$ and $\nB_T$ are disjoint, and we can estimate the total volume from below by estimating each of them individually and summing the estimates.
\begin{equation}
\nmu(\nB_K) = \mu(B(x,r) \cap K) \label{ineq1}
\end{equation}

The volume of $\nB_T$ can be estimated using Case 7 and Case 4 of Lemma \ref{estimatesrel}, and using the fact that each edge $e$ appears in $E_v$ for at most two different $v \in V$ .
\begin{eqnarray}
\nmu(\nB_T) &\geq& \frac{1}{2}\sum_{\substack{v=(g,s) \in V_G, e \in E_v \\ \nd(v,x)<r2^{-10} }} \mu_e(I_e)  \\
 &\geq& \frac{1}{2}\sum_{\substack{v=(g,s) \in V_G, e \in E_v \\ d(g,x)<r2^{-20} }} \mu(B(g,s))  \\
 &\geq&  \frac{1}{2} \mu(B(x, r2^{-26}) \cap A \setminus K). \label{ineq2}
\end{eqnarray}

In the last line we used the covering property
$$B(x, r2^{-26}) \cap A \setminus K \subset \bigcup_{\substack{ g\in G \\ d(x,g) \leq r2^{-20}}} B(g, \scl(g)2^{-10})).$$ 

To see this take an arbitrary $z \in B(x, r2^{-26}) \cap A \setminus K$. Then by the Whitney covering property \eqref{whitneycover} there is a $g \in G$ such that 
$$z \in B(g, \scl(g)2^{-10}).$$ However, $d(g,x) \geq \scl(g)2^{-1}$ by \eqref{whitneydist} and $d(g,x) \leq r2^{-26}+\scl(g)2^{-10}$ by the triangle inequality. Thus, we get $\scl(g) \leq r2^{-24}$ and $d(g,x) \leq r2^{-24}$. This is sufficient to verify the previous inclusion.

We obtain using estimates \eqref{ineq1} and \eqref{ineq2} that
\begin{eqnarray}
\nmu(\nB(x,r)) &\geq& \nmu(B_K) + \nmu(B_T) \nonumber \\
                      &\geq& \frac{1}{2}\mu(B(x, r2^{-26}) \cap A \setminus K)  + \mu(B(x,r2^{-26}) \cap A \cap K) \nonumber \\
                      &\geq&  \frac{1}{2}\mu(B(x,r2^{-26}) \cap A) \geq \frac{1}{D^{26} 4}\mu(B(x,r)). \label{case1low}
\end{eqnarray}
The uniform density of $A$ for $x \in K$ was used on the last line. \\

\item Upper bound for $x \in K, 0<r<1$: Define the following sets of edges:
$$E_2\defeq \{e \in E | \exists v=(n,s) \in V_N, e \in E_v, d(n,x) \leq r, s \leq l(n) \leq r\},$$
$$E_3\defeq \{e \in E | \exists v=(n,s) \in V_N, e \in E_v, d(n,x) \leq r, l(n) \geq r \geq s\},$$
$$E_4\defeq \{e \in E | \exists v=(g,s) \in V_G, e \in E_v, d(g,x) \leq r \}.$$

By the definition of the space and simple distance estimates from Cases 1, 4 and 6 in Lemma \ref{estimatesrel}, we can decompose the set into pieces as follows.
$$\nB(x,r) \subset A_1 \cup A_2 \cup A_3 \cup A_4,$$
with $A_1 \defeq B(x,r) \cap K=\nB(x,r) \cap K$ and 
$$A_i\defeq \bigcup_{e \in E_i} I_e, i=2,3,4.$$
% $$A_3 \defeq \bigcup_{e \in E_3} I_e$$
% and
% $$A_4 \defeq \bigcup_{e \in E_4}I_e.$$

First, we get the trivial upper bound 
\begin{equation}
\nmu(A_1) = \mu(B(x,r) \cap K) \leq \mu(B(x,r)). \label{a1up}
\end{equation}

Use Case 8 of Lemma \ref{estimatesrel} to obtain
\begin{eqnarray}
\nmu(A_2) &=& \sum_{e \in E_2} \nmu(I_e) \nonumber \\
          &\leq& \sum_{\substack{n \in N \\ d(n,x) \leq r \\ 0<l(n) \leq r}} \sum_{\substack{s \leq l(n)\\e \in E_{(n,s)}}} \nmu(I_e) \nonumber \\
          &\leq& \sum_{\substack{n \in N \\ d(n,x) \leq r \\ 0<l(n) \leq r}} 2^{4}D^{37} \mu(B(n,l(n))). \label{subsum1}
\end{eqnarray}

Next, we use the map $\rho$ from Case 9 of Lemma \ref{estimatesrel} and doubling to bound $\mu(B(n,l(n))) \leq D^{20} \mu(B(\rho(n), \scl_G(\rho(n))2^{-15}))$, and $d(x,\rho(n)) \leq 2^5 r$. Let $\overline{G} \defeq \rho\left(\{n \in N | d(n,x) \leq r, 0<l(n) \leq r\}\right)$. Then, we obtain the inclusion $B(g, \scl_G(g)2^{-15}) \subset B(x,2^6r)$ for $g \in \overline{G}$. Moreover, the balls $B(g, \scl_G(g)2^{-15}), g\in \overline{G}$  are disjoint for distinct $g \in \overline{G}$ by \eqref{whitneydisjoint}. All in all, we can apply a union bound with \eqref{subsum1}, and the bound $|\rho^{-1}(g)| \leq D^8$ from Case 9 in \ref{estimatesrel} to conclude
\begin{eqnarray}
\nmu(A_2) &\leq& \sum_{\substack{n \in N \\ d(n,x) \leq r \\ 0<l(n) \leq r}} 2^{4}D^{37} \mu(B(n,l(n))) \nonumber \\
          &\leq& 2^{4}D^{57} \sum_{\substack{n \in N \\ d(n,x) \leq r \\ 0<l(n) \leq r}}  \mu(B(\rho(n),\scl_G(\rho(n))2^{-15})) \nonumber \\
          &\leq& 2^{4}D^{57} \sum_{g \in \overline{G}} \sum_{n \in \rho^{-1}(g)}  \mu(B(\rho(n),\scl_G(\rho(n))2^{-15})) \nonumber \\
          &\leq& 2^{4}D^{65} \sum_{g \in \overline{G}}  \mu(B(g,\scl_G(g)2^{-15})) \nonumber \\
          &\leq& 2^{4}D^{65} \mu(B(x,2^6r)) \leq 2^{4}D^{71} \mu(B(x,r)).
\label{a2up}
\end{eqnarray}

Next, we estimate $\nmu(A_3)$. Use Case 8 of Lemma \ref{estimatesrel} to obtain
\begin{eqnarray}
\nmu(A_3) &=& \sum_{e \in E_3} \nmu(I_e) \nonumber \\
          &\leq& \sum_{\substack{n \in N \\ d(n,x) \leq r \\ r \leq l(n)}} \sum_{\substack{s \leq r \\ e \in E_{(n,s)}}} \nmu(I_e) \nonumber \\
          &\leq& 2^{4}D^{38} \sum_{\substack{n \in N \\ d(n,x) \leq r \\ r \leq \scl_N(n)}}  \mu(B(n,r/2)). \label{sumest2}
\end{eqnarray} 

Since $\scl_N(n)\geq r$, the net-points $n$ included in the sum are $r$-separated. Therefore, the balls $B(n,r/2)$ are disjoint. Also, $B(n,r/2) \subset B(x,2r)$. We obtain
\begin{eqnarray}
\nmu(A_3) &\leq& 2^{4}D^{38} \sum_{\substack{n \in N \\ d(n,x) \leq r \\ r \leq \scl_N(n)}}  \mu(B(n,r/2)) \nonumber \\
          &\leq& 2^{4}D^{38} \mu(B(x,2r)) \leq 2^{4}D^{39} \mu(B(x,r)). \label{a3up}
\end{eqnarray} 

Finally, we estimate $\nmu(A_4)$. Consider some $v=(g,r) \in V_G$ and $e \in E_v$ such that $d(g,x) \leq r$. By \eqref{whitneydist} we have $\scl(g)/2 \leq r$. Then by Case 7 $\nmu(I_e) \leq 2D^7 \mu(B(g,\scl(g))) \leq 2D^{25}\mu(B(g,\scl(g)2^{-15}))$. For different $v$ the balls $B(g,\scl(g)2^{-15})$ are disjoint by \eqref{whitneydisjoint} and $B(g,\scl(g)2^{-15}) \subset B(x,4r)$. Combine all these to see that 
\begin{eqnarray}
\nmu(A_4) &=&  \sum_{e \in E_3} \nmu(I_e) \nonumber \\
            &\leq& 2D^{25}\mu(B(x,4r)) \leq 2D^{30}\mu(B(x,r)). \label{a4up}
\end{eqnarray}

Combining the estimates \eqref{a1up}, \eqref{a2up}, \eqref{a3up}, \eqref{a4up}, with doubling for the different sets gives
\begin{equation}
\nmu(\nB(x,r)) \leq 2^{16} D^{71} \mu(B(x, r)). \label{case1up}
\end{equation}

\item Combine estimates \eqref{case1low} and \eqref{case1up} for $x \in K, 0<r<1/2$ to give
\begin{eqnarray}
\frac{\nmu(\nB(x,2r))}{\nmu(\nB(x,r))} &\leq& 2^{20} D^{100}. \label{case1doubl}
\end{eqnarray}
\end{itemize}

\  \\

\noindent \textbf{Case 2: $x \in e_x=((a,s),(b,t)) \in E$, $0<r<s/4$:} Denote $v=(a,s)$ and $w=(b,t)$. Let $E_v,E_w$ be the sets of edges adjacent to either vertex. It is easy to conclude that 
\begin{equation} 
\nB(x,r) \subset \nB(x,2r) \subset \bigcup_{e \in E_v \cup E_w} I_e \label{subset2}.
\end{equation}

For each $e \in E_v \cup E_w$ and any $h>0$ we have $\text{diam}(B(x,h) \cap e) \leq 2h.$ Also, $|e| \geq s$ by \eqref{edgelength} and \eqref{edgecond2}. Further, we derive from \eqref{measuredef} that in terms of densities
$$\mu_e \leq \frac{2\mu(B(a, 2^{10}s))}{|e|} \lambda \leq \frac{2\mu(B(a, 2^{10}s))}{s} \lambda.$$ 

Thus using $h=2r$, the size bound Lemma \ref{numest1} we get
\begin{eqnarray} 
\nmu(\nB(x,2r)) &=& \sum_{e \in E_v \cup E_w} \mu(I_e \cap B(x,2r)) \nonumber \\
  &\leq& \frac{8r}{s} \mu(B(a, 2^{10}s)) (|E_v|+|E_w|) \nonumber \\
  & \leq&\frac{2^{10} D^{35} r \mu(B(a, s))}{s}  \label{c2up}.
\end{eqnarray}

On the other hand for $e_x$ we have $\text{diam}(B(x,r) \cap e_x) \geq r,$ and using $|e| \leq 2^7s$ from Case 7 in Lemma \ref{estimatesrel} we see that in terms of densities
$$\mu_e \geq \frac{\mu(B(a, s))}{|e|} \lambda \geq \frac{\mu(B(a, s))}{2^7s} \lambda.$$ 

Using this we obtain the lower bound
\begin{eqnarray} 
\nmu(\nB(x,r)) &\geq& \mu_e(e \cap B(x,r)) \nonumber \\
               &\geq& \frac{r\mu(B(a, s))}{2^7s}  \label{c2low}.
\end{eqnarray}

Combine estimates  \eqref{c2low} and \eqref{c2up} to get the desired doubling bound
$$\frac{\nmu(\nB(x,2r))}{\nmu(\nB(x,r))} \leq 2^{20}D^{35}.$$ \\ 

\noindent \textbf{Case 3: $x \in e=((a,s),(b,t)) \in E$, $s/4<r<2^{12} s$ and $r < 2^{-13}$:} Let $v=(a,s)$. By the estimate \eqref{c2low} and using that for $r2^{-14}<s/4$, we get
\begin{equation}
\nmu(\nB(x,r)) \geq \nmu(\nB(x,r2^{-14})) \geq \frac{r\mu(B(a, s))}{2^{30}r} \geq 2^{-30}\mu(B(a,r2^{-12}))  \label{c3low}. 
\end{equation}

By Case 3 of Lemma \ref{estimatesrel} there is a $n \in N$ such that $\nd(v, n) \leq 2^7 s$. In particular, $\nd(n,x) \leq \nd(v,n) + |e| \leq 2^7s + 2^7 s \leq 2^{8}{s}.$ Apply the estimate \eqref{case1up} to see
\begin{equation}
\nmu(\nB(x,2r)) \leq \nmu(\nB(n,2^{11}r)) \leq 2^{16} D^{71} \mu(B(n, 2^{11}r)).\label{c3up}
\end{equation}

Using doubling with Cases 4 in Lemma \ref{estimatesrel} we get $d(n,a) \leq \nd(n,v) \leq 2^{10}r$  and $\mu(B(n,2^{11}r)) \leq \mu(B(a,2^{12}r)) \leq D^{25}\mu(B(a,2^{-12}r))$. Finally, combining this with  \eqref{c3low} we get the desired doubling bound
$$\frac{\nmu(\nB(x,2r))}{\nmu(\nB(x,r))} \leq 2^{46}D^{96}.$$ \\

\noindent \textbf{Case 4: $x \in e=((a,s),(b,t)) \in E$, $2^{12}s < r$ and $r<1/8$:} Using Cases 3 and 7 of Lemma \ref{estimatesrel} we can find a $n$
$$\nd(x,n) \leq 2^8 s \leq r2^{-4}.$$

In particular, we have
$$\nB(n,r2^{-4}) \subset \nB(x,r) \subset \nB(x,2r) \subset \nB(n, 4r).$$

Thus using \eqref{case1doubl}, we obtain
$$\frac{ \nB(x,2r) }{ \nB(x,r)} \leq \frac{\nmu(\nB(n, 4r))}{\nmu(\nB(n, r2^{-4}))} \leq 2^{200} D^{500}.$$
~\\ ~\\

These cases conclude the proof for doubling. Denote the maximum of the previous doubling constants $\nD=2^{200}D^{500}$. The cases above show that $(\nK,\nd,\nmu)$ is $(\nD,2^{-13})$-doubling, which becomes $(\nD, 2^{-33}r_0)$ with the scaling of $r_0=2^{20}$. \\ \\

\noindent \textbf{Connectivity:} Next, we move to prove connectivity. We will show that every pair of points $(x,y)$, which is sufficiently close, is $(\nC, \frac{1}{2}, \neps)$--connected for appropriately chosen $\neps$ and $\nC$. Again, we have cases depending on the positions of $x$ and $y$. The most complicated case is when both $x,y \in K$.

We always denote $r=\nd(x,y)$. For each pair of points take an arbitrary open set $E$ such that $\nmu(E \cap \nB(x,C r)) \leq \eta \nmu(\nB(x,C r))$, where $\eta,C$ are the relevant connectivity parameters in each case. Also, it is enough to verify Definition \ref{def:conbility} for open sets, because of the regularity of the measure. Here, we have already shown that $\nmu$ is a bounded measure on a compact metric space, and thus $\nmu$ is a Radon measure, and both inner and outer regular (see \cite{Mattila1999}). Note, by doubling and Lemma \ref{changeC}, we don't need to be so concerned about the parameter $C$, which represents the length of the curve fragments, being the same in each case. At the end $\nC$ will become the maximum from the different $C$'s in each case, and we will adjust $\neps$ to account for this.

For technical reasons, some cases involve seemingly stronger statements than necessary. Again, we will be somewhat liberal in our estimates. For example, we will use $\nD$ in places where $D$ would suffice. Together these cases show that each $(x,y) \in \nX \times \nX$ with $\nd(x,y) \leq 2^{-200}C^{-1}$ is $$\left(2^{40}C, 2^{-14},\epsilon (2\nD)^{-2000-2\log_2(C)}\right)\text{--connected}.$$

\noindent \textbf{Case 1: Let $\delta>0$ be arbitrary and assume $(x, y) \in T \times T$ with $\nd(x,y)<2^{-6}$, and $x \in e=((a,s),(b,t)) \in E$. If $\nd(x,y)<s2^{-4}$, then the pair $(x,y)$ is $(1,\delta, \delta (2\nD)^{-11+\log_2(\delta)})$--connected.} \\

Let $v=(a,s)$ and $w=(b,t)$. By Case 2 in \ref{estimatesrel} and \eqref{edgelength} we know that $x,y \in S \defeq \cup_{e \in E_v \cup E_w} I_e$. Lemma \ref{numest1} gives that $S$ is a connected subgraph of $T$ with at most $4D^{13}$ edges. The metric restricted to $S$ agrees with the path metric $d_T$ by similar arguments to those used in Estimate 2 in Lemma \ref{estimatesrel}. The path connecting any pair of points passes through at most three intervals. The result follows now easily from Lemmas \ref{interval} and \ref{connectivity}. \\ \\

\noindent \textbf{Case 2: Every pair $(x,y) \in K \times K$ with $d(x,y) \leq 2^{-200}C^{-1}$ is $$\left(2^{25}C, 2^{-30},\epsilon\left(2\nD{}\right)^{-1000-\log_2(C)}\right)\text{--connected}.$$} Let $E \subset \nB(x,2^{20}Cd(x,y))$ be a set such that
\begin{equation}
\nmu( E) \leq \epsilon \left(2\nD{}\right)^{-1000-\log_2(C)} \nmu\left(\nB\left(x,2^{25}Cd(x,y)\right)\right). \label{densdef}
\end{equation}

Our goal will be to first construct a curve fragment in $A$ and then to replace portions of it with a curve in $T$. This latter step is only possible, if the portions of the curve fragment in $A$ doesn't pass too close to certain ``bad'' gap points. First, we define bad bridge points.

Consider the collection $\mathbf{B}$ of ``bad'' bridge points $b \in V$ with $\scl(b)<2^{-50}$ such that
\begin{equation}
\nmu(\nB(b, 2^{10}\scl(b)) \cap E) \geq \left(2\nD\right)^{-\log_2(C)-500} \nmu\left(\nB\left(b, 2^{10}\scl(b)\right)\right). \label{badball}
\end{equation}

Define
\begin{equation}\label{eq:mcBdef}
 \mathcal{B} \defeq \bigcup_{b \in \mathbf{B}} B(\loc(b), \scl(b))
\end{equation}
and its approximant in $\nK$
$$\overline{\mathcal{B}} = \bigcup_{b \in \mathbf{B}, e \in E_b} I_e.$$

We will seek to estimate the volume of $\mathcal{B}$ from above. Each of the edges $I_e$ can appear at most twice in the union defining $\overline{\mathcal{B}}$. Use Case 7 of Lemma \ref{estimatesrel} to see that
\begin{eqnarray}
\mu(\mathcal{B}) &\leq& \sum_{b \in \mathbf{B}} \mu(B(\loc(b), \scl(b))) \nonumber \\
&\leq&  \sum_{b \in \mathbf{B}, e \in E_b} \nmu(I_e) \nonumber \\
&\leq& 2\nmu(\overline{\mathcal{B}}). \label{Best}
\end{eqnarray}

We can use the Vitali covering Lemma \cite{stein2016harmonic,federer} to choose a sub-collection $\mathbf{B}'$ of $b \in \mathbf{B}$ such that

$$\overline{\mathcal{B}} \subset \bigcup_{b \in \mathbf{B}'} \nB(b, 2^{13}\scl(b)) $$
and the balls $\nB(b, 2^{10}\scl(b))$ for $b \in \mathbf{B}'$ are disjoint. Apply doubling of $\nmu$, which we have already found, and the definition \eqref{badball} to get the following.
\begin{eqnarray}
\nmu(E) &\geq&   \sum_{b \in \mathbf{B}'} \nmu( \nB(b, 2^{10}\scl(b)) \cap E) \nonumber \\
           &\geq&       \left(2\nD\right)^{-\log_2(C)-500} \sum_{b \in \mathbf{B}'} \nmu\left( \nB(b, 2^{10}\scl(b)) \right) \nonumber \\
           &\geq&       \left(2\nD\right)^{-\log_2(C)-503} \sum_{b \in \mathbf{B}'}\nmu\left( \nB(b, 2^{13}\scl(b))\right) \nonumber \\
           &\geq&     \left( 2\nD \right)^{-\log_2(C)-503} \nmu(\overline{\mathcal{B}}) \label{volumeest}
\end{eqnarray}
Then, using \eqref{densdef} with doubling we see
\begin{equation}
\nmu(\overline{\mathcal{B}}) \leq  \epsilon\left(2\nD{}\right)^{-400} \nmu\left(\nB(x,Cd(x,y))\right). 
\end{equation}
Thus, with the estimates \eqref{densdef}, \eqref{Best} and \eqref{case1up}, we see that
\begin{equation}
\mu(\mathcal{B} \cup (K\cap E)) \leq 2\nmu(\overline{\mathcal{B}}) + \nmu(E) \leq \epsilon \mu\left(B(x,Cd(x,y))\right). \label{estimatekm}
\end{equation}

By using the connectivity condition, we can find a $1$-Lipschitz curve fragment $\gamma\co S \to A$ defined on a compact subset $S$ with $\min(S)=0,\max(S)=\len(\gamma)$, which is parametrized by length (see Lemma \ref{adjust}) and which satisfies $\gamma(0)=x, \gamma(\len(\gamma))=y$, 
\begin{equation}\label{eq:gamgapest}
|\undf(\gamma)|=\gap(\gamma) < 2^{-60} d(x,y),
\end{equation}
 $\len(\gamma)\leq Cd(x,y)$ and 
\begin{equation}\label{eq:avoidanceprop}
 \gamma^{-1}(\mathcal{B} \cup (K \cap E))) \subset \{0, \max(S)\}.
\end{equation}

The image of the curve fragment $\gamma$ may not be contained in $K$. Thus, to reach our desired conclusion we will define another curve fragment $\ngam$ in $\nK$ by replacing the portions in $A/K$ with curves in $\nK$. This is done in two steps. First we discretize the portions of the fragment in $A \setminus K$ and obtain a ``discrete path fragment'' through gap points $g$. Intuitively, the sub-segments of the curve fragment associated to gap points $g$ will be related to the corresponding bridge point $b \in V$. In the construction, we can guarantee that $b \not\in \mathbf{B}$, and thus we will be able to connect sufficiently many consecutive bridge points to each other and go from a discrete path to a continuous curve fragment. This latter part of the argument is called extension. %Finally, the curve fragment we construct won't satisfy the desired Lipschitz-bound so we will need to dilate the domain by a certain factor and thus complete the proof.  \\

\begin{enumerate}

\item \textbf{Discretization:} First, we cover the set $\gamma^{-1}(A \setminus K)=O$ by intervals as follows. Recall the Whitney covering property \eqref{whitneycover}, the length estimate, and the fact that $\gamma(O) \cap \mathcal{B} = \emptyset$. Using these, choose for each $z \in O$ a $g_z \in 
G$ with the property 
\begin{equation}\label{eq:gzdef}
 d(g_z, \gamma(z)) \leq \scl_G(g_z)2^{-10}.
\end{equation}
 Define $b_z\defeq (g_z, \scl_G(g_z))$. It is easy to see that $b_z \not\in \mathbf{B}$, since otherwise $\gamma(z) \in B(g_z, \scl_G(g_z))$, which contradicts the avoidance property \eqref{eq:avoidanceprop} and the definition of $\mathcal{B}$ in \eqref{eq:mcBdef}. 

Define the interval 
\begin{equation}
I_z\defeq (z-\scl_G(g_z)2^{-10}, z + \scl_G(g_z)2^{-10}) \label{interdef}
\end{equation} 
and the smaller interval $J_z \defeq (z-\scl_G(g_z)2^{-15}, z + \scl_G(g_z)2^{-15})$ and choose discrete centers $\nZ \subset O$ such that 

$$O \subset \bigcup_{z \in \nZ} I_{z},$$
and $J_{z}$ are pairwise disjoint for $z \in \nZ$. This is possible by the Vitali covering theorem \cite{federer, stein2016harmonic}. Because $\gamma$ is 1-Lipschitz and from the distance inequality \eqref{whitneydist}, we obtain $I_{z} \subset [0,\max(S)] \setminus \gamma^{-1}(K)$ for all $z \in \nZ$.

We will define a ``discrete curve fragment'' as follows. Consider the compact set $\nS_0 \defeq \gamma^{-1}(K) \cup \nZ$ and the curve fragment $\ngam_0 \co \nS_0 \to \nK$ which is defined by setting $\ngam_0(t)=\gamma(t)$ for $t \in S$ and $\ngam_0(z)=b_z$ for $z \in \nZ$. The crucial property we need is that $\ngam_0$ is Lipschitz with $\LIP(\ngam_0) \leq 2^{25}.$ We proceed next to showing this.

%The first estimate is immediate from the definition. 
Since the curve fragment $\gamma$ is $1$-Lipschitz on $\gamma^{-1}(K)$ (this involves Estimate 1 from Lemma \ref{estimatesrel}), it is sufficient to check the Lipschitz-bound for pairs of points $z,w \in \nZ$ and $z \in \nZ$ and $w \in \gamma^{-1}(K)$.

Consider the  case where both $z,w \in \nZ$. Then by the triangle inequality
$$d(g_z, g_w) \leq d(g_z,\gamma(z)) + d(\gamma(z), \gamma(w)) + d(g_w, \gamma(w)).$$

Further by the Lipschitz bound $d(\gamma(z), \gamma(w)) \leq |z-w|$, and by the choice of $g_z$ and $g_w$, for $t=z,w$ we get
$d(g_t, \gamma(t)) \leq \scl_G(g_t)2^{-10}$. Also, $J_z \cap J_w = \emptyset$, so $\scl_G(g_z)2^{-15} + \scl_G(g_w)2^{-15} \leq |z-w|$. All these combined give $d(g_z, g_w) \leq 2^7|z-w|.$

Use Estimate 5 from Lemma \ref{estimatesrel} to give the Lipschitz-bound.
\begin{eqnarray*}
\nd(\ngam_0(z),\ngam_0(w))&=&\nd(b_z,b_w) \\ 
&\leq& d(g_z, g_w) + 2^9(\scl_G(g_z) + \scl_G(g_w)) \\
&\leq& 2^{25}|z-w|.
\end{eqnarray*}

Consider next the case $z \in \nZ$ and $w \in \gamma^{-1}(K)$. Then by the triangle inequality
$$d(g_z, \gamma(w)) \leq d(g_z,\gamma(z)) + d(\gamma(z), \gamma(w)).$$

Further by the Lipschitz bound $d(\gamma(z), \gamma(w)) \leq |z-w|$, and by the choice of $g_z$ and Estimate \eqref{whitneydist} we get
$d(g_z, \gamma(z)) \leq \scl_G(g_z)2^{-10} \leq d(g_z, K) 2^{-9} \leq d(g_z, \gamma(w)) 2^{-9}$. Thus
$$d(g_z, \gamma(w)) \leq 2 |z-w|.$$

Also, $\scl_G(g_z) \leq 2 d(g_z, K) \leq 2d(g_z, \gamma(w)) \leq 4|z-w|$.
Use Estimate 4 from Lemma \ref{estimatesrel} for $b_z=(g_z, \scl_G(g_z))$ to give the Lipschitz-bound
\[
\nd(\ngam_0(z),\ngam_0(w))=\nd(b_z,\gamma(w)) \leq d(g_z, \gamma(w)) + 2^8\scl_G(g_z) \leq 2^{25}|z-w|.\]
~\\

\item \textbf{Extension:} The complement $[0,\max(\nS_0)] \setminus \nS_0$ is open and can be expressed as a union of maximal open intervals $I_{z,w}=(z,w)$, where $z,w \in \nS_0$. Let $\mathcal{I}\defeq\{I_{z,w}\}$ be the collection of all such intervals. Then, let $\mathcal{G}$ be the  collection of ``good'' intervals $I_{z,w}$ such that $z,w \in \nZ$ and 
\begin{equation}
d(g_z, g_w) \leq 2^{-3}(\scl_G(g_z) + \scl_G(g_w)), \label{goodinter}
\end{equation}
%or at least one of the points $z,w \in\gamma^{-1}(K)$. 
Also, let $\mathcal{F} \defeq \mathcal{I} \setminus \mathcal{G}$. 

Define $\nS \defeq \nS_0 \cup \bigcup_{I \in \mathcal{G}} I$, and the curve fragment $\ngam \co \nS \to \nK$  by $\ngam(t) = \ngam_0(t)$ for $t \in \nS_0$, and $\ngam$ is a linear parametrization of the edge $I_{(b_z,b_w)} \subset T$ for $I_{z,w} \in \mathcal{G}$. This edge exists by the following argument. 

The edge $(b_z,b_w)$ exists if the edge conditions \eqref{edgecond1} and \eqref{edgecond2} are satisfied. Firstly, from \eqref{whitneydist} we get $\scl(g_z)/2 < d(g_z, K) \leq d(g_z, g_w) + d(g_w, K)$. Thus, $\scl(g_z) 3/8 < \frac{9}{8} \scl(g_w)$ follows from estimates \eqref{goodinter} and \eqref{whitneydist}. By symmetry, we obtain $|\scl(g_z)/\scl(g_w)|<4$, and thus \eqref{edgecond2} must hold. Note, $\scl(g_z)/\scl(g_w)$ is always a power of $2$, and here we got a strict inequality. On the other hand, the condition \eqref{edgecond1} is immediate from \eqref{goodinter}. \\

\item \textbf{Estimates:} We will show that this curve fragment satisfies the following.
\begin{enumerate}
\item $\LIP(\ngam) \leq 2^{25}$
\item $\len(\ngam) \leq 2^{25} C d(x,y)$
\item $\gap(\ngam) < 2^{-34}d(x,y)$
\item $\int_{\ngam} 1_E\, ds < 2^{-300}d(x,y)$ 
\end{enumerate}

From these, the desired curve fragment satisfying the conditions of Definition \ref{def:conbility} can be obtained by restricting $\ngam$ to a large compact set in the complement of $\ngam^{-1}(E)$. Thus, we are left with proving these four estimates. 

The Lipschitz estimate is trivial, since the domain is expanded by linear parametrizations. Here we also use Estimate 2 in Lemma \ref{estimatesrel}. Similarly, the length estimate follows because $\len(\ngam) \leq 2^{25} \max(\nS_0)  \leq 2^{25}\len(\gamma) \leq 2^{25} C d(x,y)$. 

Next, estimate $\gap(\ngam)$. It is easy to see $\undf(\ngam) \subset \bigcup_{I \in \mathcal{F}} I.$ Let $I=I_{z,w} \in \mathcal{F}$ be an arbitrary interval. If $z,w \in \nZ$, define $C_{I}\defeq (z,w) \setminus (I_z \cup I_w)$. If  $z,w \not\in \nZ$, then define $C_I\defeq (z,w)$. If only $z \in \nZ$, then define $C_{I}\defeq (z,w) \setminus I_z$, and if only $w \in \nZ$, then set $C_{I}\defeq (z,w) \setminus I_w$. Clearly $C_I \subset \undf(\gamma)$, because $C_I \cap \gamma^{-1}(K) = \emptyset$ and the intervals $I_z$ cover $\gamma^{-1}(A\setminus K)$. We will show in the following paragraphs that $|I| \leq 2|C_{I}|$. Assuming this for now, we will complete the estimate. Since the intervals $I \in \mathcal{F}$ are disjoint, so are $C_I \subset I$. Also,
\begin{eqnarray*}
\gap(\ngam) &\leq& \LIP(\ngam) |\undf(\ngam)| \\
&\leq& 2^{25} \sum_{I \in \mathcal{F}} |I| \\
&\leq& 2^{26}\sum_{I \in \mathcal{F}} |C_I| \\ 
&\leq& 2^{26} |\undf(\gamma)| \leq 2^{-34} d(x,y).
\end{eqnarray*}
On the last line we used estimate \eqref{eq:gamgapest}. Now, verify $$|I| \leq 2|C_{I}|$$ for $I \in \mathcal{F}$. There are three cases:$z,w \not\in \nZ$, one of them is not in $\nZ$, or both are in $\nZ$. In the first case $I=C_I$ and thus clearly $|I| \leq 2|C_I|$. In the second case, suppose by symmetry that $z \in \gamma^{-1}(K), w \in \nZ$. Then $C_I = (z,w) \setminus I_w = (z,w-\scl_G(w)2^{-10})$, and a simple computation using the property \eqref{whitneydisjoint}, Lipschitz property of $\gamma$, the definition of $g_w$ in \eqref{eq:gzdef} and that $\gamma(z) \in K$ gives $|w-z| \geq \scl_G(g_w)/4$. The details are similar to the proof above of the Lipschitz property of $\ngam_0$. Finally, $|C_I| = w-\scl_G(w)2^{-10}-z \geq |w-z|/2 = |I|/2$. 

Now, for the third case. For each $I=I_{z,w} \in \mathcal{F}$ with $z,w \in \nZ$ we have 
\begin{equation}\label{eq:Fdist}
 d(g_z, g_w) > 2^{-3}(\scl_G(g_z) + \scl_G(g_w))
\end{equation}
and $(z,w) \cap \nS = \emptyset$.  Thus, from the triangle inequality, Lipschitz bound of $\gamma$, estimate \eqref{eq:Fdist} and the definition of $g_z,g_w$ in \eqref{eq:gzdef} we have $|z-w| \geq d(\gamma(z), \gamma(w)) \geq 2^{-4}(\scl_G(g_z) + \scl_G(g_w))$. Again, for similar arguments see the proof of the Lipschitz property of $\ngam_0$. Also, from definition \eqref{interdef} and the previous computations, we get
\begin{eqnarray*}
|C_{I}|&=&|(z,w) \setminus ( I_z \cup I_w)| \\
      &\geq& 2^{-4}(\scl_G(g_z) + \scl_G(g_w)) - 2^{-10}(\scl_G(g_z) + \scl_G(g_w)) \\
      &\geq& 2^{-5} (\scl_G(g_z) + \scl_G(g_w)).
\end{eqnarray*}

Finally, from \eqref{interdef}, 
$$|I| \leq |C_{I}| + 2^{-10}\scl_G(g_z) + 2^{-10}\scl_G(g_w) \leq 2|C_{I}|.$$

The remaining estimate concerns $\int_{\ngam}1_E\, ds$. For each $I=(z,w) \in \mathcal{G}$ denote by $e_I=(b_z, b_w)$ the edge in $T \subset \nK$ corresponding to it. Note that 
\begin{equation}
|e_I| \leq 2^{25}|z-w|=2^{25}|I| \label{Eest}
\end{equation}
by the Lipschitz estimate of $\ngam$. Assume without loss of generality that $\scl(b_w) \leq \scl(b_z)$. Since $b_z, b_w \not \in \mathcal{B},$ we have by Estimate \eqref{badball}
\begin{equation}
\nmu(\nB(b_z, 2^{10}\scl(b_z)) \cap E) < \left(2\nD\right)^{-\log_2(C)-500} \nmu\left(\nB\left(b_z, 2^{10}\scl(b_z)\right)\right). 
\end{equation}

Thus, by doubling estimates similar to Lemma \ref{interval}, estimate $\scl(b_w) \leq \scl(b_z)$ and estimate \eqref{edgelength} we have
\begin{equation}
\nmu(e_I \cap E) \leq \left(2\nD\right)^{-\log_2(C)-400} \nmu(e_I). \label{densest}
\end{equation}

Since $\ngam$ is a linear parametrization of $e_I$ with the interval $I$, we have also
\begin{equation}
\int_{\ngam|_I} 1_E \, ds \leq \left(2\nD\right)^{-\log_2(C)-300} |I|. \label{densest2}
\end{equation}

Thus, we get
\begin{eqnarray*}
\int_{\ngam}1_E \, ds &\leq& \sum_{I \in \mathcal{G}} \int_{\ngam|_I} 1_E \, ds \\
                &\leq& \sum_{I \in \mathcal{G}} \left(2\nD\right)^{-\log_2(C)-300} |I| \\
                &\leq& \left(2\nD\right)^{-\log_2(C)-300}\max(\nS_0) \\
                &\leq& \left(2\nD\right)^{-\log_2(C)-300} Cd(x,y) \leq 2^{-59} d(x,y).
\end{eqnarray*}

\end{enumerate}

\noindent \textbf{Case 3.  For any $\delta>0$ and every $x\in e = ((a,s), (b,t)) \in T$ with $s<2^{-6}$ there is a $y \in K$ such that $\nd(x,y) \leq 2^8 s$ and such that $(x,y)$ is $(2^9,\delta, \delta (2\nD)^{2\log_2(\delta)-80})$--connected.} Let $v=(a,s),w=(b,t)$. By Cases 3 and 4 in Lemma \ref{estimatesrel} there is a $n \in N$ such that $s \leq \nd(n,v) \leq 2^7 s$. Then by Case 7 and 6 in Lemma \ref{estimatesrel} $\nd(v,x) \leq 2^8s$. Further, $((n,s), v) \in E$. Set $y=n$.

Next, take $k = \lceil 20-\log_2(\delta) \rceil$. Define the points $p_1\defeq x, p_2\defeq v, p_{k+1}\defeq n=y$ and $p_i \defeq (n,s2^{2-i})$ for $i=3, \dots, k$. By Cases 4 and 6 or Lemma \ref{estimatesrel} we have
$$\nd(p_{k}, y) \leq 2^8 \scl(p_k) \leq 2^{\log_2(\delta)-5} s\leq \frac{\delta}{4} \nd(x,y).$$

Each pair $(p_i,p_{i+1})$ are edges in $T$ by condition \eqref{edgecond1} and \eqref{edgecond2}. Thus, by Lemma \ref{interval}, each pair $(p_i, p_{i+1})$, for $1 \leq i \leq k-1$, is $(1,\delta 2^{-11}, \delta (2\nD)^{-14})$--connected. Also,
$$\sum_{i=1}^{k-1} \nd(p_i, p_{i+1}) \leq 2^8s,$$
and $\nd(p_1, p_{k}) \geq s/2$. Further by Lemma \ref{connectivity} (using $L=2^8$) and distance estimates from Lemma \ref{estimatesrel}, the pair $(p_1, p_{k})$ is $(2^8,\delta/2, \delta (2\nD)^{-50+\log_2(\delta)-\log_2(k)})$--connected. We observe that $\nd(p_{k},p_{k+1}) \leq \delta/4 \nd(x,y)$. Thus, any curve fragment connecting $p_1$ to $p_k$ can be enlarged by a small gap to one connecting $p_1$ to $p_{k+1}$. This changes the length slightly, so by repeating an argument from Lemma \ref{changeC}, we see that the pair $(x,y)=(p_1, p_{k+1})$ is $(2^9,\delta, \delta (2\nD)^{2\log_2(\delta)-80})$--connected.
 \\

\noindent \textbf{Case 4. Each pair $(x,y) \in T \times T$ with $0<\nd(x,y) \leq 2^{-300}C^{-1}$ is $$\left(2^{40}C,2^{-14}, \epsilon (2\nD)^{-2000-2\log_2(C)}\right)\text{--connected}.$$}  Let $x \in e_x =(v_x,w_x)$ and $y \in e_y =(v_y, v_w)$. If $\nd(x,y) < \scl(v_x)2^{-4}$ or $\nd(x,y) < \scl(v_y) 2^{-4}$, the result follows from Case 1. Thus, assume that 

$$\nd(x,y) \geq \max(\scl(v_x)2^{-4}, \scl(v_y)2^{-4}).$$

Let $n_x, n_y$ be as in Case 3. In other words $\nd(n_x, x) \leq 2^8 \scl(v_x) \leq 2^{12}\nd(x,y)$, $\nd(n_y, y) \leq 2^{8} \scl(v_y) \leq 2^{12}\nd(x,y)$ and $(n_x,x)$ and $(n_y,y)$ are $(2^9,2^{-30}, (2\nD)^{-200})$--connected. Now, consider the discrete path $p_1\defeq x, p_2\defeq n_x, p_3\defeq n_y, p_4\defeq y$. Note, that $(p_2, p_3)$ is $(2^{25}C,2^{-30}, \epsilon (2\nD)^{-1000 - \log_2(C)})$--connected by Case 2. We have

$$\sum_{i=1}^3 \nd(p_i, p_{i+1}) \leq  2^{15}\nd(x,y).$$

Thus from Lemma \ref{connectivity}, we get that $(x,y)$ is $(2^{40}C,2^{-14}, \epsilon (2\nD)^{-2000-2\log_2(C)})$--connected. \\

\noindent \textbf{Case 5. For each $x \in T$ and $y \in K$ the pair $(x,y)$ with $d(x,y) \leq 2^{-100}C^{-1}$ is $$(2^{40}C,2^{-14}, \epsilon (2\nD)^{-2000-2\log_2(C)})\text{--connected}.$$} Define $n_x$ as before. Then we can consider the discrete path  given by $p_1\defeq x, p_2\defeq n_x$ and $p_3\defeq y$. The result follows analogously to Case 4 from Lemma \ref{connectivity}. \\

\noindent \textbf{Concluding remarks:} We have established $(\overline{D}, 2^{-5})$-doubling and $$\left(2^{40}C, 2^{-14},\epsilon (2\nD)^{-2000-2\log_2(C)},2^{-300}C^{-1}\right)\text{--connectivity}.$$ By Theorem \ref{maintheorem} we obtain that the space $(\nX,\nd,\nmu)$ is a PI-space with appropriate parameters. The constants in this theorem almost certainly could be substantially improved, but that is not relevant for us.
\QED

In the following appendix we give a result  for tangents of RNP-Lipschitz differentiability spaces which can be seen as a corollary of the previous result. However, the alternative proof is simpler.
\appendix
\section{Tangents of RNP-Lipschitz differentiability spaces} \label{appendix}

We give a shorter proof that tangents of RNP-Lipschitz differentiability spaces are almost everywhere PI-spaces. First a result on stability of the connectivity criterion. For the definition of measured Gromov-Hausdorff convergence see, for example, \cite{gigligromovhaus, keith2003modulus}.

\begin{theorem}\label{thm:stab} If $(X_i,d_i,\mu_i,x_i) \to (X,d,\mu,x)$ is a convergent sequence of proper pointed metric measure spaces in the measured Gromov-Hausdorff sense, and each $X_i$ is $(C,\delta,\epsilon)$--connected, then also the limit space $(X,d,\mu)$ is $(C,\delta',\epsilon)$--connected for every $\delta'>\delta$. 
\end{theorem}

\noindent\textbf{Proof:} It is sufficient to assume $\delta<\delta'<1$, as otherwise the claim is trivial. Consider the spaces $X_i,X$ embedded isometrically in some proper super-space $(Z,d)$, such that the induced measures converge weakly and $X_i$ converge to $X$ in the Hausdorff sense (see for similar argument \cite{keith2003modulus}). Let $x,y\in X$ be arbitrary with $d(x,y)=r$ and $E \subset B(x,Cr)$ be a set such that $$\mu(E \cap B(x,Cr)) < \epsilon \mu(B(x,Cr)).$$ There is a $0<\epsilon'<\epsilon$, such that
$$\mu(E \cap B(x,Cr)) < \epsilon' \mu(B(x,Cr)).$$ By regularity, no generality is lost by assuming that $E$ is open. Define a function on $X$ for each (fixed) $\eta>0$ by
$\rho_\eta(z) \defeq \min\left(\frac{d(z,E^c )}{\eta}, 1\right).$

Extend the function to $Z$ in such a way that it is compactly supported, Lipschitz and bounded by $1$, and denote the extension by the same symbol. By weak convergence

\begin{equation}
\int_Z \rho_\eta d\mu_i \to \int_Z \rho_\eta d\mu \leq \mu(E). \label{limit}
\end{equation}
Thus, the open sets $E_i \defeq \left\{ \rho_\eta > \dfrac{\epsilon'}{\epsilon} \right\} \cap X_i$ satisfy by \eqref{limit}

\begin{eqnarray}
\limsup_{n \to \infty}\mu_n(E_n) &\leq& \lim_{n \to \infty }\frac{\epsilon}{\epsilon'}\int_{X_n} \rho_\eta ~d\mu_i \nonumber \\
&\leq& \frac{\epsilon}{\epsilon'}\mu(E) < \epsilon \mu(B(x,Cr)). \label{Eiest}
\end{eqnarray}

Further, choose sequences $x_i,y_i \in X_i$ such that $x_i \to x$ and $y_i \to y$.  Our balls are assumed to be open, so from lower semi-continuity of the volume for open sets we see

$$\liminf_{n \to \infty} B(x_i, Cd(x_i, y_i)) \geq \mu(B(x,Cr)).$$
Finally, combining this with estimate \eqref{Eiest} gives a $N$ large enough such that for all $i > N$ it holds $$\mu_i(E_i) < \epsilon B(x_i, Cd(x_i, y_i)).$$ We can now choose, using the assumption and reparametrization from Lemma \ref{adjust}, for all $i>N$ some $1$-Lipschitz curve fragments $\gamma_i\co K_i \to X_i \subset Z$ defined on compact subsets $K_i \subset [0, 2Cd(x,y)]$, and satisfying 

\begin{itemize}
\item $0 = \min(K_i)$ and $\max(K_i) \leq Cd(x_i,y_i) \leq 2Cd(x,y)$, 
\item $$\gap(\gamma_i)=\left|[0,\max(K_i)] \setminus K_i\right| \leq \delta d(x_i,y_i), \text{and}$$
\item $$\rho_\eta(\gamma_i(t)) \leq \frac{\epsilon'}{\epsilon}<1$$ for every $t \in K \setminus \{0,\max(K_i)\}$. 
\end{itemize}

By an Arzela-Ascoli argument for the proper space $Z$ and Lemma \ref{adjustment} we can choose a sub-sequence $i_k$, a compact set $K$ such that $d_H(K_{i_{k}}, K) \to 0$, and a map $\gamma\co K \to X \subset Z$ which is a limit of $\gamma_{i_k}$. Here $d_H$ is the Hausdorff metric for compact set. Also, for any $t_k \in K_{i_k}$ such that $t_k \to t$ we have $t \in K$ and $\gamma_{i_k}(t_k) \to \gamma(t)$. From this it is easy to see $\gamma(0)=x$. Further, as $\lim_{k \to \infty} \max(K_{i_k}) = \max(K)$, we have $\gamma(\max(K))=y$. Finally, using the continuity of $\rho_\eta$ on $Z$ we get $\rho_\eta(\gamma(t)) \leq \frac{\epsilon'}{\epsilon}<1$ for any $t \in K \setminus \{0,\max(K)\}$. This means that $\gamma(t) \not \in \{\rho_\eta > \frac{\epsilon'}{\epsilon}\}$ for the same range of $t$.

By upper semi-continuity with respect to Hausdorff convergence of the volume for compact sets
$$\limsup_{k \to \infty}|K_{i_k}| \leq |K|,$$
and further
$$\left|\left[0,\max(K)\right] \setminus K \right| = \max(K)-|K| \leq \lim\inf_{k \to \infty}\max{K_{i_k}}-|K_{i_k}| \leq \delta d(x,y).$$
This shows $|\undf(\gamma)| \leq \delta d(x,y)$ and since the limit curve is 1-Lipschitz also $\gap(\gamma) \leq |\undf(\gamma)| \leq \delta d(x,y)$. 

Let us now allow $\eta>0$ to vary. For each $\eta$ we can do the previous construction and define $\gamma_\eta$ to satisfy the above conditions. Taking a sub-sequential limit again using Arzela-Ascoli and Lemma \ref{adjustment} with  $\eta \to 0$, we get a 1-Lipschitz curve fragment $\gamma$ with $|\undf(\gamma)| \leq \delta d(x,y)$ and $\len(\gamma)\leq C(d,y)$. Further, for any $\eta>0$ and any $t\in K_i \setminus \{\min(K_i), \max(K_i)\}$ we have  $\rho_\eta(\gamma_\eta(t)) \leq \frac{\epsilon'}{\epsilon}<1$. Thus, from the definition of $\rho_\eta$, we get

$$d(\gamma_\eta(t), E^c) \leq \eta \frac{\epsilon'}{\epsilon}.$$

Letting $\eta \to 0$ we get for the limiting curve that $d(\gamma(t), E^c) = 0$ for $t \in K \not\in\{\min(K),\max(K)\}$, and thus $\gamma(t)$, when defined and excluding end-points, is not in $E$. Here, we used that $E$ is open. Thus $\gamma$ is the desired curve essentially avoiding $E$ with $\gap(\gamma) \leq \delta d(x,y) < \delta' d(x,y)$ and $\len(\gamma) \leq Cd(x,y)$.

\QED

For the purposes of taking a tangent we need a slightly more general version of the previous theorem.

\begin{theorem}\label{subsetlimits} If $(X_i,d,_i,\mu_i,x_i) \to (X,d,\mu,x)$ is a convergent sequence of proper pointed metric measure spaces and $S_i \subset X_i$ is a sequence of subsets such that $X$ is $(C,\delta,\epsilon, r_i)$--connected along $S_i$ and $S_i$ is $\epsilon_i$-dense in $B(x_i, R_i)$ where $\lim_{i\to \infty} R_i=\infty$ and  $\lim_{i\to \infty} \epsilon_i=0$, $\lim_{i \to \infty} r_i = \infty$, then the limit space $(X,d,\mu)$ is $(C,\delta',\epsilon)$--connected for every $\delta'>\delta$.  
\end{theorem}

\noindent \textbf{Proof:} The proof proceeds exactly as before, and allows for constructing curve fragments $\gamma_i$ in $X$ that connect pairs of points in $S_i$, and taking their limits. The additional point we make, is that for any $x,y \in X$, there is a sequence of points $x_i, y_i \in S_i$ that converge to it.
\QED

The following corollary would be a consequence of Theorem \ref{thm:PIrect} but we provide an alternative and much simpler proof for it. For a metric doubling space $M$, let $T_x(M)$ denote the set of of measured Gromov-Hausdorff tangents of $M$ at $x$.

\begin{corollary}  Let $(X,d,\mu)$ be a RNP-differentiability space. Then $X$ can be covered, up to measure zero, by countably many positive measure subsets $V_i$, such that each $V_i$ is metric doubling, when equipped with its restricted distance,  and for $\mu$-a.e. $x \in V_i$ each space $M \in T_x(V_i)$ admits a $(1,p)$-Poincar\'e-inequality for some $p \in [1,\infty)$.
\end{corollary}

\noindent \textbf{Proof:} By \cite{bate2015geometry} we have a decomposition into sets $K_i^j \subset V_i$ with $X = \bigcup V_i \cup N = \bigcup K_i^j \cup N'$, such that each $V_i$ and $K_i^j \subset V_i$ satisfies the following,

\begin{itemize}
\item $X$ is uniformly $(D_i,r_i)$-doubling along $V_i$.
\item $X$ is uniformly $(C_i, 2^{-30}, \epsilon_i,r_i)$--connected along $V_i$.
\item $V_i$ is a uniform $(\frac{1}{2}, r^i_j)$-density set in $X$ along $K_i^j$.
\end{itemize}

Blowing up $V_i$ with its restricted measure and distance, at a point of $K_i^j$, gives the desired result by Lemma \ref{subsetlimits}. 

\QED

\noindent \textbf{Remark:} The previous proof indicates the problem of applying the iterative procedure from Theorem \ref{thm:contheorem} to conclude Poincar\'e inequalities on a RNP differentiability space. The constructed paths don't need to lie inside $K_i^j$, but instead only close by, and the closeness dictated by the density of the set and doubling. Since the curves may leave $K_i^j$, we cannot guarantee the ability to refill their gaps. However, since they still can be forced to lie in $V_i$, and $V_i \setminus K_i^j$ possesses a graphical approximation, the curves can be pushed into the filling $\nK$ constructed in the previous section. This filling consists of $K_i^j$ together with the graph $T$ constructed before.

\bibliographystyle{amsplain}\bibliography{geometric}

\end{document}